%% file: Main.tex
\documentclass[10pt,openany]{bookk}
\usepackage{amsmath,amsthm,amsfonts,amssymb,amscd}
\usepackage[verbose,includemp=false, paperwidth=15.50cm, paperheight=21cm, body={10.50cm,16cm}, rmargin= 1.5cm, twosideshift=0pt]{geometry}
\usepackage[cam,letter,center]{crop}
\usepackage[latin1]{inputenc}
\usepackage{cite}
\usepackage{color}
\usepackage{mathrsfs}
\usepackage{stmaryrd}	

\newtheorem{thm}{Theorem}[section]
\newtheorem{lemma}[thm]{Lemma}
\newtheorem{prop}[thm]{Proposition}
\newtheorem{cor}[thm]{Corollary}

\theoremstyle{definition}
\newtheorem{defn}[thm]{Definition}
\newtheorem{example}[thm]{Example}

\newtheorem{exer}{Exercise}[section]

\theoremstyle{remark}
\newtheorem*{remark}{Remark}
\numberwithin{equation}{section}

\newcommand{\ad}{\text{ad}}

\newcommand{\st}{\; | \;}             
\newcommand{\tr}{\text{tr}}

\newcommand{\hol}{\mathfrak{hol}}                         

\newcommand{\PP}{\mathbb{P}}      
\newcommand{\T}{\mathcal{T}}      

\newcommand{\Om}{\Omega}
\newcommand{\om}{\omega}

\newcommand{\C}{\mathbb{C}}       
\newcommand{\Z}{\mathbb{Z}}       
\newcommand{\R}{\mathbb{R}}       
\newcommand{\Ham}{\mathbb{H}}        

\newcommand{\CP}{\mbox{$\mathbb{CP}$}}        
\newcommand{\HP}{\mbox{$\mathbb{HP}$}}         

\newcommand{\ii}{\sqrt{-1}}

\newcommand{\gr}{\nabla}

\newcommand{\Ric}{\text{Ric}}
\newcommand{\Vol}{\text{Vol}}
\newcommand{\Ker}{\text{Ker}}
\newcommand{\ibar}{\bar {i}}
\newcommand{\kbar}{\bar {k}}
\newcommand{\jbar}{\bar {j}}
\newcommand{\vbar}{\bar {v}}
\newcommand{\zbar}{\bar {z}}
\newcommand{\lbar}{\bar {l}}
\newcommand{\pbar}{\bar {p}}

\newcommand{\kahler}{K\"ahler     }
\newcommand{\hk}{hyperk\"ahler     }
\newcommand{\KE}{K\"ahler-Einstein     }
\newcommand{\ddzi}{\frac{\partial}{\partial z^i}}

\providecommand{\ip}[1]{\ensuremath{\langle #1\rangle}}

\newcommand{\Hol}{\text{Hol}}
\newcommand{\beeq}{\begin{equation}}
\newcommand{\eeeq}{\end{equation}}
\newcommand{\beit}{\begin{itemize}}
\newcommand{\eeit}{\end{itemize}}
\newcommand{\bedes}{\begin{description}}
\newcommand{\eedes}{\end{description}}
\newcommand{\been}{\begin{enumerate}}
\newcommand{\eeen}{\end{enumerate}}
\newcommand{\beeqnar}{\begin{eqnarray*}}
\newcommand{\eeeqnar}{\end{eqnarray*}}

\def\mO{{\mathbb{O}}}
\def\cA{{\mathcal{A}}}

\def\cD{{\mathcal{D}}}
\def\cP{{\mathcal{P}}}

\def\cH{{\mathcal{H}}}
\def\cO{{\mathcal{O}}}
\def\ddt{{\frac{\partial}{\partial t}}}
\def\ddu2{{\frac{\partial^2}{\partial u^2}}}

\def\PP {{\mathbb P}}
\def\ZZ {{\mathbb Z}}

\def\RR {{\mathbb R}}

\def\HH{{\mathbb H}}
\def\CC {{\mathbb C}}

\def\La{\Lambda}
\def\De{\Delta}
\def\Om{\Omega}
\def\Ga{\Gamma}

\def\si{\sigma}

\def\om{\omega}

\def\ve{\varepsilon}
\def\vp{\varphi}
\def\rar{\rightarrow}
\newcommand{\del}{{\partial}}
\newcommand{\delbar}{\bar{\partial}}

\title{Holonomy groups in riemannian geometry}
\author{Andrew Clarke and Bianca Santoro}
\date{\today}

\begin{document}

\maketitle

\chapter*{Preface}

These are the very unpretentious lecture notes for the minicourse
{\em Holonomy Groups in Riemannian geometry}, a part of the XVII Brazilian
School of Geometry, to be held at UFAM (Amazonas, Brazil), in July of 2012.
We have aimed at providing a first introduction to the main general
ideas on the study of the holonomy groups of a Riemannian manifold.

The holonomy group  is one of the fundamental analytical objects that one can define on a riemannian manfold. It interacts with, and contains information about a great number of geometric properties of the manifold. Namely, the holonomy group  detects the local reducibility of the manifold, and also whether the metric is locally symmetric.

We have structured these notes to emphasize the principal bundle formulation. This is done because it was felt that the symmetries of a manifold, even pointwise, are best expressed in this way. We have striven though to give constructions in two different ways, on vector bundles and principal bundles, and also to stress the ways in which these perspectives relate.

We plan to keep improving these notes, and the updates will be
available at

{\bf http://www.sci.ccny.cuny.edu/$\sim$bsantoro/} and

{\bf http://www.ime.usp.br/$\sim$clarke/en/}

Meanwhile, we invite the
reader to send suggestions, comments and possible corrections to

{\it bsantoro@ccny.cuny.edu} or {\it  clarke@ime.usp.br}


\clearpage\mbox{}\clearpage

\tableofcontents
\setcounter{tocdepth}{2}

\input{Introduction.tex}

\input{BundlesandConnections.tex}

\input{ParallelTransportandHolonomy.tex}

\input{RiemannianHolonomy.tex}

\input{IrreducibleRiemannianGroups.tex}

\bibliographystyle{amsalpha}
\bibliography{Bibsantoro.bib}

\end{document}

%% file: Introduction.tex
\chapter{Introduction} 

It can be said that there are many ways of expressing the geometry of a space. The different ways can be thought of as delineating different areas of mathematics. A space in algebraic geometry is typically distinguished by its structure sheaf. That is, it is determined by a ring (rather, a sheaf of rings) that determines the local properties of the space, and carries the information about how these objects are pieced together globally. In hyperbolic geometry, that being of manifolds that admit complete metrics of constant negative curvature, much of the structure can be reduced to studying particular subgroups of the set of isometries of hyperbolic space. Some spaces can be studied by understanding a distinguished family of submanifolds, such as calibrated or isoparametric submanifolds. Perhaps the most important is of the geometry that comes from a fixed tensor field, such as a metric or symplectic form. 

Although the main aim of these notes is to study riemannian geometry, we will do so by looking at the geometry in a particular way. This is by the notion of parallelism. In euclidean space, the notion of two vectors with different base points being parallel is clear. Two vectors are parallel if we can carry one vector from its base point to the other one and map onto the other vector. 

On other manifolds, it is less clear how to make this definition. We consider a connection $\nabla$ on the tangent bundle $TM$. This is a differential operator on vector fields, that allows us to say when a vector field is constant along a curve. Two vectors, defined at different points, are parallel to one another with respect to a fixed curve if we can carry one vector to the other using the connection. That is, if $\gamma$ is the curve with endpoints $\gamma(0)=p$ and $\gamma(1)=q$, then $v\in T_pM$ is parallel to $w\in T_qM$ if there is a vector field $X$ along $\gamma$ that takes the values of $v$ and $w$ at the endpoints, and satisfies $\nabla_{\partial_t}X=0$.
We are able to define a isomorphism $P_\gamma:T_pM\to T_qM$ that identifies parallel vectors. 

The complication with this construction though, is the dependence of the parallelism on the curve that we take between $p$ and $q$. If $\gamma$ and $\tau$ are curves with the same endpoints, there is no particular reason that the parallel translation maps $P_\gamma$ and $P_\tau$ should coincide. To measure the dependence upon the curve, we consider the set of closed curves emanating from a fixed basepoint $p\in M$ and the set of automorphisms $P_\gamma$ of $T_pM$. This is the holonomy group of $\nabla$, that we will denote $Hol(\nabla,p)$. In the riemannian context, we will take the Levi-Civita connection.

The holonomy group of a riemannian manifold is one of the fundamental analytical objects that one can define on a riemannian manifold. It interacts with, and contains information about a great number of geometric properties of the manifold. The curvature tensor generates the infinitesimal holonomy transformations. One can define a surjective group homomorphism from the fundamental group of the manifold onto the quotient of the holonomy group by its connected component. The holonomy group also detects the local reducibility of the manifold, and also whether the metric is locally symmetric. These properties form part of the hypotheses for the fundamental theorem of Berger that classified the possible holonomy groups of a riemannian manifold.  

In a certain sense, the classification theorem is the climax of these notes. It shows that in the riemannian context, all of the study that we have made of holonomy groups can be reduced to a short list of cases. This theorem also demonstrates the way in which holonomy relates to a number of very central areas of riemannian geometry. K\"ahler, hyper-k\"ahler, quaternion-K\"ahler and Calabi-Yau manifolds have been greatly studied in recent years, and the holonomy perspective is a particularly fruitful way of looking at these objects. The list of Berger also contains two exceptional examples : $G_2$ and $Spin(7)$. These geometries are also very beautiful, and much has been recently done to better understand them.

We have structured these notes to emphasise the principal bundle formulation. This is done because it was felt that the symmetries of a manifold, even pointwise, are best expressed in this way. We have striven though to give constructions in two different ways, on vector bundles and principal bundles, and also to stress the ways in which these perspectives relate. 

In the chapter to follow, we have discussed vector bundles and principal bundles and the way that connections can be defined on them. A large number of examples have been given that attempt to relate these definitions to other, more familiar objects. In particular, we define the torsion of a connection and prove the existence of the Levi-Civita connection for any metric from the principal bundle point of view. 

In the next chapter we define parallel transport and holonomy and study the detailed properties of the holonomy group. We show that it is a Lie subgroup of the structure group, and we show how it relates to the curvature of the connection. In particular, we prove the Ambrose-Singer theorem that the holonomy Lie algebra is spanned by the values of the curvature form, as it takes values over the reduced sub-bundle. 

We then concentrate on riemannian holonomy. We see that the identity component of the holonomy group is a closed subgroup of $SO(n)$. We then make the elementary but important observation relating invariant tensors at a point with parallel tensors on the manifold. This is important in each example that arises from the classification theorem, and shows that the exceptional examples can be considered similarly to the more familiar ones. We consider holonomy groups that do not act irreducibly, and prove the de Rham decomposition theorem. We consider the set of abstract curvature tensors and see what space is required when the holonomy group is a proper subgroup. We then give a rapid introduction to the geometry of symmetric spaces. This leaves us with the classification theorem of Berger. 

In the final chapter we consider the various possibilities for a non-symmetric, irreducible holonomy group. We ouline the basiscs of K\"ahler and Calabi-Yau geometry, including some analytical aspects. For the holonomy groups modelled on the geometry of the quaternions, we see that one can consider auxilary bundles, called the twistor spaces, that carry much of the information about the metrics. We then discuss the exceptional groups, and discuss the geometric and curvature properties with these holonomy groups.

%% file: BundlesandConnections.tex
\chapter{Bundles and Connections}

\section{Connections on Vector Bundles}

Let $E\to M$ be a smooth vector bundle over the manifold $M$. We denote by $\Omega^0(E)$ the set of smooth sections of $E$. We say a section of the bundle $\Lambda^kT^*_M\otimes E$ is a \emph{bundle-valued} $k$-\emph{form}. We denote by $\Omega^k(E)$ the set of $E$-valued $k$-forms on $M$. The set of sections defined over some subset $U$ of $M$ is denoted $\Omega^0(U,E)$.

If $\sigma$ is a section of $E$, we would like to differentiate it. One problem that we have is that at different points in the manifold, $\sigma$ takes values in different vector spaces, so it is difficult to take the difference between them. To compare the values of $\sigma$ we have to \emph{connect} the points of $M$ in a  coherent way. 

Here is a very elegant way to deal with this problem.
\begin{defn}
A \emph{connection} on $E$ is a linear map
\begin{eqnarray}\label{eqn:connection}
\nabla:\Omega^0(E)\to \Omega^1(E)
\end{eqnarray}
that satisfies, for $f\in C^\infty(M)$, $\sigma\in\Omega^0(E)$,
\begin{eqnarray}\label{eqn:leibnitz}
\nabla(f\sigma)=df\otimes \sigma +f\nabla\sigma.
\end{eqnarray}
\end{defn}
That is, $\nabla$ satisfies a \emph{Leibnitz rule} for the pairing of smooth functions and smooth sections of $E$. We will now give some examples that will hopefully demonstrate that connections are well behaved in the standard constructions using vector bundles.

\begin{example}{\bf Levi-Civita connection} Let $E=T_M$ be the tangent bundle to the manifold $M$. If $M$ is endowed with a riemannian metric $g$, then there is a unique connection $\nabla$ that satisfies
\begin{itemize}
\item $\nabla g=0$. That is,
\begin{eqnarray*}
dg(Y,Z) =g(\nabla Y,Z)+g(Y,\nabla Z)
\end{eqnarray*}
\item $T_\nabla=0$, i.e., torsion-free, where
\begin{eqnarray*}
T_\nabla(Y,Z)=\nabla_YZ-\nabla_ZY-[Y,Z].
\end{eqnarray*}
\end{itemize}
\end{example}


\begin{example}({\bf Submanifolds and sub-bundles}) We consider a submanifold $M\subseteq \R^N$, or similarly, some other manifold. Then the tangent bundle is contained in the product
\begin{eqnarray*}
T_M\subseteq M\times \R^N.
\end{eqnarray*}
A vector field $Y$ on $M$ can then be considered as an $\R^N$-valued function on $M$. Then, $dY$ is a $1$-form on $M$ with values in $\R^N$. We then define
\begin{eqnarray*}
\nabla Y=\pi(dY)
\end{eqnarray*}
where $\pi:\R^N\to T_{x,M}$ is at each point the orthogonal projection onto the tangent bundle. Then, $\nabla$ is a connection on the tangent bundle to $M$.

Similarly, let $E\to M$ be a smooth vector bundle equipped with a fiber metric $g$. Suppose that $F$ is a sub-bundle of $E$. That is, at each $x\in M$, $F_x$ is an $n$-dimensional subspace of $E_x$ (for some fixed $n$). Let $\pi_x:E_x\to F_x$ denote the orthogonal projection onto $F_x$. Let $\nabla^E$ denote a connection on $E$. Then the operator
\begin{eqnarray*}
\nabla^F\sigma =\pi(\nabla^E\sigma)
\end{eqnarray*}
defines a connection on $F$. 
\end{example}

Before we can give another construction of connections, we will describe them locally. Suppose that the vector bundle $E$ can be trivialized over some subset $U\subseteq M$. That is, let $\{\varepsilon^i\}$ denote a set of sections of $E$ over $U$ that at every point $x$ form a basis for $E_x$. Then any section of $E$ over $U$ can be expressed as a linear combination
\begin{eqnarray*}
\sigma = \varepsilon^if_i\ \ \ \text{for}\ f_i\in C^\infty(U).
\end{eqnarray*}
Then, we can express
\begin{eqnarray*}
\nabla\sigma = \varepsilon^i\otimes df_i +f_i\nabla\varepsilon^i.
\end{eqnarray*}
Then if we express $\nabla\varepsilon^i =\varepsilon^j\otimes \theta_j^{\,i}$ where $\{\theta_j^{\,i}\}$ forms a matrix of $1$-forms on $U$ so
\begin{eqnarray*}
\nabla\sigma &=& \varepsilon^i\otimes df_i+\varepsilon^j\otimes \theta_j^{\,i}f_i\\
&=& \varepsilon^i\otimes (df_i+\theta_i^{\,k}f_k).
\end{eqnarray*}
That is, given a trivialization of the vector bundle, the connection $\nabla$ acts by the standard exterior derivative $d$ on the coefficient functions plus a linear transformation $\theta$ in the coefficients. The matrix-valued $1$-form determines the connection locally.

\begin{exer} Show that if we take a new local frame $\tilde\varepsilon=\varepsilon \cdot g$ for $g\in GL(n,\R)$ a constant (that is, $\tilde\varepsilon^i=\varepsilon^jg_j^{\,i}$), then the new connection matrix is given by
\begin{eqnarray*}
\tilde\theta =g^{-1}\theta g.
\end{eqnarray*}
\end{exer}
\begin{exer}
Show that if $g$ is not constant then
\begin{eqnarray}\label{eqn:gaugetrans}
\tilde\theta =g^{-1}dg+g^{-1}\theta g.
\end{eqnarray}
\end{exer}
Equation \ref{eqn:gaugetrans} is a ubiquitous expression in the study of connections. It is according to this rule that connections transform according to the action of the \emph{gauge group} (see \cite{Laws}). This is also the rule that gives the patching data for when we can recover a connection from its local expressions.
\begin{exer} Suppose that the vector bundle $E$ admits trivializations $\{\varepsilon_\alpha\}$ defined over the open sets $U_\alpha$ that together cover $M$ and that the trivializations are related by $\varepsilon_\alpha=\varepsilon_\beta  g_{\beta\alpha}$ for maps $g:U_\alpha\cap U_\beta\to GL(n\R)$. Suppose that on the sets $U_\alpha$ there exists matrix-valued $1$-forms $\theta_\alpha$ that together satisfy
\begin{eqnarray*}
\theta_\alpha=g^{-1}_{\beta\alpha}dg_{\beta\alpha}+g^{-1}_{\beta\alpha}\theta_\beta g_{\beta\alpha}.
\end{eqnarray*}
Show that the matrices $\{\theta_\alpha\}$ piece together to give a well-defined connection on $E$.
\end{exer}

\subsubsection{The pull back of a connection} Another geometric operation that can be performed on vector bundles and connections is the \emph{pull-back}. Let $f:M\to N$ be a smooth map and $\pi:E\to N$ be a smooth vector bundle. Then we can define the vector bundle $f^*E$ over $M$.
\begin{eqnarray*}
\begin{CD}
f^*E @. E\\
@VVV @VVV\\
f:M @>>> N
\end{CD}
\end{eqnarray*}
$f^*E$ is the set $\{(x,v)\in M\times E \ ;\ f(x)=\pi(v)\}$. Then the projection to the first factor defines this as a vector bundle where the fiber over the point $x$ is
\begin{eqnarray*}
(f^*E)_x=E_{f(x)}.
\end{eqnarray*}
If $\sigma\in\Omega^0(E)$ is a section of $E$ then $f^*\sigma=\sigma\circ f$ defines a section of $f^*E$.

Let $\nabla$ be a connection on $E$. We wish to define $f^*\nabla$ as a connection on $f^*E$. Not every section of $f^*E$ is of the form $f^*\sigma$ for some $\sigma\in \Omega^0(E)$ however, but these sections do locally generate $f^*E$.
That is, let $\{\varepsilon^i\}$ be a local frame for $E$ over $U\subseteq N$. Then $\{\tilde\varepsilon^i=f^*\varepsilon^i\}$ forms a local frame for $f^*E$ over $f^{-1}(U)\subseteq M$. That is, every section of $f^*E$ over $f^{-1}(U)$ can be expressed as $\sigma=\tilde\varepsilon^if_i$ for $f_i\in C^\infty(U)$.

With respect to the local frame $\{\varepsilon^i\}$ the  connection determines a matrix of $1$-forms $\theta$. Set $\tilde\theta=f^*\theta$ and define a connection on $f^*E$ over $f^{-1}(U)$ by
\begin{eqnarray*}
\nabla^f\sigma&=& \nabla^f(\tilde\varepsilon^if_i)\\
&=& \tilde\varepsilon^i\otimes(df_i+\tilde\theta_i^{\,j}f_j).
\end{eqnarray*}

\begin{exer} Show that if we take another local frame $\{\varepsilon^{\prime\,i}\}$ over $V\subseteq N$, the connection we obtain coincides with $\nabla^f$ defined here.
\end{exer}

We conclude that $\nabla^f=f^*\nabla$ is well defined over $f^*E$.

\section{Principal bundles}\label{sec:principalbundles}

We define and consider principal fiber bundles in this section. The examples of them are numerous, but our primary example and source of intuition should be that of the frame bundle for some vector bundle. We will emphasise a certain number of points, to make clear the relationship between principal bundles and the connections on them with vector bundles and the definition of connections that we have already seen. Our treatment of the material here closely follows that of Kobayashi and Nomizu \cite{KN}.

We first recall that an action (either left of right) of a group $G$ on a set $P$ is \emph{free} if the only element of the group that possesses a fixed point is the identity element $e$. That is, (for a right action), if $u\in P$ and $g\in G$ are such that $ug =u$ then $g=e$. The group $G$ is then in a bijective correspondence with the orbit of each element of $P$. A group action is \emph{effective} if $ug=u$ for every $u\in P$ only if $g=e$. This is to say that the homomorphism $G\to \text{Aut}(P)$ is injective.

\begin{defn}
A \emph{principal fiber bundle} over $M$ with structure group $G$ consists of the following data. $P$ is are smooth manifold that is acted upon by $G$ on the right. This satisfies the following conditions :
\begin{enumerate}
\item The action
\begin{eqnarray*}
(u,g)\mapsto ug=R_gu
\end{eqnarray*}
is free.
\item $M$ is equal to the quotient of $P$ by the action of $G$ and the projection of $\pi:P\to M$ is smooth.
\item The quotient map $\pi:P\to M$ is \emph{locally trivial} in the sense that for every $p\in M$, there exists a neighbourhood $U$ of $x$ and a diffeomorphism $f:\pi^{-1}(U)\to U\times G$ that is compatible with the fight $G$-action on $\pi^{-1}(U)$ and $U\times G$.
\end{enumerate}
\end{defn}
Then, above each $x\in M$, the set $\pi^{-1}(x)$ is a submanifold of $P$ diffeomorphic to the group $G$ so that fr any $u\in\pi^{-1}(x)$, the fiber $\pi^{-1}(x)$ corresponds exactly to the orbit $u\cdot G$ of $u$ under the right action of $G$.

\begin{example}{\bf Hopf fibrations}
The first example of a principal fiber bundle that we give is the Hopf fibration of $P=S^3$, which fibers over the $S^2=\CP^1$. The structure group $G$ is the set of unit vectors in $\C$. $\CP^1$ consists of all $1$-dimensional complex linear subspaces of $\C^2$. Any unit length vector $v\in S^3\subseteq\C^2$ spans a complex line $l=[v]$, and this same complex line is also spanned by $\lambda v$ for any unit length complex number $\lambda\in S^1$. This defines a surjective map $S^3\to \CP^1$ and the fiber of a point $l\in \CP^1$ is the set of of unit vectors $\lambda v\in l$. This is the orbit of the point $v$ under the action of the group $S^1$ on $S^3$ by multiplication. In homogeneous coordinates $[z_0;z_1]$ on $CP^1$, the principal fiber bundle can be trivialized over the sets $U_0=\{[z_0:z_1]\,; z_0\neq 0\}$ and $U_1=\{[z_0:z_1]\,; z_1\neq 0\}$. These sets cover $\CP^1$, confirming the final property that we need.
\end{example} 

\begin{exer}
Over the sets $U_0$ and $U_1$, explicitly define a trivialization for this fiber bundle.
\end{exer}
\begin{exer}
Explicitly give a diffeomorphism between $\CP^1$ and $S^2$, considered as the unit sphere bundle in $\R^3$.
\end{exer}

The same construction shows that the action of the group $U(1)=S^1$ on $S^{2n+1}$ defines a principal fiber bundle over the $n$-dimensional complex projective space $\CP^n$.

A slightly more complicated example uses the quaternions. The set of quaternions form a associative normed division algebra $\Ham$ that are generated by $i,j,k$. These elements have the square equal to $-1$ and anti-commute. $\Ham$ is in particular a $4$-dimensional real vector space. We can consider $\Ham^{n+1}$ a vector space over $\Ham$, but we must pay attention because $\Ham$ is non-commutative. Scalar multiplication by $\lambda\in\Ham$ is given on a vector by \emph{right} multiplication by $\bar\lambda$. The set of unit length quaternions $S^3$ forms a group isomorphic to $SU(2)$. We moreover have that
\begin{eqnarray*}
(v,\lambda)\mapsto v\lambda \in S^{4n+3}\subseteq\Ham^{n+1}\ \ \text{for}\ \ v\in S^{4n+3},\ \ \lambda\in S^3,
\end{eqnarray*}
defines a right action of $S^3$ on $S^{4n+3}$ that preserves the unit spheres of quaternion lines. The quotient is therefore the quaternion projective space $\HP^n$. The projection $S^{4n+3}\to \HP^n$ defines a principal fiber bundle by exactly the same argument as in the previous example.

\begin{exer}
Define an isomorphism between $S^3$, as the group of unit quaternions, and $SU(2)$, the group of $2\times 2$ hermitian matrices of unit determinant.
\end{exer}

\begin{example}{\bf Homogeneous spaces}
Our next example is given by coset spaces. Let $G$ be a Lie group and $H$ a closed subgroup. Then $H$ acts on $G$ on the right by $(g,h)\mapsto gh\in G$ and the quotient space is given by $G/H=\{gH\ ;\ g\in G\}$. If $H$ is a closed subgroup of $G$ the map
\begin{eqnarray*}
G\times H&\to& G\times G\\
(g,h)&\mapsto & (g, gh)
\end{eqnarray*}
is proper. It can therefore be seen (see \cite[Sec. 2.3]{DuiKo}) that $G/H$ is a smooth manifold and the projection $G\to G/H$ has the structure of a principal fiber bundle with structure group $H$. This is shown by constructing a transverse slice : that being a submanifold $S\subseteq G$ that is transverse to each orbit that it meets, and intersects each orbit at at most one point.

These observations can be placed more geometrically by considering some specific examples. For example, the group $SO(n)$ acts freely and transitively on the sphere $S^{n-1}$, and the isotropy group of the point $(1,0,\ldots,0)$ can be seen to equal $SO(n-1)$. We therefore obtain a diffeomorphism
\begin{eqnarray*}
SO(n)/SO(n-1)\to S^{n-1}.
\end{eqnarray*}
\end{example}

\begin{exer}
Consider the set of $2$-dimensional complex vector subspaces of $\C^n$, which we denote by $G(2,\C^n)$. Show that this is acted upon transitively by the unitary group $U(n)$ with the isotropy group of $\xi=\text{span}(\varepsilon^1,\varepsilon^2)$ equal to $U(2)\times U(n-2)$. We conclude that
\begin{eqnarray*}
U(n)/(U(2)\times U(n-2))\cong G(2,\C^n).
\end{eqnarray*}
This defines a principal bundle over $G(2,\C^n)$ with structure group $U(2)\times U(n-2)$.
\end{exer}

Note though that elements of the form $\lambda \text{Id}$ for $\lambda\in U(1)$ preserve every element of $G(2,\C^n)$, so the action is not effective. The map $U(n)\to \text{Diff}(G(2,\C^n)$ is not injective. This example will be investigated further in Section \ref{sec:SymmSpaces}.


\begin{example}{\bf Covering spaces}
Let $\tilde M\to M$ be the universal covering space of a smooth manifold $M$. Then (see Greenberg and Harper \cite{GH}, \S 5.8), the fundamental group $G=\pi_1(M,x_0)$ of $M$ acts on $\tilde M$ without fixed points and such that $M$ is equal to $\tilde M/G$. The group $G=\pi_1(M,x_0)$ is a discrete group, but can be considered a $0$-dimensional Lie group. The covering space condition implies that the quotient is locally trivial. $\tilde M$ therefore defines a principal $G$-bundle over $M$.
\end{example}

\begin{example}{\bf Principal frame bundles}
Our final example is the most important. Let $M$ be a smooth manifold and $\mathcal{F}_{GL}$ be the set of ordered bases for the tangent spaces to $M$. We define a map $\pi$ that sends a basis for $T_xM$ to the point $x\in M$. The general linear group $GL(n,\R)$ acts on the right on $\mathcal{F}_{GL}$. If $\varepsilon=\{ \varepsilon^1,\ldots,\varepsilon^n\}$ is a basis for $T_xM$, and $g=(g_i^{\ j})$,
\begin{eqnarray}
\varepsilon g&=&\{ \tilde \varepsilon^1,\ldots,\tilde \varepsilon^n\},\nonumber\\
\text{for}\ \ \tilde\varepsilon^j&=&\varepsilon^i g_i^{\ j}.\label{eqn:frameaction1}
\end{eqnarray}
By taking the local trivialization of the vector bundle $TM$ over $U\subseteq M$ we can trivialize $\mathcal{F}_{GL}$ to show that $\pi^{-1}(U)$ is diffeomorphic to $U\times GL(n,\R)$.

The tangent space $T_xM$ at a point in a manifold does not necessarily have a distinguished basis. $\R^n$, as a vector space, does have a distinguished basis. We can then identify $\mathcal{F}_{GL}$ with the set of isomorphisms
\begin{eqnarray*}
u:\R^n\to T_xM,\ \ \ \text{for}\ x\in M,
\end{eqnarray*}
in such a way that the action of $GL(n,\R)$ is by multiplication on $\R^n$ : for $v\in\R^n$,
\begin{eqnarray}\label{eqn:frameaction2}
(ug)(v)=u(gv)\in T_xM.
\end{eqnarray}
\end{example}

\subsubsection{Sub-bundles and $G$-structures} The principal fiber bundle of a manifold consists of all the bases for the tangent bundle. One can describe many geometric structures on a manifold by considering \emph{sub-bundles} of $\mathcal{F}_{GL}$.

For example, an $n$-dimensional manifold $M$ is by definition \emph{orientable} if $\Lambda^nTM\setminus\{0-section\}$ has 2 connected components. $M$ is \emph{oriented} if we have chosen one of the two components. Suppose that $M$ is an oriented manifold. Let $\mathcal{F}_{GL^+}$ be the subset of $\mathcal{F}_{GL}$ consisting of frames $\varepsilon$ such that $\varepsilon^1\wedge\ldots\wedge\varepsilon^n$ lies in the chosen connected component of $\Lambda^nTM\setminus\{0-section\}$. Then, the group $GL^+(n,\R)$ of $n\times n$ matrices of positive determinant preserves this set, and in fact acts transitively on the fibers of the projection $\pi:\mathcal{F}_{GL^+}\to M$. By similar reasoning to above, $\mathcal{F}_{GL^+}$ forms a principal bundle with structure group  $GL^+(n,\R)$.

This idea can be implemented similarly for subgroups of $GL(n,\R)$. We will describe this explicitly in two well known cases.

Suppose that $(M,g)$ is a riemannian manifold. Let $\mathcal{F}_O$ be the set of bases for fibers in $M$ that are \emph{orthonormal} with respect to the metric. Then, acting according to either of Equations \ref{eqn:frameaction1} or \ref{eqn:frameaction2}, the group $O(n)$ of orthogonal transformations preserves $\mathcal{F}_O$ and give it the structure of a principal bundle with structure group $O(n)$. If $M$ is simultaneously, one can reduce the structure group further to $SO(n)=O(n)
\cap GL(n,\R)$.

Let $M$ be $2n$-dimensional and let $J:TM\to TM$ be an endomorphism with $J^2=-Id$ on $TM$. That is, $J$ defines an almost-complex structure on $M$. We make a fixed identification of $\R^{2n}=\C^n$. Let $\mathcal{F}_{GL(\C)}$ consist of those linear maps $u:\C^n\to T_xM$ (for $x\in M$) such that
\begin{eqnarray*}
u(iv)=J_xu(v).
\end{eqnarray*}
Then, $GL(n,\C)$ ($\subseteq GL(2n,\R)$) acts on $\mathcal{F}_{GL(\C)}$ and gives it the structure of a principal bundle. $\mathcal{F}_{GL(\C)}$ consists of those frames $\{\varepsilon^1,\ldots,\varepsilon^{2n}\}$ such that
\begin{eqnarray*}
J\varepsilon^i&=& \varepsilon^{i+n}\\
J\varepsilon^{i+n}&=& -\varepsilon^i.
\end{eqnarray*}
In each of these examples, we have used the notion of a frame as a linear map from $\R^n$ to transfer a given geometric structure (orientation, metric or almost-complex structure) on the vector space $\R^n$ to one on the tangent bundle. The structure group is \emph{reduced} from $GL(n,\R)$ to the group that preserves the given structure on $\R^n$. It should be stated again that it is not necessary that the frame bundle $\mathcal{F}_{GL}$ does reduce to one with structure group $G$. Not every manifold is orientable, and not every even dimensional manifold admits an almost complex structure. A counterpoint to this statement is that one can always reduce the structure group from $GL(n,\R)$ to $O(n)$. This is related to the fact that $O(n)$ is a maximal compact subgroup of $GL(n,\R)$ and these groups have the same homotopy type.


\subsubsection{Principal bundle homomorphisms and reductions of structure group}
In these previous examples, the sets that we have considered have been submanifolds of the principal frame bundle $\mathcal{F}_{GL}$ of a manifold. We now make more concrete the idea of a subbundle.

A \emph{homomorphism} of principal bundles consists of two maps, both of which we will denote by $f$. Let $Q\to N$ be a principal $H$ bundle and $P\to M$ a principal $G$-bundle over different manifolds. Then a homomorphism of the first to the second consists of a smooth map $f:Q\to P$ and a homomorphism $f:H\to G$ such that
\begin{eqnarray*}
f(ug)=f(u)f(g)\ \ \ \text{for all }u\in Q,\ \ g\in H.
\end{eqnarray*}
In particular, a homomorphism induces a map $f:N\to M$. If $P$ and $Q$ are both bundles over $M$, we say that $Q$ is a \emph{subbundle} of $P$ if the map $Q\to P$ is an injective immersion, $H\to G$ is an injective group homomorphism, and the induced map $f:M\to M$ is the identity.

\begin{prop}\label{prop:critforsubbundle}
Let $\pi:P\to M$ be a principal $G$-bundle. Let $H$ be a Lie subgroup of $G$ and suppose that $Q$ is a subset of $P$ that satisfies the following properties :
\begin{enumerate}
\item The restriction $\pi:Q\to M$ is surjective,
\item $Q$ is preserved by right translation by elements of $H$,
\item for all $x\in M$, $H$ acts transitively on $\pi^{-1}(x)\cap Q$,
\item for all $x\in M$, there exists a local cross section of $P$, defined on a neighbourhood of $x$, that takes values in $Q$.
\end{enumerate}
Then $Q$ has the structure of a principal subbundle of $P$ with structure group $H$.
\end{prop}
\proof{
Let $\sigma:U\to P$ be a local section of $P$ that takes values in $Q$. Then, for every $u\in \pi^{-1}(U)\cap Q$, there exists $h\in H$ such that $u=\sigma(x)h$. We define a map $\psi:\pi^{-1}(U)\cap Q\to U\times H$
\begin{eqnarray*}
u\mapsto (x,h).
\end{eqnarray*}
This map is bijective, so we can define a differentiable structure on $\pi^{-1}\cap Q$ by specifying that this map is a diffeomorphism. An exercise shows that this is independent of the section $\sigma$ that we have taken.
With this topology, we can see that $Q$ satisfies the three conditions to be a principal $H$-bundle over $M$,  such that the inclusion of $Q$ into $P$, together with the inclusion of $H$ into $G$, is as a subbundle.
}

Let $\pi:P\to M$ be a principal fiber bundle with structure group $G$. The third condition in the definition of a principal bundle implies that $M$ can be covered with open sets $\{ U_\alpha\}$ on which there exist local sections $s_\alpha:U_\alpha\to P$ such that $\pi\circ s_\alpha=id$. On the intersections $U_\alpha\cap U_\beta$, the sections are equivalent under teh action of $G$. That is, they satisfy $s_\alpha=s_\beta\cdot g_{\beta\alpha}$ where $g_{\beta\alpha}: U_\alpha\cap U_\beta\to G$ is a smooth map.
\begin{exer}
Show that the collection of maps $\{ g_{\beta\alpha}\}$ satisfy :
\begin{itemize}
\item $g_{\alpha\beta}=g_{\beta\alpha}^{-1}$ on $U_\alpha\cap U_\beta$,
\item $g_{\alpha\beta}g_{\beta\gamma}g_{\gamma\alpha}=id $ on $U_\alpha\cap U_\beta\cap U_\gamma$.
\end{itemize}
\end{exer}
That is, $\{g_{\beta\alpha}\}$ defines a \v Cech cocycle with coefficients in the group $G$, with respect to the the covering $\{U_\alpha\}$.

We now describe the \emph{associated bundle construction}, to obtain a vector bundle from a principal bundle. Let $\pi:P\to M$ be a principal fiber bundle with structure group $G$. Let $V$ be a finite dimensional vector space and $\rho:G\to \text{Aut}(V)$ a representation of $G$ on $V$. We define an equivalence relation on $P\times V$ : $(u,v)\sim (ug, \rho(g)^{-1}v)$ for all $g\in G$. We then define
\begin{eqnarray*}
E=(P\times V)/\sim
\end{eqnarray*}
with the projection $\pi:P\to M$, $[u,v]\mapsto \pi(u)$.
\begin{prop}\label{prop:assocbundle}
$\pi:E\to M$ naturally has the structure of a vector bundle over $M$.
\end{prop}
\begin{exer}
Prove this proposition. In particular, use the local triviality of the fibration $\pi:P\to M$ to show that $E$ is locally trivial.
\end{exer}
It is using this relationsip between vector and principal bundles that we will compare the different notions of connections and holonomy on vector and principal bundles.
\begin{remark}
If $F$ is a smooth manifold and $\rho:G\to \text{Diff}(F)$ is a homomorphism then the same construction can be used to obtain a manifold that fibers over $M$ with all fibers diffeomorphic to $F$.
\end{remark}
\begin{prop}
There exists a one-to-one correspondence between sections $\sigma\in\Omega^0(E)$ of the vector bundle $\pi:E\to M$ and $V$-valued functions $f$ on the manifold $P$ that satisfy
\begin{eqnarray*}
f(ug)=\rho(g)^{-1}f(u),
\end{eqnarray*}
for $u\in P$ and $g\in G$.
\end{prop}
\proof{
Let $f$ be such a function on $P$. Then for $x\in $M, let $u\in\pi^{-1}(x)$. We define $\sigma$ by
\begin{eqnarray*}
\sigma(x)=[u,f(u)].
\end{eqnarray*}
Then the equivariance property of $f$ implies that the element of $E$ is independent of the choice of $u$. We also have that $\pi\circ\sigma(x)=x$.

Clearly, perhaps, every section of $E$ is obtained in such a way.
}

This identification is almost sufficient for our purposes. We also recall that if $\sigma\in\Omega^0(E)$ is a section of $E$ and $\nabla$ is a connection on $E$, then $\nabla\sigma\in\Omega^1(E)=\Omega^0(T^*M\otimes E)$. We will also describe these objects on the principal bundle.

Let $\alpha$ be a $k$-form on the manifold $P$, with values in the vector space $V$.
\begin{defn}\label{defn:horizequiv}
We say that $\alpha$ is \emph{horizontal} if $\alpha(X_1,\ldots,X_k)=0$ if any of the vectors $X_i$ lie in the kernel of the derivative $\pi_*:T_uP\to T_{\pi(u)}M$. That is,  if $X_i$ is tangent to an orbit of $G$ on $M$.

We say that $\alpha$ is a $\rho$-\emph{equivariant} $k$-form if $R_g^*\alpha=\rho(g)^{-1}\cdot\alpha$.
\end{defn}
That is, for $u\in P$ and $X_1,\ldots,X_k\in T_uP$,
\begin{eqnarray*}
(R^*_g\alpha)_u(X_1,\ldots ,X_k)&=& \alpha_{ug}(R_{g*}X_1,\ldots ,R_{g*}X_k)\\
&=& \rho(g)^{-1}\cdot \alpha_u(X_1,\ldots ,X_k).
\end{eqnarray*}
According to the terminology of Kobayashi and Nomizu (see \cite[pg. 75]{KN}), a $\rho$-equivariant $k$-form is \emph{pseudo-tensorial}. Such a form that is also horizontal is \emph{tensorial}. This terminology stems from the following proposition.

\begin{prop}\label{prop:equivfunctions}
There exists a one-to-one correspondence between $k$-forms on $M$ with values in the vector bundle $E$, that is, elements $\alpha\in\Omega^k(E)$, and horizontal, $\rho$-equivariant $k$-forms $\tilde\alpha$ on $P$.
\end{prop}
\proof{
The proof is very similar to the previous observation. Let $\tilde\alpha$ be a horizontal $\rho$-equivariant $k$-form on $P$. Let $X_1,\ldots, X_k\in T_xM$. For any $u\in\pi^{-1}(x)$, take $\tilde X_i\in T_uP$ such that $\pi_*(\tilde X_i)=X_i$. Then, the expression $\tilde\alpha_u(\tilde X_1,\ldots ,\tilde X_k)$ is independent of the choice of the vectors $\tilde X_i$. Furthermore, the element
\begin{eqnarray*}
[u,\tilde\alpha_u(\tilde X_i,\ldots ,\tilde X_k)]
\end{eqnarray*}
of $E$ is also seen to be independent of the choice of $u\in\pi^{-1}(x)$. Thus, such a form on $P$ determines a vector bundle map from $\Lambda^kTM$ to $E$, which is to say an $E$-valued $k$-form $\alpha$. Conversely, a moment's reflection shows that an $E$-valued $k$-form on $M$ similarly determines the correct form on $P$.}

It should perhaps be noted that the final statement, of independence upon $u\in P$, is true because the lifts $\tilde X_i$ at $u$ and $ug$ are mapped to one another by $R_{g*}$. They therefore both project to $X_i\in T_xM$.


\section{Connections on principal bundles}

After that relatively extensive description of principal bundles we turn to differentiation and how the notion of a connection can be given on a principal bundle. 

We first note that there is a family of distinguished vector fields on $P$ given by the action of $G$. The action of $G$ is free so for any $u\in P$ the map $\mu_u:g\mapsto ug$ is an embedding of $G$ into $P$, and the derivative at the identity
\begin{eqnarray*}
\mu_{u*} :\mathfrak{g}\to T_uP
\end{eqnarray*}
is injective, and maps onto the tangent space to the fiber $\pi^{-1}(x)$ through $u$. That is, for any $X\in\mathfrak{g}$, we can define
\begin{eqnarray}\label{eqn:Xestrella}
X^*=\frac{d}{dt}\left(u\cdot \exp(tX)\right)|_{t=0}
\end{eqnarray}
and every vector tangent to $\pi^{-1}(x)$ is equal to $X^*$ for some $X\in\mathfrak{g}$.
\begin{prop}
For $X\in\mathfrak{g}$,
\begin{eqnarray}\label{eqn:verticalequivar}
R_{g*}(X^*)=\left(\ad_g^{-1}X\right)^*
\end{eqnarray}
where here $X^*$ is evaluated at the point $u\in P$ and the term on the right at $ug\in P$.
\end{prop}
\proof{
If $X=\frac{d}{dt}g_t|_{t=0}\in\mathfrak{g}$, then
\begin{eqnarray*}
R_{g*}(X^*)=\frac{d}{dt}\left((ug)\cdot g^{-1}g_tg\right)_{t=0}=\left(\ad_g^{-1}X\right)^*.
\end{eqnarray*}
}
We will refer to the tangent space to the fiber $\pi^{-1}(x)$ at some $u\in\pi^{-1}(x)$ as the vertical subspace $\mathcal{V}_u\subseteq T_uP$. At each point it is canonically isomorphic to $\mathfrak{g}$.

We now give the definition of a connection on a principal fiber bundle. There are two equivalent ways of giving this definition, differing only in what aspects they emphasise. The first that we give emphasises the geometry of the manifold $P$ itself, and the second is more closely comparable to the definition on vector bundles.

\begin{defn}
Let $P$ be a princiapl fiber bundle with structure group $G$. A \emph{connection} $A$ on $P$ is a distribution of subspaces $A
\subseteq TP$ that satisfies the following two properties :
\begin{enumerate}
\item at every $u\in P$, the restriction of $\pi_*$ to $A_u\to T_{\pi(u)}M$ is an isomorphism,
\item the distribution is $G$-invariant. That is, $R_{g*}(A_u)=A_{ug}$.
\end{enumerate}
\end{defn}
Since the vertical space $\mathcal{V}_u$ is equal to the kernel of the projection to $T_xM$, this defines an invariant splitting of the tangent bundle of $P$
\begin{eqnarray*}
TP=\mathcal{V}\oplus A.
\end{eqnarray*}
We will sometimes refer to $A_u$ as the horizontal subspace at $u$. Furthermore, if $\pi(u)=x$, for any $Y\in T_xM$ there exists a unique
vector $Y^H\in A_u$ such that $\pi_*(Y^H)=Y$. $Y^H$ is called the \emph{horizontal lift of } $Y$. We have that $R_{g*}Y^H\in A_{ug}$ and $\pi_*R^{g*}(Y^H)=\pi^*Y^H=Y$, so the horizontal lift of a vector at one point in the fiber above $x$ determines the others, by right-translation.

We now recall two facts in order to give the second definition of a connection on a principal bundle. The first is that at each $u\in P$, we have an isomorphism
\begin{eqnarray}\label{eqn:defnphi}
\phi:\mathcal{V}_u\to \mathfrak{g}
\end{eqnarray}
 that inverts the map given in Equation \ref{eqn:Xestrella}. These second fact is that a connection $A$ defines a $G$-invariant splitting of TP
\begin{eqnarray*}
TP=\mathcal{V}\oplus A.
\end{eqnarray*}
\begin{prop}
Let $A$ be a connection on the principal $G$-bundle $P$ over $M$. Then the $\mathfrak{g}$-valued $1$-form $\phi$ on $P$ defined by extending $\phi$ to all of $TP$ by setting it identically zero on $A$ satisfies
\begin{eqnarray*}
R^*_g\phi=\ad_g^{-1}\cdot\phi.
\end{eqnarray*}
\end{prop}
$\phi$ is called the \emph{connection} $1$-\emph{form} of $A$.

\proof{
This can be seen by evaluating $(R^*_g\phi)(X)$ for $X$ vertical or horizontal respectively. In the first case, the assertion follows from Equation \ref{eqn:verticalequivar}. In the second case, it follows from the second condition in the definition of $A$, that the kernel of $\phi$ is preserved by right-translation.
}

We now consider how the connection behaves with respect to trivializations. Suppose that $s_U:U\to P$ is a local section of $P$, arising from  a local trivialization, and let $\phi$ be a connection $1$-form on $P$. Then $\phi_U=s^*_U\phi$ is a $\mathfrak{g}$-valued $1$-form on $U$. Suppose that for some other open set $V\subseteq M$, $s_V:V\to P$ is another local section. Then $s_V(x) =s_U(x)g_{UV}(x)$ where $g_{UV}:U\cap\to G$ are transition functions. If we differentiate this equation,
\begin{eqnarray}
s_{V*}(X)&=&s_{U*}(X)g_{UV}(x)+s_U(x)g_{UV*}(X)\nonumber\\
&=& s_{U*}(X)g_{UV}(x) s_V(x)g_{UV}^{-1}(x)dg_{UV}(X),\nonumber\\
\text{so}\ \ \left(s_V^*\phi\right)(X)&=& \phi(s_{U*}(X)g_{UV}) +\phi(s_V g_{UV}^{-1}dg_{UV}(X))\nonumber \\
\phi_V &=& g_{UV}^{-1}dg_{UV} +\ad_{g_{UV}}^{-1}\cdot \phi_V.\label{eqn:connformspatching}
\end{eqnarray}
This equation should be compared with Equation \ref{eqn:gaugetrans}. In the two cases, the expression $g^{-1}dg$ gives the canonical left-variant $\mathfrak{g}$-valued $1$-form on $G$, the Maurer-Cartan form. The expression $g_{UV}^{-1}dg_{UV}$ is the pull-back of this form to the set $U\cap V$.

We can conclude that if a principal bundle $P$ can be trivialized by a family of local sections $\{s_\alpha\}$ and if we locally have $\mathfrak{g}$-valued forms $\{\phi_\alpha\}$ that relate to one another via the transition functions according to Equation \ref{eqn:connformspatching}, then we can patch the data together to form a connection on $P$.

We now return to another point of view that we described earlier. Let $\pi:P\to M$ be a principal $G$-bundle and $\rho$ a representation of $G$ on the vector space $V$. Then,
\begin{eqnarray*}
E=\left(P\times V\right)/\sim
\end{eqnarray*}
is the vector bundle associated to $P$ and $\rho$. A section $\sigma$ of $E$ is defined by a function $f: P\to V$ that satisfies
\begin{eqnarray*}
f(ug)=\rho(g)^{-1}\cdot f(u)
\end{eqnarray*}
for all $u$ and $g$.
\begin{exer}
Show that if $X^*$ is a vertical vector (see Equation \ref{eqn:Xestrella}) then, evaluated at $u\in P$,
\begin{eqnarray*}
df(X^*)&=& -\rho_*(X)\cdot f(u)\\
&=& -\rho_*(\phi(X^*))\cdot f(u).
\end{eqnarray*}
\end{exer}
\begin{exer}
Show that the $V$-valued $1$-form on $P$
\begin{eqnarray*}
Df=df+(\rho_*\circ\phi)\cdot f
\end{eqnarray*}
is $\rho$-equivariant and horizontal, according to Definition \ref{defn:horizequiv}.
\end{exer}
We conclude that $Df=df+(\rho_*\circ\phi)f$ corresponds to an element of $\Omega^1(E)$. We denote this element $\nabla\sigma$.
\begin{remark}
In this context we can often omit explicit reference to the representation $\rho$ and express the above $1$-form as $Df=df+\phi f$.
\end{remark}

In the following example, we show that a connection can be thought of as a splitting of a short exact sequence. Let $\pi:P\to M$ be a principal fiber bundle with structure group $G$. Then, $\pi^*(TM)=\{(u,v)\in P\times TM\ ;\ v\in T_{\pi(u)}M$\} defines a vector bundle on $P$ that naturally admits a $G$-action, on the $u$-term. If we also consider the vector bundles $\mathcal{V}$ and $TP$, that also admit $G$-actions, then we obtain the short exact sequence of $G$-invariant vector bundles.
\begin{eqnarray*}
0\to \mathcal{V}\to TP \xrightarrow{\pi} \pi^*(TM)\to 0.
\end{eqnarray*}

\begin{exer}
Show that a connection on $P$ corresponds exactly to a $G$-invariant splitting of the short exact sequence. That is, a bundle map $\beta:\pi^*(TM)\to TP$ such that $\alpha\circ\beta=\text{id}$.
\end{exer}


\subsection{Homogeneous spaces and $G$-invariant connections}

In the following example, we return to the case where $G$ is a Lie group and $H$ is a closed subgroup. Then the set of left cosets $M=G/H=\{gH\ ;\ g\in G\}$ forms a smooth manifold. We will consider this configuration as $M$ being a space on which $G$ acts transitively on the left, and $H=\{g\in G\ ;\ g\cdot o=o\}$ being the isotropy of a point $o\in M$. The projection $\pi:G\to M$ gives $G$ the structure of a principal $H$-bundle over $M$, where $H$ acts by right translation on $G$. We consider the Lie algebra $\mathfrak{g}$ of $G$, and suppose that it decomposes as $\mathfrak{g}=\mathfrak{h}+\mathfrak{m}$ where $\mathfrak{h}$ is the subalgebra corresponding to the subgroup $H$, and $\mathfrak{m}$ is a transverse subspace that is preserved by the adjoint action of $H$ (this is to say that $M$ is a \emph{reductive} homogeneous space). That is,
\begin{eqnarray*}
\ad_H(\mathfrak{m})=\mathfrak{m}.
\end{eqnarray*}

The projection $\pi:G\to M$ maps the identity to $o\in M$ and so we have the map $\pi:\mathfrak{g}=\mathfrak{h}+\mathfrak{m}\to T_oM$. The kernel of this map is $\mathfrak{h}$ so $\mathfrak{m}$ is mapped isometrically onto $T_oM$. This observation can be made more concrete by constructing the tangent bundle via the associated bundle construction.
\begin{prop}\label{prop:tangbunfofhomspace}
There exists a vector bundle isomorphism
\begin{eqnarray*}
(G\times \mathfrak{m})/\sim &\to& TM\\
	\text{given by }\ [g,X]& \mapsto & g_*(\pi X).
	\end{eqnarray*}
\end{prop}	
The associated bundle $(G\times \mathfrak{m})/\sim$ is given by the adjoint representation of $H$ on $\mathfrak{m}$. \proof{
We first show that the map is well defined. For $h\in H$, $[gh,\ad_{h^{-1}}X]$ maps to
\begin{eqnarray*}
g_*h_*\pi(\ad_{h^{-1}}X) &=& g_*h_*\pi\left(\frac{d}{dt}(h^{-1}g_th)\right)\\
&=& g_*\pi\left(\frac{d}{dt}(g_th)\right)\\
&=& g_*\pi(X)
\end{eqnarray*}
as desired. It can be clearly seen to be surjective onto every fiber, and hence an isomorphism.
}

\begin{prop}\label{prop:invariantmetrics}
 There exists a one-to-one correspondence between $G$-invariant metrics on $M$ and $H$-invariant inner products on $\mathfrak{m}$.
\end{prop}
\proof{
We consider the projection $\pi:G\to M$ and its derivative at $e$ : $\pi:\mathfrak{g}\to T_oM$. Since the kernel of this map is $\mathfrak{h}$, the restriction to $\mathfrak{m}$ maps isomorphically onto $T_oM$. Let $\gamma_o$ be an inner product on $\mathfrak{m}$, and hence via $\pi$, also on $T_oM$. For $p=gH\in M$, we set $\gamma_p=g^{-1\ *}\gamma_o$. However, we can equally consider $g_1=gh$ for $h\in H$. $\gamma_p$ is well-defined and $G$-invariant if for $X,Y\in T_pM$,
\begin{eqnarray*}
\gamma_o(h^{-1}_*g^{-1}_*X,h^{-1}_*g^{-1}_*Y)&=&\gamma_o(g_o^{-1}X,g_o^{-1}Y)\\
\text{or } h^*\gamma_o=\gamma_o.
\end{eqnarray*}
and $\gamma$ is $H$-invariant. The converse statement is evident.

If we consider $\mathfrak{h}$ and $\mathfrak{m}$ as subspaces of $T_eG$ (rather than sets of left-invariant vector fields) we can
define a distribution on $G$ by setting, at $g\in G$, $A_g=L_{g*}\mathfrak{m}$. The first condition in the definition of a connection can be seen in this case because the subspaces $L_{g*}\mathfrak{h}$ and $L_{g*}\mathfrak{m}$ have trivial intersection. The former is the kernel of the projection
\begin{eqnarray*}
\pi_* :T_gG\to T_{gH}G/H
\end{eqnarray*}
so the latter must project isomorphically to $T_{gH}G/H$. The distribution is invariant by translation by elements in $H$ since $A_{gh}=L_{gh*}\mathfrak{m}=L_{g*}L_{h*}\mathfrak{m}$, and,
\begin{eqnarray*}
R_{h^{-1}*}A_{gh}&=& R_{h^{-1}*}L_{g*}L_{h*}\mathfrak{m}\\
&=& L_{g*}L_{h*}R_{h^{-1}*}\mathfrak{m}\\
&=&L_{g*}\ad_h\mathfrak{m}\\
&=&L_{g*}\mathfrak{m}=A_g.
\end{eqnarray*}
So we can conclude that $A$ defines a connection on the $H$-bundle over $G/H$.
\begin{exer}\label{exer:S5hasSU2str}
Show the group $SU(3)$ acts tranisitively on $S^5$, with isotropy group $H=SU(2)$. In terms of the Lie algebra of $SU(3)$, show that one can find a transverse subspace $\mathfrak{m}$, that is invariant under $SU(2)$.
\end{exer}
\begin{exer}\label{exer:pullbackprincipal}
Let $\pi:P\to N$ be a principal $G$-bundle and let $f:M\to N$ be a smooth map. Let $(A,\phi)$ be a connection on $P$. Show that one can define the bundle $f^*P$ as a fiber product over $M$, and one can define on $f^*P$ the \emph{pull-back} connection of $(A,\phi)$.
\end{exer}


\section{Curvature and torsion} \label{sec:CurvToronVBs}

\subsection{Curvature and torsion on vector bundles} We return to considering connections on vector bundles. Let $\nabla $ be a connection on the smooth vector bundle $E$. The connection defines a differential operator $\nabla:\Omega^0(E)\to \Omega^1(E)$. Then, just as we can extend the exterior derivative on functions to forms of higher degree $d:\Omega^k(\R)\to \Omega^{k+1}(\R)$, we can extend $\nabla$ to a differential operator
\begin{eqnarray*}\label{eqn:defn-dnabla}
d^\nabla &:&\Omega^k(E)\to \Omega^{k+1}(E)\\
(d^\nabla\alpha)_{X_0,\dots,X_k} &=& \sum_i(-1)^i\nabla_{X_i}(\alpha(X_0,\ldots, \hat X_i,\ldots, X_k))\\
 && +\sum_{i<j}(-1)^{i+j}\alpha([X_i,X_j],X_0,\ldots,\hat X_i,\ldots,\hat X_j,\ldots,X_k).
 \end{eqnarray*}
 This expression is highly unintuitive, but resembles very closely the explicit definition of the exterior derivative.
\begin{exer}

\begin{enumerate}
\item Show that if $f\in C^\infty(M)$ is a smooth function, $d^\nabla(f\sigma)=df\wedge\sigma +f d^\nabla\sigma$.
 \item Show that if $\alpha\in \Omega^k(\R)$ and $\sigma\in\Omega^0(E)$ is a section of $E$, then $d^\nabla(\alpha\otimes \sigma)= d\alpha\otimes \sigma +(-1)^k\alpha\wedge \nabla\sigma$.
 \end{enumerate}
 \end{exer}

Now we make use of this extension in two cases. Let $E=TM$ be the tangent bundle and $\nabla$ a connection on $TM$. We have the distinguished section being the identity transformation $Id\in \Omega^0(T^*M\otimes TM)=\Omega^1(TM)$.
\begin{defn}
The \emph{torsion} of the connection $\nabla$ is defined as $T^\nabla=d^\nabla(Id)\in\Omega^2(TM)$. That is,
\begin{eqnarray*}
T^\nabla_{X,Y}=\nabla_XY-\nabla_YX-[X,Y].
\end{eqnarray*}
\end{defn}
Our second observation is to recall that for $\sigma\in\Omega^0(E)$, $\nabla\sigma\in\Omega^1(E)$.
\begin{defn}
 Let $\nabla$ be a connection on the smooth vector bundle $E$. The \emph{curvature} of $\nabla$ is the transformation
 \begin{eqnarray*}
 R^\nabla : \Omega^0(E)&\to &\Omega^2(E)\\
 \sigma &\mapsto & d^\nabla(\nabla\sigma),\\
 R^\nabla_{X,Y}\sigma &=& \nabla_X\nabla_Y\sigma-\nabla_Y\nabla_X\sigma-\nabla_{[X,Y]}\sigma.
 \end{eqnarray*}
\end{defn}
\begin{prop}
The curvature $R^\nabla$ is linear over the space of functions
\begin{eqnarray*}
R^\nabla(f\sigma)=fR^\nabla\sigma.
\end{eqnarray*}
\end{prop}
\begin{cor} (see \cite{GKM}) For $x\in M$, $(R^\nabla\sigma)(x)$ is therefore only dependent on the value $\sigma(x)$. That is, $R^\nabla\in\Omega^2(\text{End}(E))$.
\end{cor}
These definitions are extremely well known, and these concepts are typically most useful in the forms stated above. However, when the bundle $E$ is endowed with additional structure, it is sometimes most convenient to consider the corresponding expressions on the principal bundle.

\subsection{Curvature on principal bundles} Let $\pi:P\to M$ be a principal $G$-bundle. Let $A$ be a connection on $P$ and let $\phi$ be the connection $1$-form.

\begin{defn} The \emph{curvature} of the connection $A$ is the $\mathfrak{g}$-valued $2$-form on $P$
\begin{eqnarray*}
\Omega&=&d\phi +[\phi,\phi].\\
\text{That is,}\ \ \ \Omega(X,Y)&=& d\phi(X,Y)+[\phi)X),\phi(Y)].
\end{eqnarray*}
\end{defn}
\begin{prop}\label{prop:curvequivariant} The curvature $\Om$ satisfies the invariance properties :
\begin{enumerate}
\item $R_g^*\Omega=\ad_g^{-1}\cdot\Omega$,
\item If $X^*$ be a vertical vector, then $\Omega(X^*,\cdot)=0$.
\end{enumerate}
\end{prop}
\proof{
 The connection form $\phi$ obeys the same equivariance rule, so we obtain,
\begin{eqnarray*}
R^*_g\Omega&=& R^*_gd\phi+[R^*_g\phi,R^*_g\phi]\\
&=& \ad_g^{-1}d\phi +[\ad_g^{-1}\cdot\phi,\ad_g^{-1}\cdot\phi]\\
&=& \ad_g^{-1}d\phi +\ad_g^{-1}[\phi,\phi]\\
&=& \ad_g^{-1}\Omega
\end{eqnarray*}
For the second assertion, Suppose that $X^*$ is a fundamental vertical vector field, for some $X\in \mathfrak{g}$. In the first case, we suppose that $Y^*$ is also vertical. Then,
\begin{eqnarray*}
d\phi(X^*,Y^*)&=& X^*\phi(Y^*)-Y^*\phi(X^*)-\phi([X^*,Y^*])\\
&=& -\phi([X^*,Y^*])\nonumber
\end{eqnarray*}
since the elements $\phi(X^*)$ and $\phi(Y^*)$ are the constant vectors $X$ and $Y$ respectively. Secondly, we note that the association $\mathfrak{g}\ni X\mapsto X^*$ is a Lie algebra homomorphism. Therefore,
\begin{eqnarray*}
d\phi(X^*,Y^*)=-\phi([X^*,Y^*]) &=  &-\phi([X,Y]^*)  \\ &=& -[X,Y]=-[\phi(X^*),\phi(Y^*)]
\end{eqnarray*}
as required. If $Y$ is horizontal then the equation
\begin{eqnarray*}
d\phi(X^*,Y) &+&[\phi(X^*),\phi(Y)] = X^*\phi(Y)   -Y\phi(X^*)   \\ &-&\phi([X^*,Y])+[\phi(X^*),\phi(Y)]=-\phi([X^*,Y])
\end{eqnarray*}
shows that it is sufficient to show that $\phi([X^*,Y])=0$, or that $[X^*,Y]$ is horizontal.
}

\begin{lemma}\label{lemma:bracketofhorizvertvfs}
Let $X^*$ be a fundamental vertical vector field on $P$ induced from $X\in\mathfrak{g}$. Let $Y$ be a vector field on $M$ and $Y^H$ the horizontal lift of $Y$. Then $[X^*,Y^H]=0$. If $Y^H$ is an arbitrary horizontal vector field, $[X^*,Y^H]$ is horizontal, but possibly non-zero.
\end{lemma}
\proof{(of Lemma) One can show that the flow of the the vector field $X^*$ on $P$ is given by
\begin{eqnarray*}
\varphi:P\times \R\to P,\ \ (u,t)\mapsto ug_t =R_{g_t}u
\end{eqnarray*}
where $g_\cdot:\R\to G$ is a $1$-parameter subgroup of $G$ tangent to $X$ at the identity.  The field $Y^H$ satisfies $Y^H_{ug}=R_{g*}Y^H_u$. Then, the Lie bracket is given by
\begin{eqnarray*}
[X^*,Y^H]_u&=&\lim_{t\to o}\frac{1}{t}\left( Y^H_u-R_{g_t^{-1}*}Y^H_{ug_t}\right)\\
&=& \lim_{t\to 0}\frac{1}{t}\left( Y^H_u -Y^H_u\right)\\
&=&0.
\end{eqnarray*}
If $Y^H$ is not right-invariant, the argument is similar.
}
This completes the proof of Proposition \ref{prop:curvequivariant}.

From the above calculations we can also note that if $X^H$ and $Y^H$ are horizontal vectors, then
\begin{eqnarray}\label{eqn:curvonhorizvectors}
\Omega(X^H,Y^H)=-\phi([X^H,Y^H]).
\end{eqnarray}
We can say more. If $X$ and $Y$ are vector fields on $M$, and $X^H$ and $Y^H$ are the horizontal lifts then since $\pi([X^H,Y^H])=[\pi(X^H),\pi(Y^H)]=[X,Y]$, we can see that
\begin{eqnarray*}
V^*= [X^H,Y^H]-[X,Y]^H
\end{eqnarray*}
is a fundamental vertical vector corresponding to the Lie algebra element $V=\Omega(X^H,Y^H)$. A result that follows from the above calculations is the following.

\begin{cor}
The curvature $\Omega$ of the connection $(A,\phi)$ vanishes identically if and only if the distribution $A=\ker\phi$ on $P$ is integrable.
\end{cor}

By Proposition \ref{prop:equivfunctions} these two conditions in Proposition \ref{prop:curvequivariant} show that $\Omega$ defines a $2$-form with values in a vector bundle on $M$. Let $\mathfrak{g}$ be the Lie algebra of the group $G$. Then, the adjoint representation is given by $X\mapsto \ad_gX=gXg^{-1}$. Then we will denote the vector bundle that we obtain by the associated bundle construction by $\mathfrak{g}_P$. This is a bundle of Lie algebras, with each fiber isomorphic to $\mathfrak{g}$.
\begin{cor}
The curvature $\Omega$ of a connection $A$ on a principal $G$-bundle $P$ defines a $2$-form on $M$ with values in $\mathfrak{g}_P$ that we denote $R^A$. That is,
\begin{eqnarray*}
\mathfrak{g}_P=\left(P\times \mathfrak{g}\right)/\sim,  &&   (u,X)\sim (ug,\ad_g^{-1}X),\\
R^A &\in & \Omega^2(\mathfrak{g}_P)=\Omega^0(\Lambda^2\otimes \mathfrak{g}_P).
\end{eqnarray*}
\end{cor}

\subsection{Curvature in the associated bundle construction}

We have seen that one of the primary ways of explicitly relating principal and vector bundles is by the associated bundle construction. If $\pi:P\to M$ is a principal $G$-bundle and $\rho:G\to \text{Aut}(V)$ is a representation of $G$ on a finite-dimensional vector space then $E=(P\times V)/\sim$ defines a vector bundle.
Let $A$ be a connection on $P$, with connection form $\phi$. If $\sigma\in\Omega^0(E)$ is a section of $E$ corresponding to the equivariant function $f:P\to V$, then the covariant derivative of $\sigma$ is given by, for $X\in TM$,
\begin{eqnarray*}
\nabla_X\sigma=(df(\tilde X)+\rho_*(\phi(\tilde X))f.
\end{eqnarray*}
where $\tilde X$ is a vector to $P$ that projects to $X$. (Strictly speaking the covariant derivative is given by $[u,(df_u +\rho_*(\phi)f_u)(\tilde X)]$, in the quotient space, but for notational simplicity we will not write all of this.
If $\tilde X$ is taken to be the horizontal lift of $X$, then $\phi(\tilde X)=0$ and $\nabla_X\sigma=\tilde Xf$. For two tangent vectors $X,Y$, the curvature is given by
\begin{eqnarray*}
R^\nabla_{XY}\sigma &=& \nabla_X\nabla_Y\sigma-\nabla_Y\nabla_X\sigma-\nabla_{[X,Y]}\sigma\\
&=& \tilde X(\tilde Yf)-\tilde Y(\tilde Xf)-\widetilde{[X,Y]}f\\
&=& \left([\tilde X,\tilde Y]-\widetilde{[X,Y]}\right)f
\end{eqnarray*}
In the decomposition $TP=\mathcal{V}+A$, the horizontal part of $[\tilde X,\tilde Y]$ is $\widetilde{[X,Y]}$ since they both project to the same thing in $M$, that being $[X,Y]$. This means that $([\tilde X,\tilde Y]-\widetilde{[X,Y]}$ is a vertical vector, equal to $Z^*$ for some $Z\in \mathfrak{g}$. $Z$ is given by $\phi([\tilde X,\tilde Y])$. Let $g_t$ be a $1$-parameter family in $G$ generated by $Z$. Then,
\begin{eqnarray*}
R^\nabla_{XY}\sigma &=& Z^*f\\
&=&  \frac{d}{dt}\left(f(ug_t)\right)|_{t=0}\\
&=& \frac{d}{dt}\left( \rho(g_t)^{-1}f(u)\right)_{t=0}\\
&=& -\rho_*(Z)f(u)\\
&=& -\rho_*(\phi([\tilde X,\tilde Y]))f
\end{eqnarray*}
That is,
\begin{eqnarray*}
R^\nabla\sigma =\rho_*(\Omega)\cdot f.
\end{eqnarray*}
\begin{exer}
Show that the representation $\rho:G\to \text{Aut}(V)$, and the induced map on Lie algebras $\rho_*:\mathfrak{g}\to \text{End}(V)$, induce a vector bundle map
\begin{eqnarray*}
\tilde \rho :\mathfrak{g}_P \to \text{End}(E),
\end{eqnarray*}
such that $\tilde \rho(R^A)=R^\nabla$.
\end{exer}


 \section{Torsion on the principal frame bundle}\label{sec:torsion}

The aim of this section is to study the torsion of a connection on the tangent bundle of a manifold, as originally defined much earlier,  from the point of view of principal bundles. In particular, we show that the orthonormal frame bundle to a riemannian manifold admits a unique connection for which the torsion vanishes identically.  This gives a principal-bundle interpretation of the existence of the Levi-Civita connection. To do these things, we use the fact that the frame bundle to $M$ carries the geometry of the manifold $M$ in a much greater way than an arbitrary principal bundle over $M$.

To start we give a preliminary result that has much in common with Proposition \ref{prop:curvequivariant}. We suppose that $\pi:P\to M$ is a principal $G$-bundle over $M$, and that $\rho:G\to \text{Aut}(V)$ is a representation of $G$ on the vector space $V$.  Let $\alpha$ be a $1$-form on $P$ with values in $V$ that satisfies
\begin{eqnarray*}
\alpha(X)&=& 0,\ \ \text{if}\ X \ \text{ is vertical} \\
R^*_g\alpha &=& \rho(g)^{-1}\cdot\alpha.
\end{eqnarray*}

Let $A$ be a connection on $P$ with connection form $\phi$. Then we define the $2$-form on $P$,
\begin{eqnarray}\label{eqn:defnofD}
D\alpha = d\alpha+\rho_*(\phi)\wedge\alpha
\end{eqnarray}
defined by
\begin{eqnarray*}
D\alpha(X,Y) = d\alpha(X,Y) +\rho_*(\phi(X))\alpha(Y)-\rho_*(\phi(Y))\alpha(X).
\end{eqnarray*}
\begin{lemma}
Let $\alpha$ be a $\rho$-equivariant, horizontal $1$-form on $P$. Then, $D\alpha$ satisfies
\begin{eqnarray*}
D\alpha(X,\cdot)&=& 0\ \ \text{if}\ X \ \text{ is vertical} \\
R^*_g(D\alpha)&=&\rho(g)^{-1}\cdot D\alpha.
\end{eqnarray*}
\end{lemma}
That is, $D\alpha$ is a horizontal, $\rho$-equivariant $2$-form on $P$ and so corresponds to an element of $\Omega^2(E)$.
\proof{
It should be easy to see that $D\alpha$ is $\rho$-equivariant, because $d\alpha$ is, and because of the equivariance properties of $\alpha$ and $\phi$. 
From the expression
\begin{eqnarray*}
D\alpha(X,Y) = X\alpha(Y) -Y\alpha(X)-\alpha([X,Y])+\rho_*(\phi(X))\alpha(Y)-\rho_*(\phi(Y))\alpha(X).
\end{eqnarray*}
we can see that $D\alpha(X,Y)=0$ if both $X$ and $Y$ are vertical. If $X=V^*$ is vertical, for some fixed $V\in \mathfrak{g}$, and if $Y$ is the horizontal lift of a vector field on $M$, then from Lemma \ref{lemma:bracketofhorizvertvfs}, $[X,Y]=0$. The above equation for $D\alpha$ reduces to 
\begin{eqnarray*}
D\alpha(X,Y)=X\alpha(Y)-\alpha([X,Y])+\rho(\phi(X))\cdot\alpha(Y).
\end{eqnarray*}
If $X=V^*=d/dt(ug_t)|_{t=0}$, so $\phi(X)=V$
\begin{eqnarray*}
X\alpha(Y) &=& \frac{d}{dt}\left( \alpha_{ug_t}(Y_{ug_t})\right)|_{t=0}\\
&=& \frac{d}{dt}\left(\rho(g_t)^{-1}\alpha_u(Y_u)\right)|_{t=0}\\
&=&-\rho_*(V)\cdot \alpha(Y)
\end{eqnarray*}
which implies that $D\alpha(X,Y)=0$. $D\alpha$ is therefore horizontal. 
}
It should also be noted that this lemma holds for forms of all degrees, and not solely for the derivatives of $1$-forms.

\begin{exer} Suppose that the $1$-form $\alpha$ defines an element $\sigma\in\Omega^1(E)$. Show that $D\alpha$ corresponds to the form $d^\nabla\sigma$, as defined in Equation \ref{eqn:defn-dnabla}.
\end{exer}

Just as the principal examples using the differential operator $d^\nabla$ were the curvature and torsion of a connection. we consider these examples on principal bundles.

\begin{exer} Let $f$ be a $V$-valued, $\rho$-equivariant function on $P$. Show that the curvature of the connection is given by
\begin{eqnarray*}
\Omega\cdot f=D(Df).
\end{eqnarray*}
\end{exer}

We now turn to the torsion of a connection.  Let $M$ be a smooth manifold of dimension $n$ and let $\mathcal{F}_{GL}$ be the set of orthonormal bases for $T_pM$, as $p$ varies over $M$. $\mathcal{F}_{GL}$ is the principal frame bundle of $M$. As noted earlier, a basis for $T_pM$ is equivalent to an isomorphism
\begin{eqnarray*}
u:\R^n\to T_pM
\end{eqnarray*}
and the right action of $GL(n)$ on $\mathcal{F}_{GL}$ coincides with the left action on $\R^n$. The projection of $\mathcal{F}_{GL}$ to $M$ is defined by $\pi:u\mapsto p$.

With these definitions, we can observe that $\mathcal{F}_{GL}$ canonically admits a $1$-form $\omega$ that takes values in $\R^n$ :
\begin{eqnarray*}
\omega(X)=u^{-1}(\pi_*X),\ \ \ \text{for}\ X\in\T_u\mathcal{F}_{GL}.
\end{eqnarray*}
If $X\in\mathcal{V}_u\subseteq T_u\mathcal{F}_{GL}$ is a vertical vector, then $\pi_*X=0$ so $\omega$ is a horizontal form. Furthermore,
\begin{eqnarray*}
(R^*_g\omega)(X) &=& (ug)^{-1}(\pi_*(R_{g\,*}X))\\
&=& g^{-1}u^{-1}(\pi_*X)\\
&=& g^{-1}\omega(X).
\end{eqnarray*}
We conclude that $\omega$ is $\rho$-equivariant where the representation $\rho$ is the standard multiplication of $GL)n)$ on $\R^n$. We consider the associated bundle, for this representation.

\begin{exer}
Show that the map defined by
\begin{eqnarray*}
(\mathcal{F}_{GL}\times\R^n)/\sim &\to&  TM\\
 \left[\{\varepsilon^i\},(a_i)\right] &\mapsto &\varepsilon^ia_i
\end{eqnarray*}
is well-defined and an isomorphism.

Show that the $1$-form $\omega$ corresponds to the identity, considered as a $TM$-valued $1$-form $Id\in\Omega^1(TM)$.
\end{exer}

\begin{defn}
Let $\phi$ be connection form on $\mathcal{F}_{GL}$. The \emph{torsion} of $\phi$ is the $2$-form
\begin{eqnarray*}
\Theta &=& D\omega\\
&=& d\omega +(\rho_*\phi)\wedge\omega.
\end{eqnarray*}
\end{defn}
Again, we will sometimes omit explicit reference to the representation $\rho$. A connection $\phi$ for which $\Theta=0$ is said to be \emph{torsion-free}.

The fundamental theorem of riemannian geometry states that if $(M,g)$ is a riemannian manifold, then there exists a unique torsion free connection that preserves the metric. We will now be able to prove this.

We note that if $\phi$ and $\phi^\prime$ are connection $1$-forms on $\mathcal{F}_{GL}$, then $\theta=\phi^\prime-\phi$ is a horizontal and equivariant $1$-form. The torsion of the two connections can be compared from the forumlas above to equal
\begin{eqnarray}\label{eqn:torsiondifference}
\Theta^\prime-\Theta &=& \phi^\prime\wedge\omega-\phi\wedge\omega\nonumber\\
&=& \theta\wedge\omega.
\end{eqnarray}
The differential forms $\omega$ and $\Theta$ can only defined on the frame bundle to a manifold, using the intrinsic geometry of the space. We use this intrinsic behaviour in the following construction. Let $\theta$ be a horizontal, orthogonal $1$-form on $\mathcal{F}_{GL}$. We can define an equivariant function $\tilde\theta$ on $\mathcal{F}_{GL}$ with values in $\mathfrak{gl}_n\otimes(\R^n)^*$ as follows. Let $v\in R^n$. Then,
\begin{eqnarray}\label{eqn:equivfunction}
\tilde\theta_u(v)=\theta_u(X),
\end{eqnarray}
where $X\in T_u\mathcal{F}_{GL}$ is such that $\omega(X)=v$. It can then be shown that $\tilde\theta$ is well-defined and $\rho$-equivariant, where $\rho$ is the induced representation of $GL(n)$ on $\mathfrak{gl}_n\otimes(\R^n)^*$. Then define
\begin{eqnarray*}
\mu:\mathfrak{gl}_n\otimes(\R^n)^* &\to & \R^n\otimes \Lambda^2(\R^n)^*\\
\theta&\mapsto &\theta\wedge \omega\\
\alpha\otimes\beta^*\otimes\gamma^* &\mapsto & \alpha\otimes (\beta^*\wedge\gamma^*).
\end{eqnarray*}
In particular, one can see that $\mu$ \emph{intertwines} the representations on the spaces $\mathfrak{gl}_n\otimes(\R^n)^*$ and $\R^n\otimes\Lambda^2(\R^n)^*$. We now consider the case that is relevant in riemannian geometry.
\begin{lemma}
The restriction of the map
\begin{eqnarray*}
\mu:\mathfrak{so}(n)\otimes(\R^n)^*\to \R^n\otimes\Lambda^2(\R^n)^*
\end{eqnarray*}
is an isomorphism.
\end{lemma}
\proof{
By a dimension count it is sufficient to show that $\mu$ is surjective.
\begin{eqnarray*}
\mu\Large((x\wedge y)\otimes z\Large)&=& \mu\large((x\otimes y-y\otimes x)\otimes z\large)\\
&=& x\otimes (y\otimes z-z\otimes y) -y\otimes (x\otimes z-z\otimes x)
\end{eqnarray*}
\begin{eqnarray*}
\text{so,}\ \ \mu\Large((x\wedge y)\otimes z -(y\wedge z)\otimes x +(z\wedge x)\otimes y\Large) &=& 2x\otimes (y\wedge z).
\end{eqnarray*}
$\mu$ maps onto a generating set for the codomain, so is surjective.
}

Now we return to the principal bundles. Let $(M,g)$ be a riemannian manifold and let $\mathcal{F}_O$ be the set of orthonormal bases for $T_pM$, with respect to the metric $g$.

\begin{cor}
Let $\psi$ be a horizontal, equivariant $2$-form on $\mathcal{F}_O$ with values in $\R^n$. Then, there exists a unique horizontal, equivariant $\mathfrak{so}(n)$-valued $1$-form $\theta$ such that
\begin{eqnarray*}
\mu(\theta)=\theta\wedge\omega=\psi.
\end{eqnarray*}
\end{cor}
\proof{
By the same construction as in Equation \ref{eqn:equivfunction} the form $\psi$ determines an equivariant function $\tilde\psi$. The above lemma shows that there is a function $\tilde\theta$ such that
\begin{eqnarray*}
\mu(\tilde \theta)=\tilde\psi.
\end{eqnarray*}
Moreover, since $\mu$ is at every point injective, and since it intertwines the group representations, we can conclude that $\tilde\theta$ is equivariant. Then, $\tilde\theta$ determines a $1$-form $\theta$ on $\mathcal{F}_O$, and by the construction of the map $\mu$, $\theta$ satisfies the required equation.
}
\begin{thm}
Let $(M,g)$ be a riemannian manifold. Let $\mathcal{F}_O$ be the orthonormal frame bundle. Then there exists a unique connection on $\mathcal{F}_O$ that is torsion-free.
\end{thm}
\proof{
Let $\phi$ be a connection $1$-form on $\mathcal{F}_O$. The torsion form $\Theta$ determines an equivariant function $\tilde\Theta$ on $\mathcal{F}_O$ with values in $\R^n\otimes\Lambda^2(\R^n)^*$. By the corollary, there exists a function $\tilde\theta$ such that
\begin{eqnarray*}
\tilde\Theta=-\tilde\theta\wedge\omega,
\end{eqnarray*}
or equivalently, a $1$-form $\theta$ such that $\Theta=-\theta\wedge\omega$. Therefore, by Equation \ref{eqn:torsiondifference} the connection $\phi^\prime=\phi+\theta$ has torsion
\begin{eqnarray*}
\Theta^\prime=\Theta+\theta\wedge\omega=0.
\end{eqnarray*}
}
The connection determined by the form $\phi^\prime$ is called the Levi-Civita connection of the metric $g$.

\begin{exer}
Consider the principal bundle $$SU(3)\to S^5=SU(3)/SU(2)$$ which has structure group $SU(2)$. Show that this can be considered a sub-bundle of the principal frame bundle of $S^5$. In Exercise \ref{exer:S5hasSU2str} we saw how this bundle admits an $SU(3)$-invariant connection. Calculate the form $\omega$ and the torsion $\Theta$ in terms of the decomposition of the Lie algebra $\mathfrak{su}(3)=\mathfrak{su}(2)+\mathfrak{m}$. In particular, show that the torsion is non-zero, and that the connection does not coincide with the Levi-Civita connection on $S^5$.
\end{exer}





%% file: ParallelTransportandHolonomy.tex
\chapter{Parallel Transport and Holonomy}

\section{Parallel transport on a vector bundle}

The word \emph{connection} is used in this context in differential geometry because such an object allows us to \emph{connect} the fibers of a bundle over different points of $M$. This is done by means of \emph{parallel transport}; a notion that we will discuss now.

Let $\pi:E\to M$ be a vector bundle over the manifold, and let $\nabla$ be a connection on $E$. Let $\gamma:[a,b]\to M$ be a smooth parametrized curve in $M$.

We recall the construction of the \emph{pull-back} of a connection by a smooth map and consider the bundle $\gamma^*E$ over $[a,b]$ equipped with the connection $\gamma^*\nabla=\nabla^\gamma$.
\begin{defn}
A section $\sigma\in\Omega^0(\gamma^*E)$ is \emph{parallel along} $\gamma$ if
\begin{eqnarray*}
\nabla^\gamma\sigma=0.
\end{eqnarray*}
\end{defn}
We suppose that $\gamma^*E$ is trivialized over $[a,b]$. The closed interval is contractible, so this is always possible. That is, there is a diffeomorphism
\begin{eqnarray*}
\varphi:[a,b]\times \R^n\to\gamma^*E
\end{eqnarray*}
such that $\varphi(t,\cdot):\{t\}\times\R^n\to E_{\gamma(t)}$ is a linear isomorphism for each $t\in [a,b]$. For the standard basis $\{e^i\}$ for $\R^n$,
\begin{eqnarray*}
\{\varepsilon^i(t)=\varphi(t,e^i)\}
\end{eqnarray*}
forms a basis for $\gamma^*E$ at $t\in[a,b]$. A section $\sigma$ of $\gamma^*E$ is then given by $\sigma=\varepsilon^if_i$ where $f_i:[a,b]\to \R$ are real valued functions. Then, as noted in the previous chapter,
\begin{eqnarray*}
\nabla\sigma &=& \varepsilon^i\otimes(df_i+\phi_i^{\ j}f_j)\\
\nabla_{\partial_t}\sigma&=& \varepsilon^i(\frac{df_i}{dt}+\phi_i^{\ j}(\partial_t)f_j)\\
&=& \varepsilon^i(\frac{df_i}{dt}+A_i^{\ j}f_j)
\end{eqnarray*}
where $A_i^{\ j}$ is a matrix of smooth functions. Then, the section $\sigma$ is parallel along $\gamma$ if
\begin{eqnarray}\label{eqn:ODEparallel}
\frac{df}{dt}+Af=0.
\end{eqnarray}
A standard theorem from the study of ordinary differential equations is the following.
\begin{thm}\label{thm:ODEexistence}
Let $A$ be a smooth matrix of functions on $[a,b]$. Then for any $v\in\R^n$, there exists a unique smooth function $f:[a,b]\to\R^n$ such that
\begin{itemize}
\item $\frac{df}{dt}+Af=0$,
\item $f(a)=v$.
\end{itemize}
\end{thm}
\begin{cor}
For any $v\in E_\gamma(a)$, there us a unique parallel section $\sigma\in\Omega^0(\gamma^*E)$ such that $\sigma(a)=v$.
\end{cor}
\proof{This follows directly from the proposition, after expressing $v$ in terms of the basis $\{\varepsilon^i\}$ at $t=a$. }
\begin{exer}
By proving that the operators
\begin{eqnarray*}
(T_\varepsilon f)(t)=-\int_a^tA(s)f(s)ds
\end{eqnarray*}
defined for $t\in [a,a+\varepsilon]$,  are contractions on the space of smooth functions
on  $ [a,a+\varepsilon]$ for $\ve$ sufficiently small,
prove Theorem \ref{thm:ODEexistence}.
\end{exer}

\begin{exer} Show that the solution $\sigma$ is independent of the choice of
the trivialization $\varphi$.
\end{exer}

\begin{defn} \emph{Parallel transport along} $\gamma$ is the map
\begin{eqnarray*}
P_\gamma :E_{\gamma(a)}&\to &E_{\gamma(b)},\\
P_\gamma(v)&=& \sigma(b)
\end{eqnarray*}
where $\sigma$ is the (unique) parallel section of $\gamma^*E$ such that $\sigma(a)=v$.
\end{defn}
\begin{prop}
$P_\gamma$ is a linear isomorphism.
\end{prop}
\proof{
That $P_\gamma$ is linear follows from the fact that Equation \ref{eqn:ODEparallel} is a linear differential equation. That it is well defined and injective is due to the fact that a solution to this equation is uniquely specified when a fixed vector is given.
}
\begin{exer}
Show that $P_\gamma$ is independent of the parametrization of the curve.
\end{exer}
In a certain sense, we can recover the connection from the notion of parallel transport. Let $\pi:E\to M$ be vector bundle, equipped with a connection $\nabla$. Let $\gamma$ be a curve through the point $y\in M$. Let $\sigma$ be a section of $E$ along the curve $\gamma$. We will to calculate the covariant derivative of $\sigma$ in the direction $\gamma_*(\partial_t)$. Let $\{\varepsilon^i\}$ be a basis for $E_y$. By the notion of parallel transport, we can extend each element $\varepsilon^i$ to be parallel along $\gamma$. That is, the sections $\varepsilon^i$ satisfy $\nabla_{\partial_t}\varepsilon^i=0$ for each $i$. Since parallel transport is a linear isomorphism, at each point $\gamma(t)$ in the curve the set $\{\varepsilon^i\}$ remains a basis for $E_{\gamma(t)}$. If $\sigma$ is a smooth section of $E$ along $\gamma$, it can be uniquely expressed as $\sigma(t)=\varepsilon^i(t)f_i(t)$, where the $f_i$ are smooth functions. The connection $\nabla$ therefore satisfies
\begin{eqnarray*}
\nabla_{\partial_t}\sigma &=& \nabla_{\partial_t}(\varepsilon^if_i)\\
&=& \varepsilon^i \frac{df_i}{dt}
\end{eqnarray*}
since $\nabla_{\partial_t}\varepsilon=0$.On the other hand, when we express $\sigma=\varepsilon^if_i$, we have a simple expression for the parallel transport of $\sigma$. Along the curve $\gamma$, let $P_t$ denote the parallel tranport map from $E_{\gamma(t)}$ to $E_{\gamma(0)}$. Then, $P_t\sigma(t)$ is then a curve in the vector space $E_{\gamma(0)}$.
\begin{prop}\label{prop:partranstoconnection}
Let $\sigma$ be a section of $E$ along $\gamma$. Then the covariant derivative of $\sigma$ at $y=\gamma(0)$ is given by
\begin{eqnarray*}
\nabla_{\partial_t}\sigma(0)=\lim_{t\to 0}\frac{1}{t}\left( P_t\sigma(t)-\sigma(0)\right).
\end{eqnarray*}
\end{prop}
\proof{
This follows from the above expression for $\nabla_{\partial_t}\sigma$. Since the sections $\varepsilon^i$ are parallel along $\gamma$, $P_t\varepsilon(t)=\varepsilon(0)$. Therefore, $P_t\sigma(t)=P_t(\varepsilon^i(t)f_i(t))=\varepsilon^i(0)f_i(t)$. We can then see that
\begin{eqnarray*}
\nabla_{\partial_t}\sigma(0) = \varepsilon^i(0)\frac{df_i}{dt}(0) &=& \lim_{t\to 0}\frac{1}{t}\left(\varepsilon^i(0)f_i(t)-\varepsilon^i(0)f_i(0)\right) \\
&=&\lim_{t\to 0}\frac{1}{t}\left( P_t\sigma(t)-\sigma(0)\right).
\end{eqnarray*}
}

Parallel transport can also be defined for piecewise smooth curves. That is, suppose that $\gamma:[a,b]\to M$ is a continuous map such that there exist $t_0=a<t_1<\cdots <t_k=b$ such that
\begin{eqnarray*}
\gamma_i=\gamma|_{[t_{i-1},t_i]}:[t_{i-1},t_i]\to M
\end{eqnarray*}
is a smooth curve. Then we define parallel transport along the concatenated curve $\gamma=\gamma_k*\gamma_{k-1}*\cdots *\gamma_1$ to be $P_\gamma=P_{\gamma_k}\circ P_{\gamma_{k-1}}\circ\cdots\circ P_{\gamma_1}:E_{\gamma(a)}\to E_{\gamma(b)}$.

We now consider a point $x\in M$, and the set of piecewise smooth curves in $M$ that start and end at $x$. If $\gamma$ and $\delta$ are curves in $M$ based at $x$, then $\gamma *\delta$, given by following $\delta$ and then $\gamma$, is also a piecewise smooth curve based at $x$. For the curve $\tilde\gamma$ defined by
$$
\tilde\gamma(t)=\gamma(b+a-t),$$
we have
$ P_{\tilde\gamma}=P_\gamma^{-1}.
$
\begin{defn}
Let $x\in M$. The holonomy group of the connection $\nabla$ is the group of transformations of $E_x$ given as parallel translations along piecewise smooth curves based at $x$. The group is denoted $Hol(\nabla,x)$.
\end{defn}
We will limit our study of holonomy groups for connections on vector bundles to the following proposition.
\begin{prop}
Let $x,y\in M$ and $\gamma:[a,b]\to M$ a path connecting them ; $\gamma(a)=x$, $\gamma(b)=y$. Then
\begin{eqnarray*}
Hol(\nabla,x)=P^{-1}_\gamma\cdot Hol(\nabla,y) \cdot P_\gamma.
\end{eqnarray*}
\end{prop}
The groups are naturally isomorphic. The isomorphism is not canonical though. It is dependent on the path $\gamma$.
\proof{Let $\delta$ be a piecewise smooth path based at $y$. Then $\delta^\prime=\tilde\gamma *\delta *\gamma$ is a path based at $x$ and
\begin{eqnarray*}
P_{\delta^\prime}&=&P_\gamma^{-1}\circ P_\delta\circ P_\gamma\in Hol(\nabla,x)\\
\text{so}\ \ \ Hol(\nabla,x) &\supseteq &  P_\gamma^{-1} Hol(\nabla,y) P_\gamma.
\end{eqnarray*}
The reverse inclusion is identical.
}


\section{Parallel transport and holonomy on principal bundles}\label{section:PTonPrincipalBdles}

We now consider the same situation on a principal bundle. Let $P\to M$ be a principal $G$-bundle, and let $(A,\phi)$ be a connection on $P$.

Let $\gamma:[a,b]\to M$ be a smooth curve in $M$. Suppose that $\gamma$ is an embedding.  
Let $C$ be the image of the map $\gamma$. $\pi:P\to M$ is a submersion so $\pi^{-1}(C)$ is a smooth submanifold of $P$, with two boundary components : $\pi^{-1}(\gamma(a))$ and $\pi^{-1}(\gamma(b))$. We define a vector field on $\pi^{-1}(C)$. At $u\in\pi^{-1}(C)$, let $X^*_u$ be the horizontal lift of $\gamma_*(\partial_t)$ at $x=\pi(u)$.

Then, the vector field $X^*$ points \emph{into} $\pi^{-1}(C)$ along $\pi^{-1}(\gamma(a))$ and \emph{out of} $\pi^{-1}(C)$ along $\pi^{-1}(\gamma(b))$. 
We consider the flow of the vector field $X^*$: 
$$
\varphi(u,t) \in \pi^{-1}(C)
$$
for $ x\in\pi^{-1}(C)$
such that
\begin{eqnarray*}
\frac{d}{dt}\varphi(u,t) &= &X^*_{\varphi(u,t)}\\
\varphi(u,0) &=& u.
\end{eqnarray*}
The first property of the flow is to note that if $\pi(u)=\gamma(a)$, then $\pi(\varphi(u,t))=\gamma(a+t)$. We note that the domain of the flow varies for different $u\in\pi^{-1}(C)$, but from this equation, we see that if $u\in\pi^{-1}(\gamma(a))$, $\varphi(u,t)$ is defined for $t\in[0,b-a]$.
\begin{defn}
The parallel translation along the curve $\gamma$ is the map
\begin{eqnarray*}
P_\gamma:\pi^{-1}(\gamma(a)) &\to & \pi^{-1}(\gamma(b))\\
P_\gamma(u)&=&\varphi(u,b-a).
\end{eqnarray*}
\end{defn}
We should note that one could make the above discussion slightly more elegant by considering the \emph{pull-back} of the bundle $P$ to $[a,b]$ by the map $\gamma$, with the induced pull-back connection, as discussed in Exercise \ref{exer:pullbackprincipal}. In particular we can use this to show that the flow $\varphi$ preserves the fibration over $C$, and restricts to be diffeomorphisms between fibers.
\begin{exer}
By considering the bundle $\gamma^*P$ over $[a,b]$, show that if $\pi(u)=\gamma(a)$, then
\begin{eqnarray*}
\pi(\varphi(u,t))=\gamma(a+t).
\end{eqnarray*}
\end{exer}
\begin{prop} For $g\in G$, $P_\gamma(ug)=P_\gamma(u)g$.
\end{prop}
\proof{The vector field $X^*$ is right-invariant, so the flow must also be.}

Suppose that $\gamma:[a,b]\to M$ is piecewise smooth. Then, as in the previous section, we define
\begin{eqnarray*}
P_\gamma=P_{\gamma_k}\circ P_{\gamma_{k-1}}\circ \cdots \circ P_{\gamma_1}
\end{eqnarray*}
where $\gamma_i=\gamma|_{[t_{i-1},t_i]}$ is smooth. We consider piecewise smooth curves based at $x\in M$. That is, suppose that $\gamma(a)=\gamma(b)=x$. Then, for $u\in \pi^{-1}(x)$, $P_\gamma(u)\in\pi^{-1}(x)$ and so
\begin{eqnarray}\label{eqn:elementofHol}
P_\gamma(u)=ug_\gamma\ \ \ \text{for some }\ g_\gamma\in G.
\end{eqnarray}
Let $\delta$ be some other curve based at $x$. We can concatenate $\gamma$ and $\delta$ to form $\gamma*\delta$ so
\begin{eqnarray*}
P_{\gamma*\delta}(u)=P_\gamma(P_\delta(u)) &=& P_{\gamma}(ug_\delta)\\
&=& P_\gamma(u)g_\delta\\
&=& u(g_\gamma g_\delta).
\end{eqnarray*}
\begin{defn}
Let $(A,\phi)$ be a connection on $P$, and let $u\in P$. The \emph{holonomy group} $Hol(A,u)$  is the subgroup of $G$ consisting of elements $g_\gamma$ for $\gamma$ a piecewise smooth loop based at $x=\pi(u)$.
\end{defn}

It is generally a very difficult problem to determine the holonomy group of a connection. Some of the theorems that we will encounter later in these notes will give possible restrictions and properties of holonomy groups. To begin with, we can consider a trivial example.

\begin{exer}
Take $G=S^1$ and consider $P=S^1\times S^1\to M=S^1$, where projection is on the first factor. If the coordinates on $P$ are $(\theta,\tau)$, show that
\begin{eqnarray*}
\phi:x\partial_\theta+y\partial_\tau\mapsto y-\alpha x
\end{eqnarray*}
defines a connection on $P$, where $\alpha\in\R$ is constant. Calculate the horizontal lift of $\partial_\theta$, and the flow of the vector field $X^*$. Show that the holonomy group is discrete if $2\pi\alpha$ is rational, and a dense subgroup of $S^1$ if $2\pi\alpha$ is irrational.
\end{exer}

Just as in the case for connections on vector bundles, we wish to know how $Hol(A,u)$ is dependent on the choice of $u\in P$. It is slightly more delicate in this case, although we are helped by the fact that $Hol(A,u)$ is always a subgroup of $G$, rather than of $Aut(T_xM)$, for varying $x\in M$. We first consider points in the same fiber of $P$.
\begin{prop}
Let $u\in P$ and $g\in G$. Then
\begin{eqnarray*}
Hol(A,ug)=g^{-1}\cdot Hol(A,u)\cdot g.
\end{eqnarray*}
\end{prop}
\proof{
This follows from the invariance of the parallel translation map :
\begin{eqnarray*}
P_\gamma(ug) &=& P_\gamma(u)\cdot g\\
&=& ug\cdot (g^{-1}g_\gamma g).
\end{eqnarray*}
where $g_\gamma$ is given in Equation \ref{eqn:elementofHol}.
}

\begin{prop}\label{prop:HolomBasePoint}
Let $x,y\in M$ and let $\gamma$ be a piecewise smooth path from $x$ to $y$. Let $u\in \pi^{-1}(x)$ and let $v\in P$ be that point obtained by parallel translation along $\gamma$. Then,
\begin{eqnarray*}
Hol(A,u)=Hol(A,v).
\end{eqnarray*}
\end{prop}

\proof{
Let $\tau$ be a curve based at $y$ and $g_\tau\in Hol(A,v)$ such that $P_\tau(v)=vg_\tau$. We consider the curve $\gamma^{-1}\cdot\tau \cdot\gamma$ based at $x$. Then,
\begin{eqnarray*}
P_{(\gamma^{-1}\cdot\tau \cdot\gamma)}(u) =  P_{\gamma^{-1}} P_\tau(P_\gamma u) &=&  P_{\gamma^{-1}}( P_\tau(v))\\
&=& P_{\gamma^{-1}}(vg_\tau)\\
&=& P_{\gamma^{-1}}(v)g_\tau\\
&=& ug_\tau,
\end{eqnarray*}
so $Hol(A,v)\subseteq Hol(A,u)$. Equality here is immediate.
}

We complete this section by considering a non-trivial example constructed by Hano and Ozuki \cite{HO} that demonstrates some of the properties that we will see in other results and contexts. In particular the group is a Lie subgroup of $GL(6,\R)$, but is not closed. We will later see that the connected component of the identity of the holonomy group is always a Lie subgroup, and in the riemannian context will always be closed. For that it is necessary to use a fundamental result,  the Ambrose-Singer Theorem, that will be proven later.

\begin{example}
 We consider a connection on the tangent bundle to $\R^6$. $\R^6$ has a standard frame of basis vectors $\{\varepsilon^i=\partial_{x^i}\}$, which we can consider a section of the principal frame bundle $\mathcal{F}_{GL}$ over $\R^6$. If $\phi$ is a connection form on $\mathcal{F}_{GL}$, we can pull it back and consider it as a $1$-form on $\R^6$ with values in $\mathfrak{gl}_6$. That is, $\phi(\varepsilon^k)=(\phi^{kj}_i)\in\mathfrak{gl}_6$. The torsion of $\phi$ is determined by the values
\begin{eqnarray}\label{eqn:torsionincoords}
T(\varepsilon^k,\varepsilon^l)=\varepsilon^i(\phi^{kl}_i-\phi^{lk}_i).
\end{eqnarray}
The example that we take is given by $\phi(\varepsilon^1)=0$,
\begin{eqnarray*}
\phi(\varepsilon^2) &=&
\left(\begin{array}{cccccc}
0 & 0 &      &     &              &      \\
0 & 0 &      &     &              &      \\
  &   & 0    & x^1 &              &       \\
  &   & -x^1 & 0   &              &        \\
  &   &      &     & 0            & \sqrt{2}x^1 \\
  &   &      &     & -\sqrt{2}x^1 & 0
  \end{array}\right),
\end{eqnarray*}
and for $k\geq 3$, $\phi(\varepsilon^k)=(\phi^{kj}_i)$ where $\phi^{kj}_i=0$ if $j\neq 2$, and $\phi^{k2}_i=\phi^{2k}_i$. Then, from Equation \ref{eqn:torsionincoords}, this connection is torsion-free, and the curvature tensor $\Omega$ satisfies
\begin{eqnarray*}
\Omega(\varepsilon^1,\varepsilon^2) &=&
\left(\begin{array}{cccccc}
0 & 0 &      &     &           &      \\
0 & 0 &      &     &           &      \\
  &   & 0    & 1 &             &       \\
  &   & -1 & 0   &             &        \\
  &   &      &     & 0         & \sqrt{2} \\
  &   &      &     & -\sqrt{2} & 0
  \end{array}\right),
\end{eqnarray*}
while $\Omega(\varepsilon^1,\varepsilon^k)$, $\Omega(\varepsilon^2,\varepsilon^k)$ and $\Omega(\varepsilon^k,\varepsilon^l)$ are zero, except for the positions $(2,3),\ldots,(2,6)$. Then, by a result that we have not yet seen, the theorem of Ambrose and Singer (Thm. \ref{thm:AmbroseSinger}) that very precisely relates the curvature tensor with the holonomy group, the holonomy group consists of elements of the form
\begin{eqnarray*}
\left(\begin{array}{cccccc}
1 & 0 &                &     &          &      \\
0 & 1 &  a_1           & a_2            &    a_3                 & a_4     \\
  &   & \cos(\theta)   & \sin(\theta)   &                        &       \\
  &   & -\sin(\theta)  & \cos(\theta)   &                        &        \\
  &   &                &                & \cos(\sqrt{2}\theta)   & \sin(\sqrt{2}\theta) \\
  &   &                &                & -\sin(\sqrt{2}\theta)  & \cos(\sqrt{2}\theta)
  \end{array}\right),
\end{eqnarray*}
In particular, the elements with $a_i=0$ are defined by an \emph{irrational rotation} subgroup in $S^1\times S^1\subseteq SO(4)$.
\end{example}

The group that we have just seen is somewhat pathological. It is a Lie subgroup of $GL(6,\R)$, but it is not a closed subgroup. We also note that it is not contained in any compact subgroup of $GL(6,\R)$. These properties can be well understood for holonomy groups. In particular, by Theorem \ref{thm:restrictedLiesubgroup} the connected component is a Lie subgroup. If
  $(A,\phi)$ is the Levi-Civita connection of a riemannian metric, the identity component of $Hol(A,\phi)$ is a \emph{closed} subgroup of $SO(n)$.


\section{Structure of the restricted holonomy group}

In this section and the next we give more detailed information about the holonomy group of a connection on a principal bundle. Our aim in this section is to show that $Hol(A,u)$ is a Lie subgroup of $G$. To do this is it first necessary to pay specific attention to the subgroup consisting of elements arising from parallel translation along paths homotopic to the identity. We will show that this is a path connected Lie group, and is equal to the connected component of $Hol(A,u)$ that contains the identity. We make at this point a comment on regularity. Throughout this section, for reasons of brevity we will refer to curves that are piecewise smooth.  By this we really mean \emph{piecewise of class} $C^1$. As noted in \cite[e.g. Thm. 7.2]{KN}, this distinction will not be meaningful.

Throughout this section, we assume that $\pi:P\to M$ is a principal bundle with structure group $G$ and that $(A,\phi)$ is a connection on $P$.

\begin{defn}
The \emph{restricted holonomy group} $Hol_0(A,u)$ is the subgroup of $Hol(A,u)$ consisting of group elements $g_\tau$ arising from parallel translation along $\tau$ where $\tau$ is a piecewise smooth loop in $M$, based at $x=\pi(u)$, that is homotopic to the constant curve at $x$.
\end{defn}

The first properties of the restricted holonomy group that we determine are the following.

\begin{prop} \label{prop:mapfromfundgrp} The restricted holonomy group $Hol_0(A,u)$ is a normal subgroup of $Hol(A,u)$. \\
Let $\pi_1(X,x)$ be the fundamental group of $M$ for $x\in M$. Let $u\in \pi^{-1}(x)$. There exists a surjective group homomorphism
\begin{eqnarray*}
\pi_1(M,x)\to Hol(A,u)/ Hol_0(A,u).
\end{eqnarray*}
\end{prop}
\proof{
For the first statement, observe that if $\gamma$ is any path in $M$ based at $x$, and $\sigma$ is a curve based at $x$ that is homotopic to the constant map to $x$, then $\gamma^{-1}\cdot \sigma\cdot\gamma$ is also contractible to $x$.

The homomorphism is defined as follows. For $[\gamma]\in \pi_1(M,x)$ we define
\begin{eqnarray*}
[\gamma]\mapsto g_\gamma \ \ \ \ \text{modulo}\ Hol_0(A,u).
\end{eqnarray*}
This is well-defined because if $\gamma,\sigma:[0,1]\to M$ represent the same class in $\pi_1(M,x)$, then $\gamma^{-1}\cdot \sigma$ is homotopic to zero and
\begin{eqnarray*}
g_{\gamma^{-1}\cdot \sigma}=g_{\gamma}^{-1}\cdot g_\sigma\in Hol_0(A,u)
\end{eqnarray*}
and so $g_\gamma=g_\sigma$ modulo $Hol_0(A,u)$. The homomorphism is clearly surjective.
}

We wish to prove that $Hol_0(A,u)$ is a Lie subgroup of $G$. To show this we first need a topological criterion for a subgroup of a Lie group to be of Lie type. It will then remain to show that $Hol_0(A,u)$ satisfies the necessary conditions.

In large part, the arguments in this section come from Kobayashi and Nomizu, although we will give specific references later.

\begin{thm}\label{thm:critforLiesubgroup}
Let $G$ be a Lie group and $H$ a subgroup of $G$ such that every element of $H$ can be connected by a piecewise smooth path in $H$ to the identity $e$. Then $H$ is a Lie subgroup of $G$.
\end{thm}
\proof{
We define the subspace $S\subseteq \mathfrak{g}$ to be those vectors $X$ that are tangent to smooth paths in $H$, ie for which $X=dx/dt(0)$ where $x_t$ is a smooth curve contained in $H$. Then, for $X,Y\in S$, $\lambda\in \R$,
\begin{enumerate}
\item $\lambda X\in S$ since $\lambda X=\ddt(x_{\lambda t})(0)$,
\item $X+Y\in S$ since $X+Y =\ddt(x_ty_t)(0)$.
\end{enumerate}
Moreover, $S$ is a Lie subalgebra of $\mathfrak{g}$ because $[X,Y]=dw_t/dt$ at $t=0$ where $w_{t^2}=x_ty_tx_t^{-1}y_t^{-1}$ (see \cite{Chev}). There exists a Lie subgroup $K$ of $G$ that has $S$ as its Lie algebra. $K$ is found by considering the left-invariant distribution on $G$ given at $g\in G$ by $L_{g*}S\subseteq T_gG$. This distribution in integrable because $S$ is a Lie algebra, and the maximal integral submanifold that contains $e$ is a Lie subgroup of $G$. We show that $H=K$.

To show $K\supseteq H$, let $a\in H$ and $g_t$ a curve in $H$ such that $g_0=e$ and $g_1=a$. We translate the curve slightly. That is, for $t_0\in [0,1]$ consider the curve $t\mapsto x_{t_0}^{-1}x_{t+t_0}$ contained in $H$. This curve is tangent to $L_{x_{t_0}}^{-1}(x_{t_0}^\prime)$ at the identity (for $t=0$) so this vector is contained in the subspace $S$ and $x_{t}^\prime\in L_{x_{t}}(S)$ for all $t\in [0,1]$. We therefore have that $x_t$ is contained in the maximal integral submanifold to the distribution, and hence $H\subseteq K$.

Next, we show that $K\subseteq H$. Let $\{X^i\}$ be a basis for $S$ and $g^i_t$ be curves in $H$ such that $g^i_0=e$ and $g_0^{i\prime}=X^i$. We then consider the map
\begin{eqnarray*}
f:(t_1,\ldots,t_k)\mapsto g_{t_1}^1g_{t_2}^2\cdots g_{t_k}^k
\end{eqnarray*}
maps from a neighbourhood of $0$ in $\R^k$ (for $k=\text{dim}S$) into $H$ ($\subseteq K$) but we can in particular see that the derivative at $0$
\begin{eqnarray*}
\R^k\to S
\end{eqnarray*}
is non-singular, and so $f$ defines a local diffeomorphism onto a neighbourhood of the identity in $K$. $f$ specifically takes values in $H$ however, so that neighbourhood in $K$ is contained in $H$. Since $K$ is connected, this implies that $K\subseteq H$.
}

\begin{exer}
Let $S$ be a left-invariant distribution on the Lie group $G$. Show that the subspace at the identity $S_e\subseteq \mathfrak{g}$ is a Lie subalgebra if and only if the maximal integral submanifold for $S$ that contains the identity is a Lie subgroup of $G$.
\end{exer}

The restricted holonomy group $Hol_0(A,u)\subseteq G$ consists of those elements coming via parallel translation of $u$ along piecewise smooth curves in $M$ that are bases at $x=\pi(u)$ and are contractible in $M$ to $x$. We require that  the curve is contained in a small coordinate neighbourhood, and of a very particular form. This statement is taken from Kobayashi and Nomizu \cite[pg. 285]{KN} and the result was originally due to Lichnerowicz \cite{Lich}
\begin{lemma}
Let $\mathcal{U}$ be a covering of $M$ by open sets. Let $\gamma$ be a closed, piecewise smooth curve in $M$ that is homotopic to zero. Then, $\gamma$ can be written as the concatenation of curves
\begin{eqnarray*}
\gamma =\gamma_1\cdot \gamma_2\cdots\gamma_n
\end{eqnarray*}
where each piecewise smooth curve $\gamma_i$ is of the form $\gamma_i=\mu_i^{-1}\cdot\sigma_i\cdot\mu_i$ where $\mu_i$ is a curve from $x$ to another point $y$, and $\sigma_i$ is a loop based at $y$ and contained in some element of the covering. In particular, $\sigma_i$ is given as a small rectangle, as a map from the boundary of $I\times I$ in $M$.
\end{lemma}
The proof is done by breaking the square $I\times I$ into $m^2$ small equal sqares, on each of which the homotopy of $\gamma$ to zero is smooth. The $\mu$'s and $\sigma$'s can then be taken. This is done in detail in \cite{KN}.

We take a covering $\mathcal{U}$ of $M$ by open sets, each of which is diffeomorphic to $B(0,1)\subseteq \R^n$. Consider $g_\gamma\in Hol_0(A,u)$ for $\gamma$ a curve in $M$. Then, according to the factorization lemma,
\begin{eqnarray}\label{eqn:holoneltdecomp}
g_\gamma=g_{\gamma_1}g_{\gamma_2}\cdots g_{\gamma_n}
\end{eqnarray}
where $\gamma_i=\mu^{-1}_i\cdot\sigma_i\cdot \mu_i$ is a lasso. $\mu_i$ is a piecewise smooth path from $x$ to $y$ and $\sigma_i$ is a smooth curve based at $y\in M$. From the proof of Proposition \ref{prop:HolomBasePoint}, the element of $Hol_0(A,v)$ corresponding to $\sigma_i$ is equal to the element $g_{\gamma_i}\in Hol_0(A,u)$ for the curve $\gamma_i$.
It suffices to assume that $\gamma$ is a map from the boundary of the rectangle $I\times I$ into $M$ and show that $g_\gamma$ can be connected to the identity.

We suppose that we have a smooth map $F:I\times I\to M$ and $\gamma$ is given by
\begin{eqnarray*}
\gamma(t) =\gamma_4\cdot \gamma_3\cdot \gamma_2 \cdot \gamma_1(t)
= \left\{
  \begin{array}{l l}
    F(4t,0) & \quad t\in [0,1/4]\\
    F(1,4t-1) & \quad t\in [1/4,1/2]\\
F(3-4t,1) &\quad t\in [1/2,3/4]\\
F(0,4-4t) &\quad t\in [3/4,1]
      \end{array} \right.
 \end{eqnarray*}
Then, we define a retraction of this rectangle onto the first line segment, the $s=0$ set,
\begin{eqnarray*}
G(t,s) =\left\{
\begin{array}{ll}
  F(4t,0) & \quad t\in [0,1/4]\\
    F(1,(1-s)(4t-1)) & \quad t\in [1/4,1/2]\\
F(3-4t,1-s) &\quad t\in [1/2,3/4]\\
F(0,(1-s)(4-4t)) &\quad t\in [3/4,1]
\end{array}\right.
\end{eqnarray*}
That is, for each $s$, $\gamma^s:t\mapsto G(t,s)$ is the concatenation of $4$ smooth curves. Every curve starts and ends at the point $x$. The following lemma shows that the parallel translation map defines a smooth curve in $P$, for each interval $[\frac{i-1}{4},\frac{i}{4}]$ on which $G$ is smooth.

\begin{lemma}\label{lemma:smoothcurveinP}
Let $G:I\times I\to M$ be a smooth map and let $\tau^s$ be the curve $t\mapsto G(t,s)$. Let $p(s)$ be a smooth curve in $P$ such that $\pi(p(s))=G(0,s)$. Then, define $q(s)=P_{\tau^s}$ as the parallel transport of $p(s)$ along the curve $\tau^s$. Then, $q(s)$ is a smooth path in $P$.
\end{lemma}
\proof{
This is a simple extension of the definition of parallel transport. We pull $P$ back to obtain the bundle $G^*P$ over $I\times I$ and perform parallel transport along the $\partial_t$ directions. This is a diffeomorphism over $\{0\}\times I$ to $\{1\}\times I$ and maps $p(s)$ to $q(s)$.
}

That is, the parallel translations of $u\in\pi^{-1}(x)$ define $q_1(s)$, and this determines $q_2(s)$, and so on. At the end, since all curves $\gamma^s$ start and end at $x$, the smooth curve that we end up with is contained in the fiber $\pi^{-1}(x)$. That is $q(s)=P_{\gamma^s}(u)=ug_s$. The local triviality of the bundle implies that the curve $g_s\in Hol_0(A,u)$ is smooth. Since we can see that $\gamma^1$ is the concatenation $\tilde\gamma_1\cdot\gamma_1$, the curve $g_s$ connects $g_\gamma$ to the identity, as we desired.

\begin{thm}\label{thm:restrictedLiesubgroup}
Let $A$ be a connection on the principal $G$-bundle $P$. Then the restricted holonomy group $Hol_0(A,u)$ is a Lie subgroup of $G$.
\end{thm}
\proof{
From Equation \ref{eqn:holoneltdecomp}, which follows from the factorization lemma, and from Lemma \ref{lemma:smoothcurveinP}  we see that every element of $Hol_0(A,u)$ can be connected to the identity by a piecewise smooth curve in $G$. We can then apply Theorem \ref{thm:critforLiesubgroup} to conclude that $Hol_0(A,u)$ is of Lie type.
}

\begin{cor}
The full holonomy group $Hol(A,u)$ is a Lie subgroup of $G$ and $Hol_0(A,u)$ is the connected component that contains the identity.
\end{cor}
\proof{
The restricted holonomy group $Hol_0(A,u)$ is a given as the maximal integral submanifold of the distribution $S$ on $G$. It is seen to be path connected. In the topology of the leaves of the foliation (weaker than the subset topology from $G$) this implies that $Hol_0(A,u)$ is connected.

By left translation, we can introduce a Hausdorff topology on $Hol(A,u)$ so that every point is locally homeomorphic to $Hol_0(A,u)$. $g_\tau$ and $g_\gamma$ are in the same connected component if $g_\tau^{-1}\cdot g_\gamma$ is contained in $Hol_0(A,u)$. That is, if $g_\tau=g_\gamma$ modulo $Hol_0(A,u)$. It can be seen (see \cite{somewhere}) that the fundamental group $\pi_1(M,x)$ is countable, so from Proposition \ref{prop:mapfromfundgrp} $Hol(A,u)$ has at most countably many connected components. We can therefore conclude that the topology is second-countable and $Hol(A,u)$ is a Lie group.
}



\section{Curvature and holonomy}

In this section we demonstrate the previously mentioned relationship between the holonomy group and the curvature of the connection. A priori, the connection  $(A,\phi)$ is defined on a principal fiber bundle $P$ with structure group $G$ while the holonomy group is a subgroup $Hol(A,u)$ of $G$. We can perform a \emph{reduction} of the structure group so as to be able to assume that $Hol(A,u)=G$.

In the previous section we saw that $Hol(A,u)$ was a Lie subgroup of $G$. As such, the subspace tangent to this subgroup is a Lie algebra that we denote $\hol$, leaving the connection and point of definition implicit. The connection and its curvature are $\mathfrak{g}$-valued forms on $P$.

\subsection{Values taken by the curvature form}

Our first aim is to show that the curvature, at the point $u\in P$, takes values in the subalgebra $\hol$.
The curvature of $\phi$ is defined as
\begin{eqnarray*}
\Omega=d\phi +[\phi,\phi]
\end{eqnarray*}
and so for $X^H$ and $Y^H$ horizontal vector fields on $P$, $\Omega(X^H,Y^H)=-\phi([X^H,Y^H])$. In particular, if $X^H$ and $Y^H$ are the horizontal lifts of vector fields $X$ and $Y$ on $M$, the Lie algebra element $V=\Omega(X^H,Y^H)$ corresponds exactly with the vertical vector field on $P$,
\begin{eqnarray}\label{eqn:curvdiffofbrackets}
V^*= [X^H,Y^H]-[X,Y]^H.
\end{eqnarray}
We take vectors fields $X$ and $Y$ on $M$ such that $[X,Y]=0$. We denote by $\varphi^X_t$ and $\varphi^Y_t$ the maps given by the flow of the vector fields. We consider the rectangle in $M$ given by moving $\varepsilon$ along $X$, $\varepsilon$ along the $Y$-flow, back $-\varepsilon$ along $X$ and then $-\varepsilon$ along $Y$. This is a little redundent though, because $[X,Y]=0$ implies that $\varphi^Y_{-\varepsilon}\varphi^X_{-\varepsilon}\varphi^Y_\varepsilon\varphi^X_\varepsilon(x)=x$ for all $u\in P$.

Instead, we lift to the principal bundle and define $X^H$ to be the horizontal lift of $X$ and $Y^H$ to be the horizontal lift of $Y$, and denote the flow of these vector fields by $\varphi^{X^H}_t$ and $\varphi^{Y^H}_t$ respectively. Note that since $[X,Y]=0$, $[X^H,Y^H]$ is vertical. These flow maps cover the flows of the vector fields $X$ and $Y$ and in particular
\begin{eqnarray*}
\varphi^{Y^H}_{-\varepsilon}\varphi^{X^H}_{-\varepsilon} \varphi^{Y^H}_\varepsilon\varphi^{X^H}_\varepsilon : \pi^{-1}(x)\to \pi^{-1}(x)
\end{eqnarray*}
defines the parallel transport map around the small rectangle of width $\varepsilon$ in $M$. This map can be closely related to the Lie bracket of vector fields on $P$ by the following lemma.
\begin{lemma}
Let $f$ be a smooth test function on $P$. Then $[X^H,Y^H]$, acting as a derivation, satisfies
\begin{eqnarray*}
\left([X^H,Y^H]f\right)(u)=\lim_{\varepsilon\to 0}\frac{1}{\varepsilon^2}\Big[ f\left(\varphi^{Y^H}_{-\varepsilon}\varphi^{X^H}_{-\varepsilon} \varphi^{Y^H}_\varepsilon\varphi^{X^H}_\varepsilon(u)\right)-f(u)\Big]
\end{eqnarray*}
\end{lemma}
\proof{This follows quickly from the definition of the Lie derivative, and the fact that $X^Hf=\lim_{t\to 0}\frac{1}{t}(f(\varphi^{X^H}_t(u)-f(u))$.
}


That is, $[X^H,Y^H]$ is the horizontal vector field on $P$ given at $u$ by
\begin{eqnarray*}
[X^H,Y^H](u) &=& \frac{d}{dt}\psi_t(u)|_{t=0}\\
\text{where}\ \ \ \psi_{t^2}(u) &=& \varphi^{Y^H}_{-t}\varphi^{X^H}_{-t} \varphi^{Y^H}_t \varphi^{X^H}_t(u)
\end{eqnarray*}
The map $\psi_t(u)$ defines the parallel transport of $u$ around the small rectangle based at $x=\pi(u)$ and so in particular,
\begin{eqnarray*}
\psi_t(u)&=&ug_t\ \ \ \ \text{for }\ g_t\in Hol(A,u),\\
\text{and }\ \ \ \ [X^H,Y^H](u) &=& \frac{d}{dt}(ug_t)|_{t=0}\\
&=& X^*(u)
\end{eqnarray*}
is the horizontal vector at $u$ for $X\in\hol$ according to the definition in Equation \ref{eqn:Xestrella}. Furthermore,
\begin{eqnarray*}
\Omega(X^H,Y^H)&=& -\phi([X^H,Y^H])\\
&=& -X\in \hol.
\end{eqnarray*}
We have proven the following result.
\begin{prop}
Let $(A,\phi)$ be a connection on the $G$-bundle $P$ and let $u\in P$. Then the curvature $\Omega$ of $\phi$ takes values at $u$ in the Lie subalgebra $\hol(A,u)$.
\end{prop}

In Proposition \ref{prop:HolomBasePoint}, we saw that if two points $u,v\in P$ can be joined by a horizontal curve in $P$, then the groups $Hol(A,u)$ and $Hol(A,v)$ coincide and so have the same Lie algebra $\hol$.

\begin{cor}
Let $u\in P$. Then $\Omega_v$ takes values at $v$ in $\hol$ for every $v\in P$ that can be connected to $u$ by a horizontal path.
\end{cor}

\subsection{The reduction theorem}

As previously, we suppose that $\pi:P\to M$ is a principal bundle over $M$ with structure group $G$. Suppose that $(A,\phi)$ is a connection on $P$, and that $H=Hol(A,u)$ is the holonomy group, for a point $u\in P$. If $H\subseteq G$ is a proper subgroup of $G$, then intuitively, there is a certain amount of redundance in considering the connection $A$ on the full bundle $P$. One could feel that all the information regarding the connection is contained in some set smaller than $P$.

\begin{defn}
Let $u_0\in P$. The holonomy bundle of the point $u_0$ is the set $P(u_0)$ of points $v\in P$ that can be connected to $u_0$ by a piecewise smooth horizontal path in $P$.
\end{defn}
Following Proposition \ref{prop:HolomBasePoint}, we can see that $H=Hol(A,u_0)\subseteq G$ acts on $P(u_0)$, and acts transitively on the fibers of the projection to $M$. Also, since $M$ is path-connected, the projection to $M$ is seen to be surjective. Then, in appealing to Proposition \ref{prop:critforsubbundle}, to show that $P(u_0)$ is a subbundle of $P$, we only need to construct local sections of $P$ that take values in $P(u_0)$. For $x\in M$, take a coordinate chart that identifies a small ball $B_\varepsilon(0)\subseteq \R^n$ with a neighbourhood of $x$ (for example, exponential coordinates), such that the origin maps to $x$. Let $u\in \pi^{-1}(x)\cap P(u_0)$. Then, for $y\in B_\varepsilon(0)$, let $\gamma_y$ be the line segment from $0$ to $y$ and let $\sigma(y)=P_{\gamma_y}(u)\in \pi^{-1}(y)\cap P(u_0)$. Then, $\sigma$ defines a section of $P$ that takes values in $P(u_0)$.
\begin{thm}
The holonomy bundle $P(u_0)$ is a subbundle of $P$, with structure group $H=Hol(A,u_0)$. The connection $(A,\phi)$ defines a connection on $P(u_0)$, with the same curvature form and holonomy group.
\end{thm}
\proof{All that needs to be proved is the second statement, but this is clear because the horizontal subspace $A\subseteq TP$ is obviously tangent to $P(u_0)$, since it is spanned by vectors tangent to horizontal curves, and so defines a connection on $P(u_0)$. That the curvature and holonomy do not change are also clear.
}

\subsection{The Ambrose-Singer Theorem}

We now turn to a theorem that has already been referred to, in the example at the end of Section \ref{section:PTonPrincipalBdles}. This gave a way of calculating the holonomy group in terms of the curvature tensor of a connection. Earlier in this section, we described concrete relationships between the curvature and holonomy group of a conection on a prinicpal bundle. In this section, we will state and prove this theorem, and state the equivalent version for vector bundles. We will also make note of similar theorems where we slightly change some hypotheses and definitions.


\begin{thm} \label{thm:AmbroseSinger} Let $\pi:P\to M$ be a principal fiber bundle with structure group $G$ and $(A,\phi)$ a connection on $P$. Then, the Lie algebra $\hol(A,u_0)$ of the holonomy group $Hol(A,u_0)$ of $A$ with base point $u_0\in P$ is spanned by the values $\Omega_v(X,Y)$ taken by the curvature $2$-form, as $v$ varies in the subbundle $P(u_0)$.
\end{thm}
\proof{
To begin, we note that we can we have applied the reduction theorem of the previous section, and that we have $G=Hol(A,u_0)$, and $P=P(u_0)$ so that every point in $P$ can be connected to $u_0$ by a piecewise smooth horizontal curve.
 We define $\mathfrak{h}\subseteq\mathfrak{g}$ to be the subspace spanned by the values $\Omega_v(X,Y)$, for all $v\in P(u_0)$ and $X,Y\in T_vP(u_0)$. We also define $S$ to be the distribution on $P$ spanned at $u\in P$ by
 \begin{enumerate}
\item the horizontal space $A_u$;
\item the set of fundamental vertical vectors $V^*_u$ for $V\in \mathfrak{h}$.
\end{enumerate}
We first claim that $\mathfrak{h}$ is an ideal in the algebra $\mathfrak{g}$. This follows from the invariance properties of the curvature form :
\begin{eqnarray*}
R_g^*\Omega(X,Y) = \Omega(R_{g\ *}X,R_{g\ *}Y) =g^{-1}\cdot\Omega(X,Y) g\in\mathfrak{h}
\end{eqnarray*}
for all $g\in G$. We can use this to show that $S$ is an integrable distribution on $P$. It suffices to show that if $X$ and $Y$ are horizontal vector fields, and if $V$ and $W$ lie in $\mathfrak{h}$, then $[X,Y]$, $[V^*,W^*]$ and $[X,V^*]$ all lie in $S$. The distribution is locally spanned by vector fields of these forms. In the first case, if $X^H$ and $Y^H$ are the horizontal lifts of vector fields $X$ and $Y$ on $M$, then
\begin{eqnarray*}
\Omega(X^H,Y^H)&=&V\in \mathfrak{h}\\
\text{where }\ \ V^*&=& [X^H,Y^H]-[X,Y]^H,
\end{eqnarray*}
so $[X^H,Y^H]\in S$. Secondly, $\mathfrak{h}$ is an ideal of $\mathfrak{g}$, so
\begin{eqnarray*}
[V^*,W^*]=[V,W]^*\in S.
\end{eqnarray*}
Lastly, from Lemma \ref{lemma:bracketofhorizvertvfs}, we have that $[X^H,V^*]=0$. We can then conclude that $S$ is an integrable distribution, and $P$ is foliated by a family of maximal integral submanifolds, including one, $M_{u_0}$, that contains $u_0$. We finally recall that every point in $P$ can be connected to $u_0$ by a piecewise smooth horizontal curve. That is, a curve that is necessarily contained in the integral submanifold $M_{u_0}$. This implies that $P=M_{u_0}$, dim$(TP)$=dim$(S)$ and $\mathfrak{h}=\mathfrak{g}$.
}



%% file: RiemannianHolonomy.tex
\chapter{Riemannian Holonomy}


\section{Holonomy and invariant differential forms}

There are many geometric ways to detect the holonomy group $Hol(g,x)$ of a riemannian manifold. The one most immediately connected to the definition of the holonomy group is that all endomorphisms of a fixed tangent space $T_xM$ to the manifold that arise by parallel transport around closed loops take values in a particular subgroup. Another, as we will see in the following section, is that the universal cover of the manifold admits a proper splitting into the product of lower dimensional manifolds. This fact is essentially equivalent to the holonomy group acting reducibly on $T_xM$, which is to say that it preserves two proper orthogonal subspaces $T_x=V_1+V_2$.
A third method is to consider the induced action of $Hol(g)$ on the associated vector spaces at $x$, other than $T_xM$. In the case that $Hol(g,x)$ preserves a vector $v$, it will lead to a distinguished tensor on $M$.
This section will be short, but it will lead to many ideas that will be seen in the important examples of the following chapter.

We will be primarily interested in the riemannian case, in which the connection is the Levi-Civita connection. At least to start though, we will suppose that $(A,\phi)$ is a connection on the principal frame bundle $\mathcal{F}_{GL}$. By the reduction theorem we can suppose that $(A,\phi)$ is defined on the holonomy bundle, which we denote by $\pi:P\to M$. This bundle has structure group $H=Hol(A,u_0)$.

The standard representation of $GL(n)$ on $\R^n$ is by left multiplication on column vectors. The standard algebraic operations on vector spaces lead to other representations. For example, let $(\R^n)^*$ be the dual space of linear functionals on $\R^n$. Then for $g\in GL(n)$ we define $g:(\R^n)^*\to (\R^n)^*$ by
\begin{eqnarray*}
(g\alpha)(v)=\alpha(g^{-1}v).
\end{eqnarray*}
For $T\in \mathfrak{gl}(n)\cong \R^n\otimes (\R^n)^*$ the representation of $GL(n)$ is given by $g\cdot T=gTg^{-1}$.
Some of the most important representations of $GL(n)$ that are induced from that on $\R^n$ are on $\Lambda^k(\R^n)^*$. For $g\in GL(n)$ and $\varphi\in \Lambda^k(\R^n)^*$. Then,
\begin{eqnarray*}
g: \Lambda^k(\R^n)^*& \to & \Lambda^k(\R^n)^*\\
(g\varphi)(v_1,\ldots v_k) &=& \varphi(g^{-1}v_1,\ldots,g^{-1}v_k)
\end{eqnarray*}
for $v_i\in \R^n$.

If $P\to M$ is a principal $H$-bundle on $M$, and $\rho:H\to \text{Aut}(V)$ is a representation of $H$ on the vector space $V$, the associated vector bundle is the quotient
\begin{eqnarray*}
E=(P\times V)/\sim
\end{eqnarray*}
where $(p,v)\sim (pg, \rho(g)^{-1}v)$ for all $g\in H$. If $\phi$ is a connection for on $P$, we can define a covariant derivative of sections. If the section $\sigma\in \Omega^0(E)$ defines the equivariant function $f:P\to V$, the covariant derivative $\sigma\in \Omega^1(E)$ of $\sigma$ is given by
\begin{eqnarray*}
Df = df+\rho(\phi)f.
\end{eqnarray*}
Here $\rho:\mathfrak{h}\to \text{End}(V)$ is the induced map on the Lie algebra.

\begin{prop}\label{prop:invariantelements}
There exists a one-to-one correspondence between  elements of $V$ that are invariant under the action of $H$, and sections of $E$ that are parallel with respect to the connection $(A,\phi)$.
\end{prop}
\proof{
If $f_0\in V$ satisfies $\rho(g)f_0=f_0$ then $\rho(X)f_0=0$ for all $X\in \mathfrak{h}$. We can consider $f_0$ as a constant $V$-valued function, which determines a section $\sigma$ of $E$. Then, $df_0+\rho(\phi)f_0=0$ and so $\nabla\sigma=0$.

Conversely, a section $\sigma$ of $E$ such that $\nabla\sigma=0$ can equally be thought of as a equivariant function $f:P\to V$ that satisfies $Df=0$. If $\gamma$ is a piece-wise smooth loop in $M$ that is based at $x=\pi(u_0)$ then we can lift it to a horizontal path $\tilde\gamma$ in $P$ that starts at $u_0$ and ends at $u_0g$ for some $g\in Hol(A,u_0)$. Since
\begin{eqnarray*}
\frac{d}{dt}f(\tilde\gamma(t))= Df(\partial_t)-\rho(\phi(\partial_t))f=0
\end{eqnarray*}
(the second term on the right vanishes since the lift is horizontal) we can see that $f$ is constant along the curve $\tilde\gamma$. We can then conclude that $f(u_0g)=f(u_0)$ and
\begin{eqnarray*}
\rho(g)^{-1}f(u_0)=f(u_0)
\end{eqnarray*}
for all $g\in Hol(A,u_0)$ and so the element $f(u_0)\in V$ is invariant under the action of $Hol(A,u_0)$ on $V$.
}

This is an important result in the study of holonomy groups, because many of them admit some invariant vector, and hence parallel tensor field. The usual proof of this result is in the context of connections on vector bundles. We can give this proof easily as well.

Suppose that the connection on the tangent bundle $T_M$ induces a connection on an associated vector bundle $E$. Then The holonomy group $Hol(\nabla,x)$ is a subgroup of $\text{Aut}(T_xM)$ and so acts also on fibre $E_x$ of $E$ at the point $x$. Suppose that $\varphi_0\in E_x$ satisfies $g\varphi_0=\varphi_0$ for all $g\in Hol(\nabla,x)$. For any other point $y\in M$, we take a piece-wise smooth path $\gamma_y$ in $M$ from $x$ to $y$ and define
\begin{eqnarray*}
\varphi(y)=P_{\gamma_y}(\varphi_0)
\end{eqnarray*}
as the parallel translate of $\phi_0$ to $y$. It is immediate from the invariance of $\varphi_0$ that this definition is independent of the path $\gamma_y$, so $\varphi$ is a globally defined field on $M$. It is also immediate from the definition that $\varphi$ is parallel. Conversely, for any parallel field evaluation at $x$ recovers the element $\varphi_0\in E_x$.

The above result allows us to construct an multitude of examples.

\subsubsection{Connections compatible with a metric}

If $(A,\phi)$ is a connection on the frame bundle $\mathcal{F}_{GL}$ such that $H=Hol(A,u_0)$ is compact, then it can be shown that $H$ preserves a positive definite inner product on $\R^n$. Proposition \ref{prop:invariantelements} then shows that there exists a metric on $T_M$ such that $\nabla g=0$. In this case, for any local vector fields $X,Y,Z$ on $M$,
\begin{eqnarray*}
Zg(X,Y)=g(\nabla_ZX,Y)+g(X,\nabla_ZY).
\end{eqnarray*}
If $\nabla$ is also torsion-free, then it is the Levi-Civita connection of the metric $g$. If $u_0$ is orthonormal with respect to $g$, then the holonomy bundle is contained in the principal bundle $\mathcal{F}_O$ of $g$-orthonormal frames.

\subsubsection{Complex and Almost Complex Structures}

If we suppose that $\R^{2n}$ has the basis $\{e^1,\ldots, e^{2n}\}$ then we can define an endomorphism $J_0:\R^{2n}\to \R^{2n}$ by
\begin{eqnarray*}
J_0e^j &=& e^{j+n}\\
J_0e^{j+n} &=& -e^j
\end{eqnarray*}
for $j=1,\ldots, n$. Then $J_0^2=-\text{Id}$ and $J_0$ defines the complex structure on $\C^n\subseteq \R^{2n}\otimes \C$ given by multiplication by $i=\sqrt{-1}$. $\C^n$ has basis $\{\varepsilon^j=1/2(e^j-iJe^j)\}$ for $j=1,\ldots,n$.
The subgroup of $GL(2n,\R)$ that preserves $J_0$ (i.e. every element satisfies $gJ_0g^{-1}=J_0$) is $GL(n,\C)$.

If $\nabla$ is a connection on $TM$ such that $Hol(\nabla,u_0)\subseteq GL(n,\C)$ then from the above construction there exists a tensor field $J$ of endomorphisms of $TM$ such that $\nabla J=0$. This means that for $v\in TM$,
\begin{eqnarray*}
\nabla(Jv)=J\nabla v.
\end{eqnarray*}
It can easily be seen that $J^2=-\text{Id}$. We can therefore consider the complexification of the tangent bundle $T_\C M=TM\otimes \C$ and see that it decomposes into the subbundles $T^{(1,0)}_M$ and $T^{(0,1)}_M$ of $+i$ and $-i$ eigenvectors of $J$ respectively.
\begin{exer}
Show that if $\nabla J=0$, then the connection $\nabla$ extends to a connection on the bundle $T^{(1,0)}_M$.
\end{exer}
The existence of a connection on $T^{(1,0)}_M$ is not particularly interesting, but if we take into account  the torsion of $\nabla$ we can obtain new information. Given an almost-complex structure $J$ on the manifold $M$, we define the \emph{Nijenhuis tensor} of $J$ to be
\begin{eqnarray*}
N_J(X,Y)=\frac{1}{4}([JX, JY ]- J[JX, Y ]-J[X, JY ]-[X, Y ]).
\end{eqnarray*}
Then, by expressing, for example, $[JX,Y]=\nabla_{JX}Y-\nabla_YJX$ one can easily see that if a $2n$-dimensional manifold $M$ admits a torsion-free connection on $TM$ with holonomy group contained in $GL(n,\C)$, then the Nijenhuis tensor of the induced almost-complex structure vanishes identically. It is a famous theorem of Newlander and Nirenberg (see \cite{NN}) that this implies that the almost complex structure is integrable, and arises from an atlas of maps from $\C^n$ with holomorphic transition functions.

\subsubsection{Hermitian Manifolds}

The terminology for almost-complex manifolds that admit a compatible riemannian metric is varied. In the symplectic context, they are sometimes referred to as \emph{almost-symplectic} structures and in the case that the almost complex structure is integrable, the are called hermitian structures. We will discuss these now.

If we suppose that $M$ supports an almost complex structure $J$, a riemannian metric is hermitian (with respect to $J$) if for any $X,y\in TM$,
\begin{eqnarray*}
g(JX,JY)=g(X,Y).
\end{eqnarray*}
In this case we can define the bilinear form $\omega$ by
\begin{eqnarray}\label{eqn:symp-metric}
\omega(X,Y)=g(JX,Y)
\end{eqnarray}
and it can easily be seen that $\omega$ is skew-symmetric, since $\omega(X,X)=0$. On $\R^{2n}$, the subgroup of $GL(2n, \R)$ that preserves the euclidean inner product $g_0$ is $SO(2n)$ and the subgroup that preserves the standard almost complex structure $J_0$ is $GL(n,\C)$. The subgroup that preserves the standard symplectic form $\omega_0$ (without necessarily preserving $g_0$ or $J_0$ is the symplectic group $Sp(2n,\R)$. The group of transformations of $\R^{2n}$ that preserve any two of these three is the group $U(n)$ of unitary matrices, since
\begin{eqnarray*}
U(n)=SO(2n)\cap GL(n,\C)\cap Sp(2n,\R),
\end{eqnarray*}
with the same intersection if we take any two of these three groups.
A manifold $M$ that supports an almost complex structure $J$ and a compatible hermitian metric $g$ then admits a reduction of its principal frame bundle to $\mathcal{F}_{U}$, which has structure group $U(n)$. To do this we  only consider those linear maps $u:\R^{2n}\to T_xM$ such that
\begin{eqnarray*}
u(J_0v)=Ju(v),\ \ \ \ &&\ \ \ u^*g=g_0.
\end{eqnarray*}
Conversely, if $\mathcal{F}_{GL}$ admits a connection $(A,\phi)$, or equivalently we have $\nabla$ on $TM$, with holonomy group contained in $U(n)$, then from the holonomy bundle, we obtain an almost complex structure $J$ and hermitian metric $g$ on $M$ that satisfy $\nabla g=0$ and $\nabla J=0$. The metric and endomorphism $J$ also determine a non-degenerate $2$-form $\omega$ that satisfies $\nabla\omega=0$.

If $\nabla$ is torsion-free, then $\nabla$ is necessarily the Levi-Civita connection of the metric $g$. From the discussion above, the almost complex structure is integrable and $d\omega=0$. In this situation, we say that the $g$ is a K\"ahler metric on the complex manifold $(M,J)$.

This example is typical among the riemannian holonomy groups. In each case, the holonomy group $Hol(g)$ of the Levi-Civita connection preserves one or more tensors, in addition to the metric, when we consider the representations on $\Lambda^*(\R^n)^*$. In the K\"ahler case this is the sympletic form $\omega_0$. In the Calabi-Yau case it is a complex volume form $\Omega_0\in \Lambda^{n,0}$. In te case of  quaternion-K\"ahler manifolds, where $Hol(g)\subseteq Sp(n)Sp(1)\subseteq SO(4n)$, the holonomy group preserves the $4$-form 
\begin{eqnarray*}
\Phi_0=\frac{1}{2}\left( \omega_I^2+\omega_J^2+\omega_K^2\right)
\end{eqnarray*}
where $\omega_I$, $\omega_J$ and $\omega_K$ are $2$-forms on $\Ham^n$ that correspond to the almost-complex structures given by right multiplication by $i$, $j$ and $k$ respectively. This idea can also be used to understand the other important holonomy groups, but we will leave these to the following section. 


\section{Reducible holonomy groups}\label{sec:reducible}

In this section, we develop some ideas that appeared in the previous one. In that case, we considered the representation of $H=Hol(g)$ on $\Lambda^*(\R^n)^*$ and supposed that there was some element $\alpha_0$ that was left invariant : $g\cdot \alpha_0=\alpha_0$ for all $g\in H$. We were then easily able to construct a tensor $\alpha $ on $M$ was parallel and coincided with $\alpha_0$ at at a fixed point.  In this section we consider the standard representation of $H$ on $\R^n$ and suppose that there is a subspace $V\subseteq \R^n$ that is preserved by $H$. Since $H\subseteq SO(n)$, the orthogonal complement of $V$ is also preserved so we can suppose that $\R^n=V_1+V_2$ and $H\subseteq SO(V_1)\times SO(V_2)$. The conclusion that we arrive at is that if $M$ is simply-connected, then $M$ globally splits as a riemannian product $M=M_1\times M_2$ with metric $g=g_1+g_2$ and $Hol(g)=Hol(g_1)\times Hol(g_2)$. This is the result that we wish to prove in this section. As in other sections, the discussion that we give follows that from \cite{KN}, although the result is originally due to de Rham \cite{deRham}.

The most lengthy part of this section is the proof of Theorem \ref{thm:dRprojwelldefined}, and most of of our time will be spent proving this result. This theorem says that we can give a well-defined map $p:M\to M_1\times M_2$. We will then quickly be able to see that it is an isometry.

We suppose in this section that $(M,g)$ is a connected and simply-connected riemannian manifold, of dimension $n$, and such that the metric is complete. We will also suppose that the holonomy group $Hol(g,x_0)\subseteq SO(T_{x_0}M)$ preserves two orthogonal subspaces $V_1(x_0)$ and $V_2(x_0)=V_1(x_0)^\perp\subseteq T_{x_0}M$.  Let $V_1$ be of rank $k$ and $V_2$ be of rank $n-k$.

We define two distributions on $M$ as follows. For $y\in M$, let $\gamma$ be a path in $M$ from $x_0$ to $y$. We set $V_1(y)=P_\gamma(V_1(x_0))$ where $P_\gamma$ is the parallel translation map from $T_{x_0}M$ to $T_yM$. Similarly, $V_2(y)=P_\gamma(V_2(x_0))$. We obtain some preliminary information from the following sequence of lemmas..
\begin{prop}
\begin{enumerate}
\item The distributions $V_1$ and $V_2$ are well-defined and smooth.
\item The distributions $V_i$ are involutive.
\item The maximal integral submanifolds of $V_1$ and $V_2$ are totally geodesic and complete in the induced metric.
\end{enumerate}
\end{prop}
We denote by $M_1(y)$ and $M_2(y)$ the maximal integral submanifolds of the two distributions that pass through the point $y\in M$.

\proof{We show that the distributions do not depend upon the paths chosen. For any two paths $\gamma$ and $\gamma$ from $x_0$ to $y$, $\gamma^{-1}\cdot\sigma$ is a closed loop based at $x_0$, and so we recall that $P_{\gamma^{-1}\cdot\sigma}$ preserves the subspace $V_1(x_0)\subseteq T_{x_0}M$. Hence,
\begin{eqnarray*}
P_\gamma(V_1(x_0))=P_\gamma\cdot P_{\gamma^{-1}\cdot\sigma}(V_1(x_0))=P_{\sigma}(V_1(x_0)).
\end{eqnarray*}
Secondly, if $X$ and $Y$ are local sections of $V_1$, we wish to show that $[X,Y]$ is a local section of $V_1$ as well. Since $\nabla$ is torsion free, it is sufficient to show that $\nabla_XY$ is contained in $V_1$ as well, but this follows from Proposition \ref{prop:partranstoconnection}, which expresses the connection as an infinitesimal parallel transport.

Finally, let $\gamma$ be a geodesic in $M$, with $\gamma^\prime(0)$ contained in $V_1(y)$. Then, the tangent vector $\gamma^\prime(t)$ is given by parallel translation along $\gamma$, of the initial tangent vector $\gamma^\prime(0)$ and is hence then contained in $V_1(\gamma(t))$. The curve $\gamma$ is then an integral submanifold of the distribution and so contained in $M_1(y)$. This submanifold is hence totally geodesic. Since the geodesics in the induced metric coincide exactly with those in the ambient space, that start in $M_1(y)$, they are defined for all values of $t$. The space is therefore complete.
}

\begin{lemma}
On a neighbourhood of any point $y\in M$ there exists a local coordinate system $(x^1,\ldots, x^n)$ such that $(\partial/\partial x^1,\ldots,\partial/\partial x^k)$ locally spans the distribution $V_1$ and $(\partial/\partial x^{k+1},\ldots,\partial/\partial x^n)$ locally spans $V_2$. The maximal integral submanifolds $M_1$ are then locally given as the common level sets $x^i=c^i$ for $k+1\leq i\leq n$ and the $M_2$ are the levels sets of the coordinates $x^j$ for $1\leq j\leq k$.
\end{lemma}
\proof{
Since the distributions are involutive, there exist local coordinates $(u^i)$ and $(v^j)$ such that $\partial/\partial u^1,\ldots,\partial/\partial u^k$ locally span $V_1$ and such that $\partial/\partial v^{k+1},\ldots,\partial/\partial v^n$ locally span $V_2$. Then, since $V_1$ and $V_2$ are transverse at a point $y$, the map
\begin{eqnarray*}
(x^1,\ldots,x^n)=(u^1,\ldots,u^k,v^{k+1},\ldots,v^n):U\to M
\end{eqnarray*}
has derivative to $T_yM$ an isomorphism, so locally defines the coordinates that we need.
}

\begin{prop} Let $y\in M$ and let $M_1(y)$ and $M_2(y)$ be the maximal integral submanifolds for the two distributions that contain $y$. Then there is a neighbourhood $U$ of $y$ diffeomorphic to $U_1\times U_2$, for $U_i\in M_i(y)$, and such that the metric $g$ on $U$ is isometric to the product metric on $U_1\times U_2\subseteq M_1(y)\times M_2(y)$.
\end{prop}
\proof{
The fact that there exist nighbourhoods of the form $U=U_1\times U_2$ where the slices are integral submanifolds of the two distributions follows from the previous lemma. To show that the metric on the local neighbourhood is the product we show that for $i\leq i,j\leq k$, $g(\partial_i,\partial_j)$ does not depend on $x^l$, where $k+1\leq l\leq n$. We have,
\begin{eqnarray*}
\frac{\partial}{\partial x^l} g(\partial_i,\partial_j) &=& g(\nabla_{\partial_l}\partial_i,\partial_j)+g(\partial_i,\nabla_{\partial_l}\partial_j)\\
&=&g(\nabla_{\partial_i}\partial_l,\partial_j)+g(\partial_i,\nabla_{\partial_j}\partial_l)\\
&=& \frac{\partial}{\partial x^i}g(\partial_l,\partial_j)-g(\partial_l,\nabla_{\partial_i}\partial_j) \\
&&\ \ \ \ \ \ +\frac{\partial}{\partial x^j}g(\partial_i,\partial_l)- g(\nabla_{\partial_j}\partial_i,\partial_l)
\end{eqnarray*}
and all terms on the right hand side of this equation vanish because $\partial_l$ is everywhere orthogonal to $V_1$, and $\nabla_{\partial_i}\partial_j$ is contained in $V_1$ when $\partial_i$ and $\partial_j$ are.
}

\begin{defn}\label{defn:derhamprojection}
Let $\gamma$ be a smooth curve in $M$ starting at $x\in M$. We will define the projection of $\gamma$ onto $M_1(x)$ to be a curve $\tilde \gamma$ in the integral submanifold $M_1(x)$ defined as follows. Given $\gamma:I\to M$, define $\gamma_t(s)=\gamma(ts)$. $\gamma_t$ is a curve in $M$ from $x$ to $\gamma(t)$. Let $P_{\gamma_t^{-1}}$ be the parallel transport map from $T_{\gamma(t)}M$ to $T_xM$, and set
\begin{eqnarray*}
\sigma(t)=P_{\gamma_t^{-1}}\left(\pi^{V_1}({\gamma}^\prime(t)\right)
\end{eqnarray*}
where $\pi^{V_1}({\gamma}^\prime(t))$ is the component of the derivative of the curve in the $V_1$-subspace. Then $\sigma$ defines a curve in the subspace $V_1(x)\subseteq T_xM$. We define the curve $\tilde\gamma$ starting at $x$ in the integral submanifold $M_1(x)$ to be defined by the fact that
\begin{eqnarray}\label{eqn:development}
\tilde\gamma^\prime(t)=P_{\tilde\gamma_t}(\sigma).
\end{eqnarray}
\end{defn}
That is, parallel translation along $\tilde\gamma$ to $\tilde\gamma(t)$ of the vector $\sigma(t)$ gives the velocity vector of the curve. This notion is related to the \emph{development} of a curve, which arises from the parallel translation in the bundle of \emph{affine} frames to $M$. This is studied in detail in \cite[Sec. III.4]{KN}. In particular, if $(M,g)$ is a complete riemannian manifold, and $\sigma$ is a curve in $T_xM$, there is a curve $\tilde\gamma$ in $M$ that satisfies Equation \ref{eqn:development} and is defined for the same values of $t$ that $\sigma$ is (see \cite[Thm. IV.4.1]{KN}).

\begin{prop}
If $M$ is equal to the product $M=M_1\times M_2$ then the projection as defined above coincides with the usual projection onto the first factor.
\end{prop}

\proof{
This is reasonably obvious, intuitively. If a path $\mu$ is given in the product as $\mu=(\mu_1,\mu_2)$, then parallel transport along $\mu$ of vectors tangent to the first factor $M_1$ is given by parallel translation in along $\mu_1$, together with trivial translation in $M_2$. The projection map is taken by parallel translation along $\mu$, and then back along $\mu_1$.
}

\begin{thm}\label{thm:dRprojwelldefined}
If $\gamma_1$ and $\gamma_2$ are curves in $M$ from $x_0$ to $y$ that are homotopic, then the projection map to $M_1(x_0)$ defined using the two curves have the same endpoint in $M_1(x_0)$.
\end{thm}
This in particular means that the projection $p_1:M\to M_1(x_0)$ is well-defined. We can equally well define the map $p_2:M\to M_2(x_0)$ and hence the two maps together $p=(p_1,p_2):M\to M_1(x_0)\times M_2(x_0)$. This is the principal result of this section and we will spend much of the rest of the section proving it.


The important local phenomenon that we will use is the local decomposition $U=U_1\times U_2$ in a neighbourhood of any $y\in M$. We suppose that $\gamma_1$ and $\gamma_2$ are homotopic. Then, using a homotopy $F:I\times I\to M$ we can make a series of adjustments to $\gamma_1$, with each taking place in a good neighbourhood $U$. The curves $\gamma^1,\ldots, \gamma^N$ can be taken so that $\gamma^1=\gamma_1$ and $\gamma^N=\gamma_2$. and so that $\gamma^i$ and $\gamma^{i+1}$ coincide except on a small interval with images contained in $U=U_1\times U_2$. We can therefore suppose that $\gamma_1$ and $\gamma_2$ almost coincide :
\begin{eqnarray*}
\gamma_1 &=& \kappa\cdot \mu\cdot \tau\\
\gamma_2 &=& \kappa\cdot \nu\cdot \tau
\end{eqnarray*}
where $\tau$ is a curve from $x_0$ to $z_1$, $\mu$ and $\nu$ are curves that connect $z_1$ and $z_2$, and $\kappa$ is a curve from $z_2$ to $y$. We suppose that $\mu$ and $\nu$ are contained in $U\cong U_1\times U_2$ where $U_i\subseteq M_i(z_1)$.

\begin{lemma}
Given a curve $\gamma:I\to M$ in $M$ that starts at $x_0$, and $a\in I$, define the curve $\gamma_a$ by $\gamma_a(t)=\gamma(t)$ for $t\in [0,a]$ and $\gamma_a(t)$ for $t\in [a,1]$ is given by the projection, according to the Definition \ref{defn:derhamprojection}, of $\gamma|_{[a,1]}$ onto the integral submanifold $M_1(\gamma(a))$.

Then, the projection of $\gamma$ onto $M_1(x_0)$ coincides with the projection of $\gamma_a$ on $M_1(x_0)$.
\end{lemma}
That is, given a curve $\gamma$ the projection onto $M_1(x_0)$ coincides with the composition of two intermediate projections, one for part of the curve onto another integral submanifold, and then the resulting curve onto $M_1(x_0)$.

We omit the proof of this lemma, and refer the interested reader to \cite{KN}.

Following this lemma we suppose that the curve $\kappa$ given above is actually contained in the integral submanifold $M_1(z_2)$. In essence, the projection of $\gamma_1$ onto $M_1(x_0)$ is by a sequence of intermediate projections, and we suppose that $\kappa$ is already the result of a projection.

In $U$, the curve $\mu$ is given by $(\mu_1,\mu_2)\in U_1\times U_2$, and the projection of $\mu$ onto $M_1(z_1)$ is $\mu_1$.

By taking $U$ sufficiently small, we can suppose that both $U_1$ and $U_2$ are geodesically convex. Let $\mu^*$ be the geodesic in $U$ from $z_2$ to the endpoint of the curve $\mu$. Then $\mu^*$ is contained in $M_2(z_2)$ and parallel transport along $\mu^*$ of $V_1$-vectors is the identity, in the trivialisation of $V_1$ along $M_2(z)$. The curves $\mu^{-1}=(\mu_1^{-1},\mu_2^{-1})$ and $\mu_1^{-1}\cdot \mu^*$ from $z_2$ to $z_1$ therefore determine the same parallel transport of vectors in the distribution $V_1$.

We now wish to lower the curve $\kappa$ from the submanifold $M_1(z_2)$ to the submanifold $M_1(z_1)$. To do this we need the following lemmas. We leave the proofs to later.

\begin{lemma}
Let $\tau:I\to M$ be a curve in $M_1(x)$, and let $\sigma:[0,s_0]\to M$ be an arc-parameterised geodesic in $M_2(x)$ with $\sigma(0)=\tau(0)$ and $Y_0=\sigma^\prime(0)$ . Let $Y_t$ be the vector field along $\tau$ obtained by parallel of $Y_0$. Then, the homotopy $f(t,s)=\exp_{\tau(t)}(sY_t)$ satisfies the following :
\begin{enumerate}
\item $f(t,0)=\tau(t)$, $f(0,s)=\sigma(s)$,
\item parallel translation along the loop given as $f$ restricted to the boundary of $I\times [0,s_0]$ is trivial,
\item the vector $f_t(t,s)$ is parallel to $\tau^\prime(t)$ along the curve $s\mapsto f(t,s)$,
\item the vector $f_s(t,s)$ is parallel to $\sigma^\prime(s)$ along the curve $t\mapsto f(t,s)$.
\end{enumerate}
Moreover, $f(t,s)=\exp(sY_t)$ is the unique homotopy that satisfies these properties.
\end{lemma}

\begin{lemma}
The projection, in the sense considered above, of the curve $\tau\cdot \sigma^{-1}$ onto the integral submanifold $M_1(\sigma(s_0))$ is given by $t\mapsto f(t,s_0)$.
\end{lemma}

We apply these lemmas, in the case that $\tau$ is the curve $\kappa$, and $\sigma$ is the geodesic $\mu^*$. We can conclude that the projection of $\kappa\cdot\mu$ is given by $\mu_1$, followed by the curve $\kappa^\prime$ that is obtained from the homotopy above. We return to also consider the curve $\nu$ between $z_1$ and $z_2$. Then $\nu=(\nu_1,\nu_2)$, and since the endpoints coincide with those of $\mu$, the endpoint of $\nu_1\in M_1(z_1)$ coincides with that of $\mu_1$, so we can take the same geodesic $\mu^*$ as we did above. It follows that the homotopy that carried the curve $\kappa$ to $\kappa^\prime$ is the same, and so the projection of $\kappa\cdot\nu$ onto $M_1(z_1)$ is given by $\kappa^\prime\cdot\nu_1$.

We now deal with the curve $\tau$ from $x_0$ to $z_1$. We can decompose it as the product of curves $\tau^1,\ldots,\tau^k$ where each is contained in an open set of the form $U=U_1\times U_2$. Then, as in the previous case, we can project $\kappa^\prime\cdot\mu_1$ and $\kappa^\prime\cdot\nu_1$ onto the maximal integral submanifold $M_1$ that passes through the initial point of $\tau^k$. From the homotopy lemma, this is dependent on the parallel transport of a $V_2$-vector (given there as $Y_0$) along the $V_1$ path. Parallel transport of $V_2$-vectors along curves tangent to $V_1$ is only dependent on the endpoints. Since the curves $\kappa^\prime\cdot\mu_1$ and $\kappa^\prime\cdot\nu_1$ eventually coincide, we can conclude the same thing, for the projection of $\kappa^\prime\cdot\mu_1\cdot\tau^k$ and $\kappa^\prime\cdot\nu_1\cdot\tau^k$ on the maximal integral submanifold $M_1$.

 We can continue this argument, for the other curves $\tau^i$, and conclude that the projections of $\gamma_1=\kappa\cdot\mu\cdot\tau$ and $\gamma_2=\kappa\cdot\nu\cdot\tau$ onto the submanifold $M_1(x_0)$ have the same endpoints. This concludes the proof of Theorem \ref{thm:dRprojwelldefined}.

It goes without saying that there is no distinguished position in this arguement of the first factor over the second. We can therefore give a well-defined projection $p_2:M\to M_2(x_0)$, and hence two together $p=(p_1,p_2):M\to M_1(x_0)\times M_2(x_0)$.

\begin{thm}
\label{t.splitting}
The map
$p=(p_1,p_2):M\to M_1(x_0)\times M_2(x_0)$ is an isometry.
\end{thm}
\proof{One can show that if $f:M\to N$ is an isometric immersion of riemannian manifolds of the same dimension, and if $M$ is complete, then $f$ is a covering map. In the current context though, if $h$ is a homotopy in $M$ of paths contained in $M_1(x_0)$, then $p_1\circ h$ is a homotopy in $M_1(x_0)$. We conclude that $M_1(x_0)$ and $M_2(x_0)$ are also both simply connected. If we can show that $p=(p_1,p_2)$ is an preserves lengths of tangent vectors, then we can conclude that it is an isometry.

At $y\in M$, let $X=Y+Z\in V_1+v_2=T_yM$. Then, $p_{1*}(X)$ is given by the parallel translation of $Y$ along a path $\gamma$ from $y$ to $x_0$, followed by parallel translation along $\tilde{\gamma}$ to $p_1(y)$. This preserves the length of $Y$, and $p_{2*}$ preserves the length of $Z$. $Y$ and $Z$ are orthogonal, so we can conclude that $p$ preserves the lengths of tangent vectors.
}

\begin{thm}
\label{t.Hol0iscpt}
Let $(M, g)$ be a riemannian manifold of dimension $n$.
Then $Hol^0(g)$ is a closed, connected Lie subgroup of $SO(n)$.
\end{thm}

This theorem is a consequence of two results. The first, Theorem \ref{t.splitting},
which tells us that if $(M, g)$ is a complete, simply-connected Riemannian manifold,
then it is isometric to  a product manifold $(M_1 \times \cdots \times M_k, g_1 \times \cdots \times g_k)$
such that the holonomy representation  of $Hol(g_j)$, for each metric $g_j$, is 
irreducible. 

The second result we need is the following theorem, whose proof can be found at 
\cite{KN}, Appendix 5.
\begin{thm}
\label{t.KNappendix}
Let $G$ be a connected Lie subgroup of $SO(n)$
which acts irreducibly in $\RR^n$. Then, $G$ is closed in $SO(n)$.
\end{thm}
Combining Theorems \ref{t.Hol0iscpt} and \ref{t.KNappendix}, we see that
$Hol^0(g)$ is a compact subgroup of $SO(n)$. Cheeger and Gromoll \cite{CheegerGromoll}
proved the following result.
\begin{thm}
\label{thm.Hol0iscpt}
Let $(M, g)$ be a compact irreducible riemannian manifold of dimension $n$.
Then $Hol(g)$ is a compact  Lie subgroup of $SO(n)$.
\end{thm}
Their proof relies on the study of the fundamental group of a compact, irreducible riemannian manifold.


\section{Decompositions of the curvature tensor}

It is a well known fact in riemannian geometry that there are a number of different types of curvature of which one can speak on a riemannian manifold. The strongest and most restrictive is the sectional curvature, which is defined in terms of two-dimensional planes tangent to the manifold. The Ricci curvature, in a given direction, is ($n-1$-times) the average of the sectional curvatures of planes that contain the direction. The scalar curvature is ($n$-times) the average of the Ricci curvatures at that point. We won't deal here with the sectional curvature, but rather the curvature operator $R_{\cdot\cdot}\cdot$. The Ricci and scalar curvature can be considered components of the full curvature operator, in a decomposition that we will describe, and this leads us to the remaining term, the Weyl curvature.

The frame bundle perspective allows us to consider these tensors as functions with values in a vector space, and we will primarily consider the linear algebra of that vector space. Throughout this section we should recall the equivalence between sections of a vector bundle, and equivariant functions on the correct principal bundle.

\subsection{Abstract curvature tensors}

We suppose that $(M,g)$ is a riemannian manifold and $\mathcal{F}_O$ is the orthogonal frame bundle. In Section \ref{sec:torsion} we saw that there exists a unique connection $(A,\phi)$ on $\mathcal{F}_O$ that is torsion-free. It is uniquely determined by the equation
\begin{eqnarray*}
d\omega=-\phi\wedge\omega
\end{eqnarray*}
where $\omega:T\mathcal{F}_O\to \R^n$ is the canonical $\R^n$-valued $1$-form on $\mathcal{F}_O$.

The riemannian curvature tensor is without exageration the fundamental object of study in riemannian geometry. The symmetries of the curvature tensor are well-known, but we will discuss them briefly here. We will consider the curvature in one of two ways. Firstly, we will consider $\nabla$ as the Levi-Civita connection on $TM$, with curvature $R\in\Omega^2(End(TM))$
\begin{eqnarray*}
R_{XY}Z=\nabla_{[X,Y]}Z-[\nabla_X,\nabla_Y]Z.
\end{eqnarray*}
This differs from the definition that we made in earlier sections by a factor of $-1$. We make this definition to coincide with usual conventions (see \cite{Besse}).
Otherwise we will consider the curvature as an equivariant, horizontal $2$-form on $\mathcal{F}_O$ with values in $\mathfrak{so}(n)$, defined by
\begin{eqnarray*}
\Omega=d\phi+\phi\wedge\phi.
\end{eqnarray*}
These two viewpoints will be related according to the associated bundle construction considered in Section \ref{sec:principalbundles}. We start with the first and secong {\em Bianch identities}.
\begin{prop}
let $(A,\phi)$ be a connection on $\mathcal{F}_O$, with curvature form $\Omega$ and torsion $\Theta$. Then,
\begin{enumerate}
\item $D\Theta=\Omega\wedge\omega$,
\item$D\Omega=0$.
\end{enumerate}
\end{prop}
Here $D=D^\phi$ is the differential operator define by $\phi$ in Equation \ref{eqn:defnofD}. In particular in the first statement we have that,
\begin{eqnarray*}
D\Theta(X,Y,Z) =\Omega(X,Y)\omega(Z)+\Omega(Y,Z)\omega(X)+\Omega(Z,X)\omega(Y).
\end{eqnarray*}
\begin{cor}
Let $\nabla$ be the Levi-Civita connection on $TM$ (that is, such that $T^\nabla=0$). Then,
\begin{enumerate}
\item $R_{XY}Z+R_{YZ}X+R_{ZX}Y=0$,
\item $(\nabla_XR)_{YZ}T+(\nabla_YR)_{ZX}T+(\nabla_ZR)_{XY}T=0$.
\end{enumerate}
\end{cor}
The second Bianchi identity, as given in each of the above results, holds for an arbitrary conection on a bundle (see \cite{Laws}). This can be expressed as $d^\nabla R^\nabla=0$, in the notation of Section \ref{sec:CurvToronVBs}.

\proof{
We have that $\Theta=d\omega+\phi\wedge\omega$, so
\begin{eqnarray*}
D\Theta &=& d\Theta+\phi\wedge\Theta\\
&=& d\phi \wedge\omega -\phi\wedge d\omega+\phi\wedge d\omega +\phi\wedge\phi\wedge\omega\\
&=& -\phi\wedge\phi\omega +\Omega\wedge\omega+\phi\wedge\phi\wedge\omega\\
&=&\Omega\wedge\omega.
\end{eqnarray*}
For the second equation we have that
$D\omega = d\Omega +(\rho_*\phi)\wedge\Omega$ where $\rho$ is the adjoint representation. That is,
\begin{eqnarray*}
D\Omega &=& d\Omega +\phi\wedge\Omega-\Omega\wedge\phi\\
&=& d\phi\wedge\phi -\phi\wedge d\phi +\phi\wedge(d\phi +\phi\wedge\phi) -(d\phi+\phi\wedge\phi)\wedge\phi\\
&=& 0.
\end{eqnarray*}
The corollary is immediate, in the terminology of the associated bundle construction.
}

It should be noted from the proof that in the first case, we do not at any point use the fact that the connection preserves the metric; it holds for any connection on $\mathcal{F}_{GL}$. In the second case, the proof shows that this holds for any connection, as stated above. The simple proofs to these results also demonstrate one of the advantages of using the principal bundle formalism. Instead of a potentially more complicated calculation on a vector bundle, we can simply use exterior calculus to derive the same result.

The riemannian curvature tensor $R$ can be seen to be a $2$-form on $M$ with values in the bundle $\mathfrak{so} (T_M)$ of skew-symmetric transformations of $T_M$. According to the correspondence that we have seen in previous section, this corresponds to a $2$-form on $\mathcal{F}_O$ with values in $\mathfrak{so}(n)$ that is equivariant and horizontal in a certain sense. It can also be considered as a section of the bundle $\Lambda^2T_M^*\otimes \mathfrak{so}(T_M)$. According to the same correspondence, this determines an equivariant function $\mathcal{R}$ on $\mathcal{F}_O$ with values in $\Lambda^2\R^{n*}\otimes \mathfrak{so}(n)$. This relates to the $2$-form $\Omega$ by
\begin{eqnarray}\label{eqn:curvfunction}
\mathcal{R}_{\pi(X)\pi(Y)}=-\Omega(X,Y)\in\mathfrak{so}(n)
\end{eqnarray}
for $X,Y\in T\mathcal{F}_O$. By making the identification $\Lambda^2\cong \mathfrak{so}(n)$, we will study the vector space $\Lambda^2\otimes\Lambda^2$ and determine the subset that $\mathcal{R}$ takes values in.

\begin{prop}
$\mathcal{R}$ takes values in the symmetric product $S^2(\Lambda^2)\subseteq \Lambda^2\otimes\Lambda^2$.
\end{prop}
\proof{ The first Bianchi implies that $\mathcal{R}$ satisfies
\begin{eqnarray*}
\mathcal{R}_{xy}v+\mathcal{R}_{yv}x+\mathcal{R}_{vx}y=0
\end{eqnarray*}
for $x,y,v\in \R^n$. Therefore, together that $\mathcal{R}$ is anti-symmetric in the respective arguments,
\begin{eqnarray*}
\langle \mathcal{R}_{xy}v,w\rangle &=& -\langle \mathcal{R}_{yv}x+\mathcal{R}_{vx}y,w\rangle\\
 &=& \langle \mathcal{R}_{yv}w,x\rangle +\langle \mathcal{R}_{vx}w,y\rangle\\
 &=& -\langle \mathcal{R}_{vw}y+\mathcal{R}_{wy}v,x\rangle -\langle \mathcal{R}_{xw}v+\mathcal{R}_{wv}x,y\rangle\\
 &=& 2\langle \mathcal{R}_{vw}x,y\rangle +\langle \mathcal{R}_{wy}x+\mathcal{R}_{xw}y,v\rangle\\
 &=& 2\langle \mathcal{R}_{vw}x,y\rangle -\langle \mathcal{R}_{yx}w,v\rangle,
 \end{eqnarray*}
 so $\langle \mathcal{R}_{vw}x,y\rangle =\langle \mathcal{R}_{xy}v,w\rangle$ as desired.
}

We consider the map
\begin{eqnarray*}
\beta:S^2(\Lambda^2)&\to &\otimes^4(\R^n)^*\\
\beta(A)_{vwxy}&=&\langle A_{vw}x+A_{wx}v+A_{xv}w,y\rangle.
\end{eqnarray*}
If $w=v$, then $\beta(A)_{vwwy}=0$ so we can see that $\beta$ takes values in $\Lambda^4$.
\begin{prop}
$\beta$ is given by
\begin{eqnarray*}
\beta:\phi\otimes\phi\mapsto \frac{1}{6}\phi\wedge\phi.
\end{eqnarray*}
\end{prop}
\proof{
We have that
\begin{eqnarray*}
\beta(\phi\otimes\phi)_{vwxy}=\phi(v,w)\phi(x,y)+\phi(w,v)\phi(v,y)+\phi(x,v)\phi(w,y)
\end{eqnarray*}
which equals $\phi\wedge\phi$ up to scale.
}

\begin{defn}
We define the space of abstract curvature tensors in $n$-dimensions to be
\begin{eqnarray*}
{\bf R}_n=S^2(\Lambda^2)\cap\ker \{\beta:S^2(\Lambda^2)\to\Lambda^4\}.
\end{eqnarray*}
\end{defn}
Then, if $(M,g)$ is a riemannian manifold and $\mathcal{F}_O$ is the orthonormal frame bundle over $M$, then the curvature defines a $O(n)$-equivariant function $\mathcal{R}$ on $\mathcal{F}_O$ with values in ${\bf R}_n$.

\begin{exer}
Show that ${\bf R}_n$ is a vector space of dimension $\frac{n^2(n^2-1)}{12}$.
\end{exer}
We can therefore see that ${\bf R}_2$ is $1$-dimensional, which is seen in the fact that the curvature tensor on a surface is determined by the gaussian curvature. We also have that ${\bf R}_3$ is $6$-dimensional and ${\bf R}_4$ is $10$-dimensional. Coincidentally, perhaps, this equals the dimension of the space of symmetric bilinear forms on $\R^3$ and $\R^4$, respectively. This numerical coincidence exhibits itself in the study of Einstein metrics in these dimensions. In $5$-dimensions, we don't have such a nice result. ${\bf R}_5$ is $50$-dimensional. For more information on Einstein metrics, see \cite{Besse}. From this point on, we suppose that $n\geq 3$.

We now consider the second Bianchi identity, and the space of covariant derivatives of abstract curvature tensors. The riemannian curvature $R$ defines a function $\mathcal{R}$ that is defined on $\mathcal{F}_O$ and takes values in ${\bf R}_n$.  The covariant derivative $\nabla R$ is determined by the equivariant, horizontal $1$-form on $\mathcal{F}_O$ given by $D\mathcal{R}=d\mathcal{R}+\rho(\phi)\mathcal{R}$ where $\rho$ denotes the representation of $O(n)$ on ${\bf R}_n$. Again using the $1$-form $\omega$, this is equivalent to a function $\mathcal{DR}$ on $\mathcal{F}_O$ that takes values in the vector space
\begin{eqnarray*}
\Lambda^1\otimes {\bf R}_n\subseteq \Lambda^1\otimes\Lambda^2\odot\Lambda^2.
\end{eqnarray*}
As in the case of the first Bianchi identity, we make the following definition.
\begin{eqnarray*}
\gamma:\Lambda^1\otimes {\bf R}_n&\to &\Lambda^3\otimes \Lambda^2\\
\gamma(A)_{uvwxy}& =&A_{uvwxy}+A_{vwuxy}+A_{wuvxy}.
\end{eqnarray*}
\begin{prop}\label{prop:secondbianchi}
$\gamma$ is given, on elements of $\Lambda^1\otimes\Lambda^2\otimes\Lambda^2$, by
\begin{eqnarray*}
\gamma: \alpha\otimes\sigma\otimes\tau\mapsto \alpha\wedge\sigma\otimes\tau.
\end{eqnarray*}
\end{prop}
The proof is identical to that of the first Bianchi identity. Then, by the second Bianchi identity, the function $\mathcal{DR}$ on $\mathcal{F}_O$ takes values in the linear vector space
\begin{eqnarray*}
{\bf DR}_n=\ker\gamma\subseteq \Lambda^1\otimes{\bf R}_n.
\end{eqnarray*}

\subsection{The decomposition of ${\bf R}_n$ into irreducible factors}

We now consider some elementary linear algebra to elucidate the decomposition of curvature tensors. Suppose that $V$ and $W$  are finite dimensional vector spaces equipped with inner products, and $c:V\to W$ is a surjective linear map. Then $V$ decomposes as $V=V_1+V_2$ where $V_1=\ker c$ and $V_2=V_1^\perp=\text{Im}c^*$ where $c^*:W\to V$ is the adjoint of $c$. Then we can decompose an element $v=v_1+v_2$ where $c(v)=c(v_2)=r\in W$. We wish to have an expression for $v_2$. $cc^*$ is an isomorphism of $W$ so if $r=cc^*(r^\prime)$ then, $v_2=c^*(r^\prime)$. Explicitly, the projection onto the $V_2$-factor in the decomposition of $V$ is given by $v\mapsto c^*(cc^*)^{-1}c(v)$.

\begin{defn}
The \emph{Ricci contraction} is the map
\begin{eqnarray*}
c:{\bf R}_n&\to& Sym^2(\R^n)\\
\langle c(R)v,w\rangle &=&\tr\{x\mapsto R_{xv}w\}\\
&=& \sum_i \langle R_{ve_i}w,e_i\rangle
\end{eqnarray*}
\end{defn}
It can then be seen that $c$ is $O(n)$-equivariant and, slightly less evidently, is surjective onto $Sym^2(\R^n)$. We apply the above linear algebra to ${\bf R}_n$ and the contraction $c$. The inner products on ${\bf R}_n\subseteq Sym^2(\Lambda^2)$ and $Sym^2(\R^n)$ are given by $\tr(\alpha\circ\beta)$ in each case. The inner product on $\Lambda^2$ is given by $-\tr(\phi\circ\psi)$, when the elements are considered skew-symemtric transformations. The adjoint is $c$ will be seen  in the following exercise.

\begin{exer}
Let $S:\R^n\to \R^n$ be a symmetric linear transformation.
\begin{enumerate}
\item Show that if $\phi$ is a skew-symmetric map, then $\{S,\phi\}=S\phi+\phi S:\R^n\to \R^n$ is skew-symmetric.
\item Show that $\{S,\cdot\}:\Lambda^2\to\Lambda^2$ is symmetric, which is to say that $\{S,\cdot\}\in S^2(\Lambda^2)$.
\item Show that $\beta(\{S,\cdot\})=0$, or that $\{S,\cdot\}\in {\bf R}_n$.
 \item Show that if $R\in{\bf R}_n$, then $\langle \{S,\cdot\},R\rangle =\langle S,c(R)\rangle$.
 \end{enumerate}
\end{exer}
Each of these statements are elementary calculations. In the final case the quantity is equal to
\begin{eqnarray*}
\sum_I\langle \{S,\varepsilon_I\},R_{\varepsilon_I}\rangle,
\end{eqnarray*}
summing over an orthonormal basis for $\Lambda^2$. In particular, this shows that the adjoint of the Ricci contraction $c^*(S)\in S^2(\Lambda^2)$ is given by
\begin{eqnarray*}
c^*(S)\phi=\{S,\phi\}.
\end{eqnarray*}
This coincides with the negative of the Kulkarni-Nomizu product $S\owedge g$ of the symmetric bilinear form $\langle S\cdot,\cdot\rangle$ with the metric, as explained in \cite{Besse}, when we consider bilinear forms instead of endomorphisms. The negative sign corresponds to the different definition of the Ricci contraction. Expressed in this way, it can easily be seen that $c^*$ is injective. This in particular implies that $cc^*$ is an isomorphism of $S^2(\R^n)$.
\begin{prop}
The composition $cc^*:S^2(\R^n)\to S^2(\R^n)$ and its inverse $(cc^*)^{-1}$ are given by
\begin{eqnarray*}
cc^*(S)&=& (n-2)S +(\tr S)I\\
(cc^*)^{-1}(S)&=& \frac{1}{n-2}S-\frac{\tr S}{2(n-1)(n-2)}I.
\end{eqnarray*}
where $I$ is the identity map.
\end{prop}
\proof{
\begin{eqnarray*}
\langle (cc^*)(S)v,w\rangle &=& \sum_i\langle (c^*S)_{e_iv}w,e_i\rangle\\
&=& \sum_i\langle S\left(\langle v,w\rangle e_i-\langle e_i,w\rangle v\right),e_i\rangle
 + \sum_i\langle\langle v,Sw\rangle e_i-\langle e_i,Sw\rangle v,e_i\rangle\\
 &=& \langle v,w\rangle \sum_i\langle S_i,e_i\rangle -\sum_i\langle Sv,e_i\rangle\langle e_i,w\rangle\\
 && \ \ \  +\langle Sv,w\rangle\sum\langle e_i,e_i\rangle -\sum_i\langle Sw,e_i\rangle\langle e_i,v\rangle\\
  &=& (n-2) \langle Sv,w\rangle +(\tr S)\langle v,w\rangle.
\end{eqnarray*}
The inverse can be found by supposing that $(cc^*)^{-1}S=\frac{1}{n-2}S+\lambda I$ and solving for $\lambda$.
}
If we define $S^2_0(\R^0)$ to be the symmetric transformations with trace equal to $0$, then any symmetric map $S$ on $\R^n$ can be uniquely expressed as
\begin{eqnarray*}
S=S_0+\frac{\tr S}{n}I
\end{eqnarray*}
where $S_0\in S^2_0(\R^n)$ so $(cc^*)^{-1}$ can be expressed as
\begin{eqnarray*}
(cc^*)^{-1}:S\mapsto \frac{1}{n-2}S_0 +\frac{\tr S}{2n(n-1)}I.
\end{eqnarray*}
We finally apply this to the decomposition of abstract curvature tensors.
\begin{thm}
There exists a unique $O(n)$-invariant decomposition
\begin{eqnarray*}
{\bf R}_n ={\bf S}_n +{\bf Z}_n+ {\bf W}_n
\end{eqnarray*}
where
\begin{eqnarray*}
{\bf W}_n &=& \ker c,\\
{\bf Z}_n &=& \{ c^*(S_0)\ ;\ S_0\in S^2_0(\R^n)\}\\
{\bf S}_n &=& \{\lambda I|_{\Lambda^2}\ ;\ \text{for } \lambda\in \R\}.
\end{eqnarray*}
\end{thm}
\proof{
That the decomposition exists is clear from the above. The detailed proof that each factor is irreducible is given in \cite{BGM}.
}

According to this decomposition, any $R\in {\bf R}_n$ can be expressed as
\begin{eqnarray*}
R=\frac{s}{n(n-1)}I|_{\Lambda^2} +\frac{1}{n-2}\{r_0,\cdot\} +W_R
\end{eqnarray*}
where $r=c(R)\in S^2(\R^n)$ has trace-free part $r_0$ and trace $s=\tr r$, and where $W_R\in {\bf W}_n$. This decomposition is equivariant with respect to the action of $O(n)$ on ${\bf R}_n$.

Let $(M^n,g)$ be a riemannian manifold with orthonormal frame bundle $\mathcal{F}_O$. Consider the curvature as an equivariant function $\mathcal{R}$ on $\mathcal{F}_O$ with values in ${\bf R}_n$, as in Equation \ref{eqn:curvfunction}.
\begin{prop}
The function $c(\mathcal{R})$ on $\mathcal{F}_O$ defines a symmetric $2$-form $r_g$ on $T_M$. $\tr(c(\mathcal{R}))$ defines a function $s_g$ on $M$. The component of $\mathcal{R}$ in ${\bf W}_n$ determines a section of $S^2(\Lambda^2T^*_M)$. These tensors are related to the riemannian curvature operator $R$ by the relation
\begin{eqnarray}\label{eqn:CurvTensorDecomp}
R_{XY} = \frac{s_g}{n(n-1)}X\wedge Y +\frac{1}{n-2}\{r_{g0},X\wedge Y\} +W_g(X,Y)
\end{eqnarray}
\end{prop}
This follows from the equivariance of the functions, which is because $c$ intertwines the representations of $O(n)$ on the respective spaces, and from the relationship between functions on $\mathcal{F}_O$ and tensor fields on $M$ that we have repeatedly used.

\begin{defn}
The function $s_g$ is the \emph{scalar curvature}, the bilinear form $r_g$ is called the \emph{Ricci curvature} and $W_g$ is the \emph{Weyl curvature} of the metric $g$.
\end{defn}

\subsection{Curvature decomposition for reduced holonomy metrics}

In this section we make some elementary observations on the curvature tensor in the case that $H=Hol(A,u)$ is a proper subgroup of $SO(n)$. We will also assume that $M$ is oriented, to immediately reduce the holonomy group from $O(n)$ to $SO(n)$.

To start with, we note that the Levi-Civita connection reduces to the holonomy bundle. This is a principal subbundle $\mathcal{F}_H\subseteq \mathcal{F}_{SO}$ with structure group $H$. In particular, the curvature form $\Omega$ on $\mathcal{F}_H$ takes values in the Lie algebra $\mathfrak{h}=\mathfrak{hol}(A,u)\subseteq \mathfrak{so}(n)\cong \Lambda^2$. Further more, we have the following theorem.
\begin{prop}\label{prop:absttensorsreduced}
Let $\mathcal{R}$ be the curvature tensor, considered as a function on the bundle $\mathcal{F}_H$. Then $\mathcal{R}$ takes values in the subset
\begin{eqnarray*}
{\bf R}_H=\ker \beta\cap S^2(\mathfrak{h})\subseteq {\bf R}_n
\end{eqnarray*}
\end{prop}
This statement just follows from the fact that $\Omega$ takes values in $\mathfrak{h}$, and the curvature function $\mathcal{R}$ takes values in the symmetric product.

Interesting and important properties are of how this space behaves with respect to the decomposition of the curvature tensor according to Equation \ref{eqn:CurvTensorDecomp} in the case that the manifold has holonomy in $H$. For some possible holonomy groups, the Ricci curvature must vanish identitically. For others, the trace-free part vanishes (which is to say Einstein), while for others little can be said. We will discuss this later.

In the classification theorem that we will describe, it is also necessary to consider the second Bianchi identity. We consider the restriction
\begin{eqnarray*}
\gamma:\Lambda^1\otimes \mathfrak{h}\odot\mathfrak{h}&\to &\Lambda^3\otimes \mathfrak{h}\\
\alpha\otimes \sigma\odot\tau &\mapsto &\frac{1}{2}\left((\alpha\wedge\sigma)\otimes\tau +(\alpha\wedge\tau)\otimes\sigma\right).
\end{eqnarray*}

If $(M,g)$ is a riemannian manifold and the holonomy group $H=Hol(g,u)$ is a proper subgroup of $SO(n)$ then the covariant derivative of the curvature tensor $\mathcal{DR}$ takes values in the subset  $\ker\gamma\cap\Lambda^1\otimes{\bf R}_H$. The important observation of Berger about this space is that for almost all subgroups of $SO(n)$, the space $\ker\gamma$ is empty. Which implies that $\nabla R=0$ or, as we will see in the following section, $(M,g)$ is locally symmetric.


\section{Symmetric Spaces}\label{sec:SymmSpaces}


On a neighbourhood of a point, say $x$, in any manifold it is possible to define a diffeomorphism that fixes $x$, and has minus the identity as its derivative at $x$. In the presence of a riemannian metric, such a map can be given by reversing the geodesics through $x$. In this section we study a class of riemannian manifolds for which there exists such a map for each point $x$ that is also a local isometry. The study of riemannian symmetric spaces is undertaken in great detail in a range of classical references, including \cite{Chev,Helg,KN}. \cite{Salamon} has a shorter summary of the theory, but one that is particularly relevant for the study of holonomy groups.

\begin{defn}
The riemannian manifold $(M,g)$ is a \emph{riemannian symmetric space} if for every $x\in M$ there exists an isometry $s_x:M\to M$ such that
\begin{eqnarray*}
s_x(x)&=&x\\
s_{x*}&=& -\text{id}_x.
\end{eqnarray*}
Necessarily, the square of this map must be the identity : $s_x^2=\text{id}$. $(M,g)$ is a locally symmetric space if for any point, such a isometry can be defined on a neighbourhood of the point.
\end{defn}
We can summarise some immediate consequences of this definition.
\begin{prop} Let $(M,g)$ be a riemannian symmetric space.
\begin{enumerate}
\item The curvature tensor satisfies $\nabla R=0$.
\item The metric $g$ is complete.
\item The group of isometries of $g$ acts transitively on $M$.
\end{enumerate}
\end{prop}
\proof{
At a point $x\in M$ we consider the isometry $s_x$. As an isometry, this must preserve the curvature tensor and its covariant derivative. At $x$ however, $s_{x*}=-1$ so
\begin{eqnarray*}
(\nabla R)(X,Y,Z,W)&=&s_x^*(\nabla R)(X,Y,Z,W)\\
&=&s_{x*}\left(\nabla_{s_{x*}X}R\right)_{s_{x*}Ys_{x*}Z}(s_{x*}W)\\
&=& (-1)^5(\nabla R)(X,Y,Z,W)
\end{eqnarray*}
which implies that $\nabla R=0$.

For $x\in M$, let $y=\exp_x(\varepsilon X)$ be a point on a geodesic $\gamma$ emanating from $x$. Then,
\begin{eqnarray*}
\gamma(t)=
\left\{
\begin{array}{ll}
  \exp_x(tX) & \quad t\in [0,\varepsilon]\\
    s_y\left(\exp_x((\varepsilon-t)X)\right) & \quad t\in [\varepsilon,2\varepsilon].
\end{array}\right.
\end{eqnarray*}
That is, $\gamma$ can be extended by reflection through $y$. This means that every geodesic is infinitely extensible and so $g$ is complete.

Finally, for $x,y\in M$, by completeness there exists a geodesic $\gamma$ between them such that $\gamma(0)=x$ and $\gamma(2T)=y$. Then, consider the symmetry $s_z$ where $z=\exp_x(T)$. This map preserves the geodesic and  distances along it so $s_z(x)=y$.
}

This means that $M$ can be expressed as $M=G/H$, where $G$ is a connected group of isometries of $M$ and $H=G_o=\{g\in G\ ;\ g\cdot o=o\}$ is the isotropy group of a fixed point $o\in M$. In particular, $H$ is closed and $G$ can be considered a principal $H$-bundle over $M$.

Now mention thing about $H$ nec being contained in $SO(T_oM)$.

 For much of the remainder of this section, we will consider the algebra associated to $G$ and $H$. To start, we observe that we can define an automorphism of $G$.
\begin{eqnarray}\label{eqn:symspacesymmetry}
\sigma: G &\to& G\nonumber\\
\sigma(g)&=& s_o\circ g\circ s_o^{-1}:M\to M.
\end{eqnarray}
Then, we can easily see that $\sigma^2=\text{Id}$. Furthermore, if $h\in H$, $\sigma(h)(o)=o$ and $\sigma(h)_*=s_{o*}h_*s_{o*}^{-1}=h_*$ at $o$. This implies (see Exercise \ref{exer:IsomCoinciding}) that $\sigma(h)=h$ and in particular $H$ is contained in the fixed point set $G^\sigma=\{g\in G\ ;\ \sigma(g)=g\}$ of the automorphism $\sigma$.

We can almost say that these subgroups do in fact coincide.
\begin{prop}
Let $G$ be the connected component of the identity of the isometry group of a riemannian symmetric space $M=G/H$ where $H$ is the subgroup that fixes a chosen point $o$. Let $\sigma$ be as defined in Equation \ref{eqn:symspacesymmetry}. Then $H$ is contained in $G^\sigma$, and contains the connected component of the identity of $G^\sigma$.
\end{prop}
\proof{We have already shown the first of these statements. Let $g_t$ be a one-parameter subgroup of $G^\sigma$ and define $o_t=g_t\cdot o\in M$. Then, $g_t\in G^\sigma$ is equivalent to saying that $s_o\circ g_t=g_t\circ s_o$. Therefore,
\begin{eqnarray*}
s_o(o_t)=g_ts_o(o)=o_t
\end{eqnarray*}
and the elements $g_t$ permute the fixed points of the symmetry map $s_o$. However, since $s_{o*}=-\text{id}$, these fixed points are isolated from one another, and necessarily $o_t=o$ and $g_t\in H$. This implies that $G^\sigma_e\subseteq H$.
}

\begin{exer}\label{exer:IsomCoinciding}
Let $(M,g)$ be a connected riemannian manifold. Use an open/closedness argument to show that if $g$ and $h$ are isometries of $M$ such that $g(p)=h(p)$ and $g_{*p}=h_{*p}$ then $g=h$.
\end{exer}

The fact that a symmetric space $M$ is homogeneous $M=G/H$ reminds us of some examples that we encountered in the first section of these notes.  The above result in particular means that $H$ and $G^\sigma$ share the same Lie algebra $\mathfrak{h}$.
 We consider the linearisation of $\sigma$. Specifically, the automorphism $\sigma:G\to G$ induces an automorphism $\sigma:\mathfrak{g}\to \mathfrak{g}$ that satisfies $\sigma^2=\text{Id}$.

\begin{prop}
The linearisation $\sigma:\mathfrak{g}\to \mathfrak{g}$ induces a decomposition $\mathfrak{g}=\mathfrak{h}+\mathfrak{m}$ that is related to the Lie bracket on $\mathfrak{g}$ by
\begin{eqnarray*}
[\mathfrak{h},\mathfrak{h}]\subseteq \mathfrak{h},\ \ \ [\mathfrak{h},\mathfrak{m}]\subseteq\mathfrak{m},\ \ \ [\mathfrak{m},\mathfrak{m}]\subseteq \mathfrak{h}.
\end{eqnarray*}
\end{prop}
\proof{
Since $\sigma^2=\text{Id}$, we take $\mathfrak{h}$ and $\mathfrak{m}$ to be the $+1$ and $-1$ eigenspaces of $\sigma$. Then the relations in the proposition then follow from the fact that $(-1)^2=1$, for instance.
It can furthermore be seen that $\mathfrak{h}$ is the Lie algebra of the group $H$.
}

In particular, a riemannian symmetric space is a reductive homogeneous space, as we have considered in the first chapter of these notes. Additionally however we have the condition that $[\mathfrak{m},\mathfrak{m}]\subseteq \mathfrak{h}$, which will help us greatly. We will from this point largely study the riemannian geometry of $M=G/H$ by considering the algebra of this decomposition.
If we consider $\mathfrak{g}$ to be the set of left-invariant vector fields on $G$, then the splitting $\mathfrak{g}=\mathfrak{h}+\mathfrak{m}$, defines two distributions on $G$. The distribution $\mathfrak{h}$ is the kernel of the tangent projection $TG\to TM$ and $\pi:\mathfrak{m}\to T_{\pi(g)}M$ is at each point an isomorphism.

We have seen that $\mathcal{F}_{GL}$ is the set of isomorphisms $u:\R^n\to T_xM$ for all $x\in M$. If we identify $\mathfrak{m}\cong \R^n$ then we can identify $G$ with a subbundle of $\mathcal{F}_{GL}$ with structure group $H$, or as a subset of isomorphisms $\pi_g:\mathfrak{m}\to T_xM$ where $x=\pi(g)\in M$.
\begin{prop}
The projection $\pi$ which at each point $g$ defines the isomorphism
\begin{eqnarray*}
\pi_g:\mathfrak{m}\to T_xM
\end{eqnarray*}
defines $G$ as a subbundle of $\mathcal{F}_{GL}$ with structure group $H$, where $H$ is considered as a subgroup of $GL(\mathfrak{m})$.
\end{prop}
\proof{
This follows from Proposition \ref{prop:critforsubbundle}. We can check each of the conditions for this result. The only one that is not evident is the second one. $\pi_g$ is defined by $\pi_g(v)=g_*(\pi v)$ where $\pi:\mathfrak{m}\to T_oM$, and so one can easily see that
\begin{eqnarray*}
\pi_{gh}(v)=\pi_g(\text{Ad}_hv)
\end{eqnarray*}
from which the second and third conditions follow. The local section comes from the local triviality of the fibration $G\to M$.
}


In the first  chapter of these notes we saw that there are two important forms on the principal frame bundle to a manifold. The first is the $\R^n$-valued form $\omega$ given by $\omega(X)=u^{-1}(\pi X)$. The second is the $\mathfrak{h}$-valued $1$-form given as a connection form. In the following result result we see that these can be calculated easily on a symmetric space.

We will collect some elemenatary facts about the geometry of $M$.
\begin{thm}\label{thm:symspaceproperties}
Let $M=G/H$ be a riemannian symmetric space of real dimension $n$, where $G$ is a connected group of isometries. Then,
\begin{enumerate}
\item Then tangent bundle of $M$ can be expressed as the associated bundle $T_M=(G\times \mathfrak{m})/\sim$, arising from the adjoint representation of $H$ on $\mathfrak{m}$.
\item The metric on $M$ is induced by an $H$-invariant inner product on $\mathfrak{m}$.
\item The fundamental $\R^n=\mathfrak{m}$ valued $1$-form on $G$ is given by
\begin{eqnarray*}
\omega :TG &\to &\R^n\cong \mathfrak{m}\\
\omega(X)&=& \pi^{\mathfrak{m}}(g^{-1}_*X)\ \ \ \text{for }X\in T_gG,
\end{eqnarray*}
and where $\pi^{\mathfrak{m}}:\mathfrak{h}+\mathfrak{m}\to \mathfrak{m}$ is the projection.
\item The Levi-Civita connection form $\phi$ is given by projection to the $\mathfrak{h}$-component in the same decomposition. That is,
\begin{eqnarray}\label{eqn:LConSymSpace}
d\omega =-\rho(\phi)\wedge\omega
\end{eqnarray}
where $\rho$ is the adjoint representation of $H$ on $\mathfrak{m}$.
\item The curvature, considered as an $\mathfrak{h}$-valued $2$-form on $G$ is given by
\begin{eqnarray*}
\Omega(X,Y)=-[X,Y]
\end{eqnarray*}
where $X$ and $Y$ are left invariant elements of $\mathfrak{m}$.
\end{enumerate}
\end{thm}
\proof{
The associated bundle construction, to obtain a vector bundle with typical fibre $\R^n$ from a principal $G$-bundle via a representation of $G$ on $\R^n$, is given in Proposition \ref{prop:assocbundle}. In this case, the isomorphism
\begin{eqnarray*}
(G\times \mathfrak{m})/\sim\  \to TM
\end{eqnarray*}
is given in Proposition \ref{prop:tangbunfofhomspace} for reductive homogeneous spaces, as in this case.

By assumption, $G$ acts on $M$ by isometries so we can apply Proposition \ref{prop:invariantmetrics}. This result states that for a reductive homogeneous space $M=G/H$ with decomposition $\mathfrak{g}=\mathfrak{h}+\mathfrak{m}$, there exists a correspondence between $G$-invariant metrics on $M$ and $\ad H$-invariant inner products on $\mathfrak{m}$.

That the formula given here gives the fundamental $\R^n\cong\mathfrak{m}$ valued $1$-form on $G$ can be seen from the above embedding of $G$ into $\mathcal{F}_{GL}$. A vector is translated back to $e\in G$ and then projected to $\mathfrak{m}$, using the inverse of the map $\pi_g$.

We show that Equation \ref{eqn:LConSymSpace} holds. Suppose that $X,Y$ are vectors tangent to $G$. We can suppose that they are both left-invariant. If either both $X$ and $Y$ lie in $\mathfrak{h}$, or both lie in $\mathfrak{m}$, then both sides of Equation \ref{eqn:LConSymSpace} are equal to zero. This in particular uses the fact that $\omega([X,Y])=0$, or that $[X,Y]\in \mathfrak{h}$ in these cases. If $X\in \mathfrak{h}$ and $Y\in \mathfrak{m}$, then
\begin{eqnarray*}
d\omega(X,Y)&=& -\omega([X,Y])\\
&=& -[X,Y]\\
&=& -\rho(\phi(X))\omega(Y).
\end{eqnarray*}
In the case of the curvature, if $X,Y\in\mathfrak{m}$ are left-invariant,
\begin{eqnarray*}
\Omega(X,Y)=-\phi([X,Y])=-[X,Y].
\end{eqnarray*}
}

We can make some comments about these results. From the second statement we see that, via the adjoint representation, $H\subseteq SO(\mathfrak{m})$ for some positive definite inner product on $\mathfrak{m}$. This fact is also related to the inclusion of $H$ in $SO(T_oM)$. Since $h\cdot o=o$, $h$ must act, infinitesimally, on $T_oM$. From \ref{exer:IsomCoinciding}, this inclusion must be injective. Via the inclusion of $G$ in $\mathcal{F}_{GL}$, we can assume that $G\subseteq\mathcal{F}_O$ and any connection of $G$ takes values in the Lie algebra $\mathfrak{h}\subseteq\mathfrak{so(m)}$. In particular, the connection $\phi$ that is shown to be torsion-free according to Equation \ref{eqn:LConSymSpace}.
\begin{exer}
Show that the connection $\phi$ given above on the bundle $G\to M$ satisfies $\nabla R=0$.
\end{exer}
\begin{example} {\bf The complex grassmannian}
 We now consider a particular example in more detail.
$G(2,C)$ is the space of $2$-dimensional complex linear subspaces in $\C^n$. It can easily be seen that the group $U(n)$ acts transitively, by considering hermitian bases for $\C^n$. The subgroup $U(1)$ acts trivially so we consider \begin{eqnarray*}
G=U(n)/U(1),
\end{eqnarray*}
 which acts effectively of $G(2,\C^n)$. The isotropy group of $\xi=\text{span}\{e^1,e^2\}$ fixes $\xi$ and the orthogonal complement $\xi^\perp$ and equals
 \begin{eqnarray*}
 H=U(2)U(n-2)=(U(2)\times U(n-2))/U(1)
 \end{eqnarray*}
 where again we quotient by the subgroup that acts trivially. The symmetries can be understood at the Lie algebra level. The algebra $\mathfrak{u}(n)$ decomposes as elements of the form
 \begin{eqnarray*}
\left(\begin{array}{cc}
A & X  \\
-\bar{X}^t & B  \\
  \end{array}\right)
  \end{eqnarray*}
 for $A\in \mathfrak{u}(2)$, $B\in \mathfrak{u}(n-2)$ and $X\in M_\C(2,n-2)$. The subalgebra $\mathfrak{h}^\prime=\mathfrak{u}(2)+\mathfrak{u}(n-2)$ consists of elements with $X=0$ and the subspace $\mathfrak{m}^\prime$ is of elements with $A=B=0$. We define $\mathfrak{h}=(\mathfrak{u}(2)+\mathfrak{u}(n-2))/\mathfrak{u}(1)$ and $\mathfrak{m}$ to be the images of these spaces in the quotient $\mathfrak{u}(n)/\mathfrak{u}(1)$. We can define an automorphism $\sigma:\mathfrak{u}(n)\to \mathfrak{u}(n)$ by
 \begin{eqnarray*}
 \sigma &\equiv & 1\ \ \text{on }\mathfrak{h}^\prime\\
 \sigma &\equiv & -1\ \ \text{on } \mathfrak{m}^\prime.
 \end{eqnarray*}
This automorphism descends to the algebra $\mathfrak{u}(n)/\mathfrak{u}(1)$ and defines an automorphism of $U(n)/U(1)$. In particular, since $H\subseteq \{g\ ;\ \sigma(g)=g\}$, $\sigma$ descends to define a difeomorphism of $G(2,\C^n)=G/H$.
Then, for any point $x=gH\in G/H$ we can define a involution $s_x:G(2,\C^n)\to G(2,\C^n)$ that fixes $x$ and acts as $-1$ on $T_xG(2,\C^n)$. We act on cosets by
\begin{eqnarray*}
s_x&=&g\circ\sigma\circ g^{-1}\\
s_x(g^\prime H)&=& g\sigma(g^{-1}g^\prime)H.
\end{eqnarray*}
Then, as we have seen in the more general case, any inner product on $\mathfrak{m}$ that is preserved by the action of $U(2)U(n-2)$ extends to define a metric on $G(2,\C^n)$ that is invariant under the action of $U(n)/U(1)$ and for which each of the maps $s_x$ are isometries.
\end{example}

This example illustrates the method with which we can recover the riemannian symmetric space $(M,g)$ from purely algebraic data. The proof of the following is found in \cite{KN}.
\begin{thm}
Let $G$ be a Lie group and $H\subseteq G$ is a closed Lie subgroup such that there exists a decomposition of Lie algebra $\mathfrak{g}=\mathfrak{h}+\mathfrak{m}$ such that
\begin{eqnarray*}
[\mathfrak{h},\mathfrak{h}]\subseteq \mathfrak{h},\ \ \ [\mathfrak{h},\mathfrak{m}]\subseteq \mathfrak{m},\ \ \ [\mathfrak{m},\mathfrak{m}]\subseteq \mathfrak{h}.
\end{eqnarray*}
Suppose also that $\text{Ad}H$ is a compact subgroup of $GL(\mathfrak{m})$. Then $M=G/H$ admits a $G$-invariant metric that gives it the structure of a riemannian symmetric space.
\end{thm}

We now return to study the holonomy groups of symmetric spaces. From Theorem \ref{thm:symspaceproperties} we can see that the Levi-Civita connection on $M$ is defined on the subbundle given by $G\mathcal{F}_{GL}$. The holonomy bundle $P$ is therefore a subbundle of $G$ and $Hol(g)$ is a subgroup of $H$. By the Ambrose-Singer theorem (Theorem \ref{thm:AmbroseSinger}) the holonomy Lie algebra is spanned by elements $\Omega_v(X,Y)$ as $v$ varies in $P$. In the symmetric space case we see that if $X,Y$ are left-invariant vector fields in $\mathfrak{m}$, then $[X,Y]$ is left-invariant and so the value of $\Omega(X,Y)=-[X,Y]\in\mathfrak{h}$ is independent of $g\in G$ (Note that this follows also from the fact that $\nabla R=0$). It suffices to evaluate the curvature at $e\in G$. We have the following result.
\begin{prop}
Let $(M,g)$ be a riemannian symmetric space where $M=G/H$ as above. Then the holonomy Lie algebra is equal to the vector space generated by $[\mathfrak{m},\mathfrak{m}]\subseteq \mathfrak{h}$.
\end{prop}

We can use this result to further understand the previous example, although we will now consider the group $G=SU(n)$ acting transitively on $G(2,\C^n)$ with isotropy $H=S(U(2)\times U(n-2))$. This is to say, the subgroup of the product $U(2)\times U(n-2)\subseteq SU(n)$ with unit determinant. Then $\mathfrak{su}(n)=\mathfrak{h}+\mathfrak{m}$ where $\mathfrak{h}=\mathfrak{s}(\mathfrak{u}(2)+\mathfrak{u}(n-2))$. Then, since $SU(n)$ is compact and semi-simple, one can see that $[\mathfrak{m},\mathfrak{m}]=\mathfrak{h}$.
\begin{exer}
Use the fact that if $G$ is compact and semi-simple then $\mathfrak{g}$ admits a bi-invariant inner product, together with the fact that the representation of $\mathfrak{h}$ on $\mathfrak{m}$ is effective, to show that if $M=G/H$ is symmetric and $G$ is compact and semi-simple, then $\mathfrak{h}=[\mathfrak{m},\mathfrak{m}]$.
\end{exer}

This example, though somewhat non-trivial, is useful because it describes the properties of a large family of riemannian symmetric spaces. We will not describe the classification in more detail, but it is was achieved by Cartan in the 1920's (see \cite{Cartan, Helg, KN}). In particular, the holonomy groups for this class of riemannian manifolds can be fully understood. If $G$ is a simple Lie group of compact type, there are $7$ infinite families of riemannian symmetric spaces and $12$ exceptional examples. The holonomy groups of the non-exceptional cases consist of the groups $SO(n)$, $Sp(n)$ and  $U(n)$, together with certain products of these (as in the case of $G(2,\C^n)$). These manifolds can be interpreted as grassmannians of particular types. The exceptional examples are more complicated, and arise from representations of the exceptional Lie algebras. Many of these can be seen to parametrize certain geometric objects.

\section{Classification Theorem of  Berger}
\label{ss.Berger}
At this point, it seems natural to ask which subgroups of $O(n)$ can be realized as the
holonomy group of an $n$-dimensional Riemannian manifold $(M, g)$.

Before we proceed with the answer, it is helpful to further precise this question. By passing
to the study of the restricted holonomy $Hol_0(g)$, we may
assume that $M$ is simply connected.
Also, we have seen in Section \ref{sec:reducible} that if $g$ is locally reducible,  then $Hol_0(g)$ also splits as a product.
So, we will consider irreducible metrics. Finally, if $g$ is locally symmetric, then the restricted holonomy
falls in the list of holonomy groups of symmetric spaces.

Therefore, our reformulated question becomes: which subgroups of $SO(n)$ can be realized as
the holonomy group of  an irreducible, non-symmetric metric $g$ on a simply connected
Riemannian manifold $M$?

This question was addressed by Berger \cite{Berger}.

\begin{thm}
\label{thm:Berger} Let $M$ be an $n$-dimensional 
simply-connected riemannian manifold, with an irreducible, non-symmetric Riemannian metric
on $M$. Then, exactly one of the seven cases hold:
\been
\item
$Hol(g) = SO(n)$;
\item
$n= 2m$, $m >1$, and $Hol(g) = U(m) \subseteq SO(2m)$. The metric $g$ is called a {\em \kahler metric};
\item
$n= 2m$, $m >1$, and $Hol(g) = SU(m) \subseteq SO(2m)$. The metric $g$ is a {\em Calabi-Yau metric};
\item
$n= 4m$, $m >1$, and $Hol(g) = Sp(m) \subseteq SO(4m)$, and $g$ is a {\em \hk metric};
\item $n= 4m$, $m >1$, and $Hol(g) = Sp(m)Sp(1) \subseteq SO(4m)$. The metric $g$ is called a {\em quaternionic \kahler metric};
\item
$n=7$, and $Hol(g) = G_2 \subseteq SO(7)$;
\item
$n= 8$,  and $Hol(g) = Spin(7) \subseteq SO(8)$;
\eeen
\end{thm}

In fact, Berger showed that those groups were the only possibilities for 
the holonomy group of a Riemannian manifold.
It took another thirty years to confirm that all of the groups in Berger's list
do occur.

Berger's original list also included the possibility of $Spin(9)$ as a subgroup of $SO(16)$. 
Riemannian manifolds with such holonomy were later shown independently by D. Alekseevski \cite{Alexseevski} and Brown-Gray 
\cite{Brown-Gray} to be necessarily locally symmetric, and hence should be excluded from the list above.

In \cite{J}, Joyce describes the groups in Berger's list making a nice analogy with 
the four division algebras: $\RR$, $\CC$, $\HH$ and $\mO$. The group $O(n)$ is the automorphism 
group of $\RR^n$, 
$U(n)$ is the automorphism
group of $\CC^n$, $Sp(n)$ is the automorphism
group of $\HH^n$, and the groups $SO(n)$, $SU(n)$ and $Sp(n) Sp(1)$ 
are the "determinant one" subgroups of $O(n)$, $U(n)$ and $Sp(n)$, respectively.

As for the exceptional groups, we can still make an analogy to the octonions, even though
not so directly. If we write $\mO = \RR \oplus \Im \mO$, we will see in Section \ref{s.G2}
that $G_2$ is the automorphism group of $\Im \mO$, while the group $Spin(7)$ is the automorphism
group of $\RR^8$ which preserve a certain part of the multiplicative structure of $\mO$ (see Section \ref{s.Spin7}).

From Berger's classification, we can split the  holonomy groups of riemannian metrics 
into two main categories: \kahler  and non-\kahler holonomy groups.
The groups $U(n)$, $SU(n)$ and $Sp(n)$ represent \kahler metrics on $M$, so their 
inclusion in the \kahler group is clear. Even though metrics  with holonomy group contained
in $Sp(n) Sp(1)$ 
are not necessarily \kahler, we will see in Section \ref{s.SpnSp1} that a manifold
$(M, g)$ with this holonomy is such that there exists a $\CC\PP^1$-bundle $Z$ over $M$,
called the {\em twistor space} of $M$, such that $Z$
 is a complex manifold, and all of the geometric structure of $(M, g)$ can be uniquely recovered
 from holomorphic data in $Z$. If $(M, g)$ has positive scalar curvature, $Z$ is in fact \kahler.
 Because of this, metrics with holonomy group contained
in $Sp(n) Sp(1)$  will also be included in the \kahler category, as it 
can be understood with tools from complex and \kahler geometry.
The remaining interesting groups are $G_2$ and $Spin(7)$, the exceptional holonomy groups, 
that will form a class on their own.

Here is an outline of the main ideas in the proof of Theorem \ref{thm:Berger}. First, 
note that since $M$ is assumed to be simply-connected, Theorem \ref{thm.Hol0iscpt}
 guarantees that $Hol(g) \subset SO(n)$ is a closed, connected subgroup. Since
 $g$
 is irreducible, that implies that the representation of $Hol(g)$ in $\RR^n$ is also irreducible.
 Berger's proof now consists of testing, for each closed, connected 
subgroup $H \subset SO(n)$ with irreducible representation in $\RR^n$, 
if such $H$ could be a holonomy group.

For $H$ to be a holonomy group, it needs to pass two tests: Theorem \ref{prop:absttensorsreduced}
tells us that if ${\mathfrak{h} = \text{Lie} H}$, then the riemann curvature 
tensor belongs to ${\bf R}_H=S^2 \mathfrak{h}\cap \ker(\beta)$, and the first Bianchi identity holds.
So, in one hand, if  ${\mathfrak{h}} $ is small, the subspace ${\bf R}_H \subset \mathfrak{h}$ will be  
even smaller. On the other hand, Ambrose-Singer Theorem (Theorem \ref{thm:AmbroseSinger}) tells us that
${\bf R}_H $ needs to be big enough to generate $\mathfrak{h}$. Many subgroups $H$ of $SO(n)$ will 
fail this test.  The second test is to satisfy the second Bianchi identity (see Proposition \ref{prop:secondbianchi}), without having the curvature tensor
be covariant constant, as this last condition would mean that $g$ is locally symmetric.

\section{Holonomy groups and cohomology}
\label{s.stilltobemade}

Let $(M, g)$ be a compact, riemannian manifold. 
The goal of this section is to prove that a $G$-structure on $M$
divides $\La^k(T^*M)$ into a  sum of vector bundles, which correspond
to irreducible representations of $G$.

Let $G \subset GL(n,\RR)$ be  a subgroup, and let $Q$ be a $G$-structure
on $M$. Then, to each representation $\rho$ of $G$ on a vector space $V$,
we can associate a vector bundle $\rho(Q)$ over $M$, with fiber $V$. The important observation
here is that if $\rho $ is the restriction to $G$ of the representation of $GL(n, \RR)$ in 
$\RR^n$, then $\rho(Q) = TM$.

This representation induces representations of $G$ on the dual $V^*$, and its
exterior powers $\La^kV^*$. We write $\rho^k$ for the representation of 
$G$ on $\La^kV^*$. Then, due to the observation on the above paragraph, 
$\rho(Q) = \La^k T^*M$.

Now, if $G$ is a proper subgroup of $GL(n, \RR)$, $\rho^k$ of $G$ in $\La^kV^*$ may be reducible.
Then, we can split $\La^k V^*$ and $\rho^k$ into a sum of irreducible representations of $G$.
The following proposition 
summarizes the conclusions we can derive from the above, and that will be used widely in the 
next chapter.

\begin{prop}
\label{p.redexteriorform}
Let $G \subset GL(n,\RR)$ be  a Lie subgroup, and let $M$ be an n-manifold
with $G$-structure $Q$. Then, there is a natural decomposition 
$$
\La^k T^*M = \sum_{i \in I_k} \La^k_i,
$$
where $\La^k_i$ is a vector subbundle of $\La^k T*M$  with fiber $W_i^k$, where
the representation 
$(\rho^k, \La^k T^*M) $ splits into a direct sum $  \sum_{i \in I_k} \La^k_i$ of 
irreducible representations of $G$ on $\La^k T^*M$.

Moreover, if $\phi: \La^k V^* \rar \La^\ell V^*$ is a $G$-equivariant linear map,
then there is a corresponding map $\Phi: \La^k T^* M \rar  \La^\ell T^* M$.
\end{prop}

\begin{example}
\label{ex.HodgeStar}
The {\em Hodge Star Operator $*$} is defined as follows. For  a 
riemannian manifold $(M,g)$, and for any
any $\beta \in \La^kT^*M$, the hodge dual $*\beta$ is the unique
$(n-k)$-form such that 
$$
\ip{\alpha, \beta}_{L^2}  d\Vol_g = \alpha \wedge * \beta, \, \, \, {\text{for all }} \,   \beta \in \La^kT^*M,
$$
where $\ip{\alpha, \beta}_{L^2} $
is the $L^2$-inner product on $(M, g)$.

We want to relate the Hodge star with Proposition \ref{p.redexteriorform}.
Let $G = SO(n)$. Then, a $G$-structure on an $n$-manifold $M$ is equivalent to 
a riemannian metric $g$ and an orientation. Then, there is an isomorphism 
$$
* : \La^k V^* \rar \La^{n-k} V^*,
$$
defined in such a way that, for  $\om \in \La^k V^*$, 
$\om \wedge * \om = d{\text{vol}}_M$.

Thanks to Proposition \ref{p.redexteriorform}, this map
induces an isomorphism 
$$
* : \La^k T^*M \rar \La^{n-k} T^*M,
$$
which is the { Hodge Star operator}.
\end{example}

The $d$-Laplacian can be defined as 
\beeq
\label{e.Laplacian}
\De = d d^* + d^*d,
\eeeq
where the map $d^*: C^\infty(\La^k T^*M) \rar C^\infty(\La^{k-1} T^*M)$,
the {\em formal adjoint of $d$}\footnote{The name 'formal adjoint' refers to the fact that $d^*$ is indeed the adjoint of 
$d$ with respect to the $L^2$-inner product on $M$}, is defined by 
$$d^* \beta = (-1)^{kn + n + 1} * d (*\beta).
$$
for $\beta \in C^\infty(\La^k T^*M)$.
Note that the definition of the Laplacian depends on the metric $g$, as the Hodge star 
also depends on it.

Now, we summarize the main results to be used in the following chapter.
\begin{thm}
Let $M$ be a compact $n$-manifold with a torsion-free $G$-structure $Q$ on $M$,
where $G$ is a Lie subgroup of $O(n)$, and let $g$ be the riemannian metric associated
to $Q$. Let $\La^k_i$ be defined as in Proposition~\ref{p.redexteriorform}.

The Laplacian $\De = d d^* + d^*d$ of $g$ maps $C^\infty(\La^k_i)$ to itself.
Define $\cH_i^k = \Ker \De|_{\La^k_i}$, and let $H_i^k(M, \RR) \subset H^k(M, \RR)$
be the subspace of $H^k(M, \RR)$ with representatives in $\cH^k_i$. 

Then, $H_i^k(M, \RR)  \cong \cH^k_i$, and we have the direct sum decomposition
\beeq
H^k(M, \RR) = \oplus_{i \in I_k} H_i^k(M, \RR)
\eeeq
induced by the irreducible representation decomposition of the representation $G$ in $\RR^n.$
\end{thm}

%% file: IrreducibleRiemannianGroups.tex
\chapter{Irreducible Riemannian Groups}
\label{ch.irred}
This chapter is devoted to the study of each one of the groups on Berger's Classification.

Since $Hol(g) \subset SO(n)$ for any Riemannian metric $g$, we will skip the situation when
all we know is that the holonomy is contained in $SO(n)$, since this is not interesting.

\section{$U(n)$}
\label{s.Un}
According to Berger's classification, the riemannian metrics
$g$ on even-dimensional manifolds with $\Hol(g) = U(n)$
are the \kahler metrics. Such condition imposes serious restrictions
on the curvature tensor and in the topology of the manifold.
We now define such manifolds, and list a few of their properties.
The treatment follows closely \cite{Santoro}.

\begin{defn}
Let $M$ be an $n$-dimensional manifold. We say that
$M$ is a {\em complex manifold} if it admits a system of
{\em holomorphic coordinate charts}, that is, charts such that the
transition functions are biholomorphisms.
\end{defn}

 At each point $p$ in a complex manifold $M$, we can define a map
 $$
 J: T_p M \rightarrow T_pM, \hspace{1cm} J = d(z^{-1}\circ \sqrt{-1}\circ z),
 $$
where $z$ is a holomorphic coordinate defined on a neighborhood of $p$.

It is simple to check that $J^2 = -Id$, and we call $J$ an {\em almost-complex structure}.
In fact, the definition of an almost-complex structure is more general: it
is simply a map $J \in End(TM)$ such that $J^2 = -Id$.
We see that a complex structure induced an almost-complex structure, but the converse is
not true. We say that the almost-complex structure $J$ is {\em integrable} if there
exists an underlying complex structure which generates it.

The {\em Newlander-Nirenberg Theorem} states that an almost-complex structure is
integrable if the {\em Nijenhuis tensor}
$$
N_J(X, Y) = [JX, JY] - J[JX, JY]- J[X, JY] - [X,Y]
$$
vanishes identically.

Let $M$ be a compact, complex manifold of complex dimension
$n$, and consider $g$,
a hermitian metric defined on $M$.
Note that $g$ is a complex-valued sesquilinear form acting on
$TM \times TM$, and can therefore be written as
$$
g = S - 2\sqrt{-1} \om_g,
$$
where $S$ and $-\om$ are real bilinear forms.

If $(z_1, \dots, z_n)$ are local coordinates around a point $p\in M$,
we can write the metric $g$ as $\sum g_{i\jbar}dz^i \otimes d\zbar^j$.
Then, it is easy to see that in these coordinates
$$
\om_g = \frac{\sqrt{-1}}{2} \sum_{i, j = 1}^n g_{i\jbar}dz^i \wedge d\zbar^j.
$$
The form $\om_g$ is a real $2-$form of type $(1,1)$, and is called
the {\it fundamental form} of the metric $g$.

\begin{defn}
We say that a hermitian metric on a complex manifold is
{\em \kahler} if its associated fundamental form $\om_g$
is closed, {\em i.e.}, $d \om_g = 0$.
A complex manifold equipped with a \kahler metric is
called a {\em \kahler manifold}.
\end{defn}

Another characterization of a \kahler manifold $M$ is a manifold
equipped with an almost-complex structure $J$ and a metric $g$ such that
$g$ is $J$-invariant ({\em i.e.}, $g(JX,JY) = g(X, Y)$) and
$J$ is parallel with respect to the Levi-Civita connection of $g$.

The reader can check that the conditions on $J$ and $g$
of a \kahler manifold implies that the Nijenhuis tensor vanishes, so
any \kahler manifold is necessarily complex.

On a \kahler manifold \footnote{In fact, this whole paragraph holds for an
almost-complex manifold.} $M$, the complexified tangent bundle $TM_\CC = TM \otimes \CC$
has a natural splitting. If $p$ is a point in $M$, the extension of map $J$ to
$T_pM_\CC$ (as a complex-linear map) has $\sqrt{-1}$ and $-\sqrt{-1}$ as
eigenvalues, and we define the associated eigenspaces as $T^{1,0}_pM$
and $T^{0,1}_pM$. Hence, we have a decomposition
$$
TM_\CC = T^{1, 0} M \oplus T^{0,1}M.
$$

\begin{defn}
A differential $(p+q)-$form $\om$ is  {\em of type $(p,q)$} if it is a
section of
$$
\La^{p,q}M = (\La^p T^{1, 0}) \otimes (\La^q T^{0, 1}).
$$
\end{defn}

Let $z^i = x^i + \sqrt{-1} y^i$, $i = 1, \cdots, n$
be complex coordinates. Then, we define
$$
dz^i = dx^i + \ii dy^i  \hspace{0.5cm} \text{and} \hspace{0.5cm} d\zbar^j = dx^j - \ii dy^j,
$$
and $\frac{\partial}{\partial z^i}$, $\frac{\partial}{\partial \zbar^j}$, the dual of
$dz^i$ and $d\zbar^j$. From this definition, it is easy to see that
\begin{eqnarray*}
T^{1,0}M &=& \text{span}\left\{ \frac{\partial}{\partial z^i} \right\}^n_{i=1} \\
T^{0,1}M &=& \text{span} \left\{ \frac{\partial}{\partial \zbar^i}  \right\}^n_{i=1}
\end{eqnarray*}

\begin{exer}
We consider the decomposition of forms into $(p,q)$-type. Show that $N_J=0$ if and only if $(d\alpha)^(0,2)=0$ for every form $\alpha$ of type $(1,0)$.
\end{exer} 

Exterior differentiation
$$d: \Omega^{p,q}M \rightarrow \Omega^{p+1,q}M \oplus \Omega^{p,q+1}M $$
also splits according to this decomposition:
\begin{eqnarray*}
\partial:&& \Omega^{p,q}M \rightarrow \Omega^{p+1,q}M \\
\bar{\partial}:&& \Omega^{p,q}M \rightarrow \Omega^{p,q+1}M
\end{eqnarray*}
so that $d=\partial+\bar{\partial}$. 

The identity $d^2 = 0 $ implies that $\del^2 = 0 $, $\delbar^2 = 0$, and
$\del \delbar + \delbar \del  = 0$. Using the second identity, we may define
the {\em Dolbeault Cohomology groups $H^{p,q}_{\delbar}$} of a complex
manifold by
\beeq
H^{p,q}_{\delbar} = \frac{\Ker[\delbar: C^\infty(\La^{p,q}M) \rightarrow C^\infty(\La^{p,q+1}M)]}{\text{Im}[\delbar: C^\infty(\La^{p-1,q}M) \rightarrow C^\infty(\La^{p,q}M)]}
\eeeq

If $M$ is \kahler, we recall that  the \kahler (closed) form
$\om_g$ is given, in local coordinates, by
$$
\frac{\sqrt{-1}}{2} \sum_{i, j = 1}^n g_{i\jbar}dz^i \wedge d\zbar^j,
$$
where $g_{i\jbar} = g(\frac{\partial}{\partial z^i}, \frac{\partial}{\partial \zbar^j})$ 
\footnote{Here, we are abusing notation and writing $g$ for both the Riemannian metric
and its complex extension to $TM_\CC$}.

A last remark we make is that, due to $J$-invariance, the coefficients of type
$g_{ij}$ and $g_{\jbar\lbar}$ vanish identically. This cancelation phenomenon also
happens for some of the coefficients of the Riemann curvature tensor. This is the subject
of the next section.

\subsection{Curvature and its contractions on a \kahler manifold}

Many of the results in this section will be presented without a proof,
mainly because they involve direct calculations using the definition.

The {\em Complex Christoffel symbols} are defined in analogy with the Riemannian
version. We denote by $\gr_\CC$ (for simplicity, just $\gr$)  the Levi-Civita
connection of the hermitian metric $g$.
\begin{eqnarray*}
\gr_{\ddzi}\frac{\partial}{\partial z^j} & = &
\left(
\Ga^l_{ij} \frac{\partial}{\partial z^l}
+ \Ga^{\lbar}_{ij} \frac{\partial}{\partial \zbar^l}
\right)\\
\gr_{\ddzi}\frac{\partial}{\partial \zbar^j} & = &
\left(
\Ga^l_{i\jbar} \frac{\partial}{\partial z^l}
+ \Ga^{\lbar}_{i\jbar} \frac{\partial}{\partial \zbar^l}
\right)
\end{eqnarray*}

A more useful expression for the Christoffel symbols is given by the
following lemma.

\begin{lemma}
In holomorphic coordinates, the Christoffel symbols are
$$
\Ga^k_{ij} = \frac{1}{2}g^{k\lbar}
\left(
\ddzi g_{j\lbar} + \frac{\partial}{\partial \zbar^j}g_{i\lbar}
-\frac{\partial}{\partial \zbar^j} \frac{\partial}{\partial \zbar^l}g_{ij}
\right)
=
g^{k\lbar} \ddzi g_{j\lbar}.
$$
and $\Ga^k_{ij} = \Ga^k_{ji}$

All the coefficients are zero, except for the ones of the form $\Ga^k_{ij}$
or $\Ga^{\kbar}_{\ibar \jbar}$.
\end{lemma}

Let $R(g) = R_{i\jbar k \lbar}$ be the coordinates of the $(4,0)$-Riemann curvature tensor of the metric $g$
written in holomorphic coordinates. It is useful to know the following expression, writing $R(g)$
in terms of the metric $g$.

\begin{lemma}
The components of the \kahler Riemann tensor are given by
$$
R_{i\jbar k \lbar} =
-\frac{\partial^2}{\partial z^i \partial \zbar^j} g_{k\lbar} + g^{u\vbar}\ddzi g_{k\vbar}
\frac{\partial}{\partial d\zbar^j} g_{u\lbar}.
$$
\end{lemma}

Note that the only non-vanishing terms have two barred indices exactly. The vanishing of the
others has to do with the fact that $R(X,Y)$ is $J$-invariant on a \kahler manifold.

We define the {\em (complex) Ricci curvature tensor} of the metric $g$ as being the trace of the
Riemann curvature tensor. Its components in local coordinates can be written
as
\begin{equation}
\label{e.Riccidef}
\Ric_{k\lbar} =  \sum_{i, j = 1}^n g^{i\jbar} R_{i\jbar k \lbar} =
-\frac{\partial^2}{\partial z_k \partial \zbar_l} \log \det (g_{i\jbar}).
\end{equation}

The {\em Ricci form} associated to $g$ can then be defined by setting
\beeq
\label{e.ricci}
\Ric = \sum_{i, j = 1}^n Ric_{i\jbar} dz^i \wedge d\zbar^j.
\eeeq
in local coordinates.

Finally, we define the {\em Laplacian} acting on functions  to
be given by
$$
\De  = g^{i{\jbar}}\gr_i \gr_{\jbar} = g^{i\jbar} \frac{\partial^2}{\partial z^i \partial z^j}.
$$

\subsection{Hodge Theory for \kahler manifolds}

Recall from Section \ref{s.stilltobemade} that the space $\cH^k$ of Hodge $k$-forms
decomposes into a direct sum of irreducible representations of $\Hol(g)$.  We will work out
a few of the details for the case
when $\Hol(g) = U(n)$, and derive a few constraints in the topology of a compact \kahler manifold.
Let $M$ be a compact \kahler manifold, and define the vector space $\cH^{p,q}$ of {\em harmonic forms
of type $(p,q)$} by
$$
\cH^{p,q} = \Ker\left(     \De: C^\infty(\wedge^{p,q}M) \rightarrow  \wedge^{p,q}M)\right).
$$

\noindent
{\bf Exercise:} Show that $\alpha \in \cH^{p,q}$ if and only if $\del \alpha =
\delbar \alpha = \del^* \alpha = \delbar \alpha  = \delbar^*\alpha = 0$.

The following result is the \kahler analog of the Hodge Decomposition Theorem,
and its proof can be found in \cite{Wells}.
\begin{thm}
Let $M$ be a compact \kahler manifold. Then
$$
C^\infty(\wedge^{p,q}M) =
\cH^{p,q} \oplus \del[ C^\infty(\wedge^{p-1,q}M)] \oplus\delbar[ C^\infty(\wedge^{p,q-1}M)] .
$$
\end{thm}
As a consequence of the theorem above, we see that the Dolbeault groups  $H_{\delbar}^{p,q}$
are isomorphic to $\cH^{p,q}$. We now define, for the complex Laplacian $\De$,
$$
\cH^{k} = \Ker\left(     \De: C^\infty(\wedge^{k}M \otimes \CC) \rightarrow  \wedge^{k}M \otimes \CC)\right).
$$

Since the real and complex laplacian differ by a constant factor, there is an isomorphism
between $\cH^k$ and the complex cohomology $H^k(M, \CC)$.  We then define
$H^{p,q}(M) \subset H^{p+q}(M, \CC)$ to be the subspace with representatives in $\cH^{p,q}$.
Then,
\begin{thm}
\label{t.hodge}
Let $M$ be a compact $2n$-dimensional \kahler manifold. Then,
$H^k(M,\CC) \bigoplus^k_{j=0} H^{k, j-k}(M)$, and
any element in $H^{p,q}(M)$ can be represented by a harmonic
$(p,q)$-form. Furthermore, we have $H^{p,q}(M) \cong H^{p,q}_{\delbar}(M)$ and
\begin{eqnarray}
\label{e.hodge}
H^{p,q}(M) &\cong &  \overline{H^{q,p}(M)}\\
H^{p,q}(M) &\cong &(H^{n-p,n-q}(M))^*.
\end{eqnarray}

\end{thm}

Note that the definition of $H^{p,q}(M)$ does not depend on the choice of \kahler
structure, only on the complex structure. However,  the relationships of the cohomology groups
in (\ref{t.hodge}) are only true for \kahler manifolds.

We define the {\em Hodge numbers}, in analogy with the Betti numbers $b_k$, by
$h^{p,q} = \dim H^{p,q}(M)$. So, equation (\ref{e.hodge}) implies that
\beeq
\label{e.hodgerestr}
h^{p,q} = h^{q,p} = h^{n-p,n-q} = h^{n-q, n-p}.
\eeeq
Because of the relations above, we see that some compact, complex
manifolds may not admit a \kahler metric. For example, note that
the $k$-th Betti number $b^k$
satisfy $b^k = \sum_{j=0}^k h^{j, k-j}$. For  $k$ odd,
the Betti number $b^k$ needs to be even.

\begin{exer}
Check that the Hopf surface $S$, obtained as a quotient of $\CC^2 \setminus \{ 0\}$
by a free action of $\ZZ$ (any discrete group works too),
 is a compact manifold that admits a complex
structure, but no \kahler structure.
\end{exer}

\section{$SU(n)$}
\label{s.SUn}
The metrics with holonomy group $SU(n)$ are the {\em Calabi-Yau} metrics,
\kahler metrics with zero Ricci curvature. Let's recall the main result about the existence of
such metrics in compact Calabi-Yau manifolds.

Recall that the coordinates of the Ricci tensor are given by (\ref{e.ricci}).
The {\em Ricci form} associated to $g$ can then be defined by setting
$$
\Ric = \sum_{i, j = 1}^n Ric_{i\jbar} dz^i \wedge d\zbar^j.
$$
in local coordinates. In fact, a computation shows that
the Ricci form is given by
$$
\Ric = -\frac{\ii}{2\pi}\frac{\partial^2}{\partial z^i \partial z^j}(\log \det(g)).
$$

Now, given a metric $g$, we can define a matrix-valued $2$-form $\Om$ by
writing its expression in local coordinates, as follows
\begin{equation}
\label{e.curvform}
\Om_i^j = \sum_{i, p = 1}^n g^{j\pbar}  R_{i\pbar k \lbar} dz^k \wedge d\zbar^l.
\end{equation}
This expression for $\Om$ gives a well-defined matrix of $(1,1)$-forms, to
be called the {\em  curvature form } of the metric $g$.

Following Chern-Weil Theory, we want to look at  the following expression
$$
\det\left(\text{Id} + \frac{t\sqrt{-1}}{2\pi} \Om \right) = 1 + t\phi_1(g) + t^2\phi_2(g) + 
\dots,
$$
where each $\phi_i(g)$ denotes the $i$-th homogeneous component of the left-hand side,
considered as a polynomial in the variable $t$.

Each of the forms $\phi_i(g)$ is a $(i,i)$-form, and is called the
{\em $i$-th Chern form} of the metric $g$. It is a fact (see for example
\cite{Wells} for further explanations) that the cohomology class
represented by each  $\phi_i(g)$ is independent on the metric $g$, and hence it is
a topological invariant of the manifold $M$.
These cohomology classes are called the {\em Chern classes} of $M$
and they are going to be denoted by $c_i(M)$.

{\bf Remark:} We can define more generally the curvature $\Om(E)$ of a
hermitian metric $h$ on
a general complex vector bundle $E$ on a complex manifold $M$.

Let $\nabla = \nabla(h)$ be a connection on a vector bundle $E \rightarrow M$.
Then the {\em curvature} $\Om_E(\nabla)$ is defined to be the element
$\Om \in \Om^2(M, \text{End}(E, E))$ such that the $\CC$-linear mapping
$$
\Om : \Gamma(M, E) \rightarrow \Om^2(M, E)
$$
has the following representation with respect to a frame $f$:
$$
\Om(f) = \Om(\nabla, f) = d\theta(f) + \theta(f) \wedge \theta(f).
$$
Here, $\Gamma(M, E)$ is the set of sections of the vector bundle $E$,
$\Om^2(M, E)$ is the set of $E-$valued $2$-forms, and
$\theta(f)$ is the connection matrix associated with $\nabla$ and $f$
(with respect to $f$, we can write $\nabla = d + \theta(f)$).

Similarly one defines
the Chern class $c_i(M, E)$ of a vector bundle and these are also independent on the
choice of the connection.
In fact, we use the expression ``Chern classes $c_i(M)$ of the manifold $M$''
meaning the Chern classes $c_i(M, TM)$
of the tangent bundle of M.

\bigskip

We will restrict our attention to the first Chern class $c_1(M)$ of the manifold $M$.
Note that the form $\phi_1(g)$
represents  the class $c_1(M)$ (by definition),
and that $\phi_1(M)$ is simply the trace of the
curvature form:
\begin{equation}
\label{e.phi1}
\phi_1(g) =  \frac{\sqrt{-1}}{2\pi} \sum_{i=1}^n \Om_i^i
= \frac{\sqrt{-1}}{2\pi} \sum_{i, p = 1}^n g^{i\pbar}  R_{i\pbar k \lbar} dz^k \wedge d\zbar^l.
\end{equation}

On the other hand, notice that the right-hand side of (\ref{e.phi1}) is equal to
$\frac{\sqrt{-1}}{2\pi} \Ric_{k\lbar}$, in view of (\ref{e.Riccidef}).
Therefore, we conclude that the Ricci form of a \kahler metric represents
the first Chern class of the manifold $M$.
A natural question that arises is: given a \kahler 
class
$[\om] \in H^2(M, \RR) \cap H^{1,1}(M, \CC)$ in a compact, complex manifold
$M$, and any $(1, 1)$-form $\Om$ representing $c_1(M)$, is that possible to
find a metric $g$ on $M$ such that $\Ric(g) = \Om$?

This question was addressed to by Calabi in 1960, and it was answered by Yau \cite{Yau} almost
$20$ years later.

\begin{thm}
{\bf (Yau, 1978)}
If the  manifold $M$ is compact and \kahler, then there exists a unique
\kahler metric $g$ on $M$ satisfying $\Ric(g) = \Om$.
\end{thm}

This theorem has a large number of applications in different
areas of Mathematics and Physics. Its proof is based on translating
the geometric statement into a non-linear
partial differential equation, as follows.

First fix a \kahler form $\om \in [\om]$ representing the previously given
\kahler class in $H^2(M, \RR) \cap H^{1,1}(M, \CC)$.
In local coordinates, we can write
$\om$ as $\om = g_{i\jbar}dz^i \wedge d\zbar^j$.

The $(1, 1)$-form $\Om$ is a representative for  $c_1(M)$, and
we have seen that $\Ric(\om)$
represents the same cohomology class as $\Om$.
Therefore, since $\Ric(\om)$  is also a  $(1, 1)$-form,
we have that,  due to the famous $\partial \bar \partial$-Lemma,
there exists a function $f$ on $M$ such that
$$
\Ric(\om) - \Om =  \frac{\sqrt{-1}}{2\pi} \partial \bar\partial f,
$$
where $f$ is uniquely determined after imposing the normalization
\begin{equation}
\label{e.integrability}
\int_M \left( e^f - 1 \right) \om^n = 0.
\end{equation}
Notice that $f$ is fixed once we have fixed $\om$ and $\Om$.

Again by the $\partial \bar \partial$-Lemma,
any other $(1, 1)-$form in the same cohomology class $[\om]$
will be written as
$\om + \frac{\sqrt{-1}}{2\pi} \partial \bar\partial \phi$, for some
function $\phi \in C^\infty(M, \RR)$.

Therefore, our goal is to find a representative
$\om + \frac{\sqrt{-1}}{2\pi} \partial \bar\partial \phi$ of the
class $[\om]$ that satisfies
\begin{equation}
\label{e.ric}
\Ric\left( \om + \frac{\sqrt{-1}}{2\pi} \partial \bar\partial \phi \right)
=
\Om =
\Ric(\om) - \frac{\sqrt{-1}}{2\pi} \partial \bar\partial f.
\end{equation}

Rewriting (\ref{e.ric}) in local coordinates, we have
$$
- \partial \bar\partial \log \det \left( g_{i\jbar} +
\frac{\partial^2 \phi}{\partial z_i \partial \zbar_j}
                                                  \right)
=
- \partial \bar\partial \log \det \left( g_{i\jbar}\right) - \partial \bar\partial f,
$$
or
\begin{equation}
\label{e.welldef}
 \partial \bar\partial \log
\frac{\det \left( g_{i\jbar} +
\frac{\partial^2 \phi}{\partial z_i \partial \zbar_j}
                                                  \right)}{\det \left( g_{i\jbar}\right)}
= \partial \bar\partial f.
\end{equation}

Notice that, despite of the fact that this is an expression given in
local coordinates, the term at the right-hand side of (\ref{e.welldef})
is defined globally.
Therefore, we obtain an equation
well-defined on all of $M$. In turn, this equation gives rise to the
following (global) equation
\begin{equation}
\label{e.MA100}
\left(\om + \frac{\sqrt{-1}}{2\pi}\partial \bar\partial \phi \right)^n = e^f \om^n.
\end{equation}

We shall also require positivity of the resulting \kahler form:
$$\left(\om + \partial \bar\partial \phi \right) > 0  \hspace{.2cm}
\text{ on} \hspace{.2cm}     M. $$
This equation is a non-linear partial differential equation of
Monge-Amp\`ere type, that is going to be referred to from now on as the Complex
Monge-Amp\`ere Equation.

We remark that, if $\phi$ is a solution
to (\ref{e.MA100}),  $\om + \partial \bar\partial \phi$
is the \kahler form of our target metric $g$, {\em ie},
$\Ric(g) = \Om$.
Therefore,  in order
to find metrics that are solutions to Calabi's problem, it suffices to determine a solution
$\phi$ to (\ref{e.MA100}).

The celebrated Yau's Theorem in \cite{Yau} determines a unique solution to (\ref{e.MA100})
when $f$ satisfies the integrability condition
(\ref{e.integrability}). The proof of this result is based on the continuity method,
and we sketch here a brief outline of the proof.

The uniqueness part of Calabi Conjecture was proved by Calabi in the $50$'s.
Let $\om', \om'' \in [\om]$ be representatives of the \kahler class $[\om]$
such that $\Ric(\om') = \Ric(\om'') = \Om$. Without loss of generality, we may assume
that $\om'' = \om$, and hence $\om' = \om + \partial \bar \partial u$.

Notice that
\begin{eqnarray}
\label{e.blah}
0 &=& \frac{1}{\Vol_\om(M)} \int_M u((\om')^n - \om^n) \\
&=&\frac{1}{\Vol_\om(M)} \int_M -u \partial \bar \partial u \wedge
\left[  (\om')^{n-1} + \right. \\
&& \left.(\om')^{n-2}\wedge \om + \cdots + \om^{n-1}   \right].
\end{eqnarray}

\noindent However $\om'$ is a \kahler form, so that $\om'>0$. We then conclude that
the right-hand side of (\ref{e.blah}) is bounded from below by
$$\frac{1}{\Vol_\om(M)} \int_M -u \partial \bar \partial u \wedge w^{n-1}.$$
Therefore,
\begin{eqnarray}
0 &\geq &\frac{1}{\Vol_\om(M)} \int_M -u \partial \bar \partial u \wedge w^{n-1}
\\
&=& \frac{1}{n\Vol_\om(M)} \int_M |\partial u|^2 w^{n}
\\
&=& \frac{1}{2n\Vol_\om(M)} \int_M |\nabla u|^2 w^{n},
\end{eqnarray}
implying that $|\nabla u| = 0$, hence $u$ is constant, proving the uniqueness
of solution to (\ref{e.MA100}).

\bigskip

Let us now consider the existence of solution to  (\ref{e.MA100}).
Define, for all $s \in [0,1]$, $f_s = sf + cs$, where
the constant $c_s$ is defined by the requirement that $f_s$ satisfies the
integrability condition $\int_M[e^{f_s} - 1 ]\om^n = 0 $.

Consider the family of equations
\begin{equation}
\label{e.fs}
\left( \om + \partial \bar \partial u_s\right)^n = e^{f_s} \om^n.
\end{equation}
We already prove that the solution $u_s$ to (\ref{e.fs}) is unique, if
it exists.

Let $A = \{s \in [0,1];$ (\ref{e.fs}) is solvable for all $t \leq s \}$.
Since $A \neq \emptyset$, we just need to show that $A$ is open and closed.

{\bf Openness:} Let $s \in A$, and let $t$ be close to $s$. We want to show that
$t \in A$. In order to do so, let $\om_s = \om + \partial \bar \partial u_s$, for $u_s$ a
solution to  (\ref{e.fs}).
We define the operator $\Psi = \Psi_s$ by
$$
\Psi: X \rightarrow Y; \hspace{2cm}
\Psi(g) = \log\left(\frac{(\om_s + \partial \bar \partial g)^n}{w_s^n}\right),
$$
where $X$ and $Y$ are subsets (not subspaces) of
$C^{2, 1/2}(M.\RR)$ and $C^{0, 1/2}(M.\RR)$
satisfying some extra non-linear conditions.

The linearization of $\Psi$ about $g = 0$ is simply the
metric laplacian with respect to the metric $\om_s$.
By the Implicit Function Theorem, the invertibility of the
laplacian (a result that can be found in \cite{GT}, for example)
establishes the claim.

{\bf Closedness:}
The proof that $A$ is closed is a deep result, involving complicated
{\em a priori} estimates. A reference for this proof is Yau's paper itself \cite{Yau},
or for a more detailed proof, the books \cite{Tianbook} and \cite{Asterisque}.

\section{$Sp(n)$}
\label{s.Spn}

The metrics $g$
with holonomy group $Hol(g) \subset Sp(n)$
are the {\em \hk} metrics. Since $Sp(n) \subset SU(2n)$,
all the \hk manifolds are necessarily \kahler and Ricci-flat.
But more is true: we will see that a \hk metric is not only
\kahler with respect to one complex structure on the manifold,
but to a whole $2$-sphere of complex structures.

Let us start by describing $Sp(n)$.
In order to do so, we now  recall some basic facts about  the quaternions
$\HH$.

\begin{defn}
The {\em quaternions $\HH$} are defined by
\beeq
\HH = \{ x_0 + x_1 i + x_2 j + x_3 k \st x_j \in \RR\} \cong \RR^4.
\eeeq
The {\em imaginary quaternions} are given by $\Im \HH = <i, j, k> \cong \RR^3$,
with multiplication given by
$$
ij = -ji = k , \, \, \,  jk = -kj = i   , \, \, \, ki = -ik = j , \, \, \, i^2 = j^2=k^2 = 1.
$$
\end{defn}

We can define a metric $g$ and $2$-forms $\om_1, \om_2$ and $\om_3$ on $\HH^n$ as follows.
If $q_1, \cdots, q_n$ are coordinates on $\HH^n$ such that
$q_\ell = x^0_\ell +   x^1_\ell i +  x^2_\ell j +  x^3_\ell k $, then
\begin{eqnarray*}\label{e.omhk}
g = \sum_{\ell = 1}^n \sum_{p=0}^3 (dx_\ell^p)^2    \, \, \,\, \, \,\, \, \,  \, \, \,\, \, \,\, \, \,&   &
\om_1 =  \sum_{\ell = 1}^n dx_\ell^0 \wedge  dx_\ell^1 + dx_\ell^2 \wedge  dx_\ell^3 \\
\om_2 =  \sum_{\ell = 1}^n dx_\ell^0 \wedge  dx_\ell^2 + dx_\ell^1 \wedge  dx_\ell^3 &  \, &
\om_3 =  \sum_{\ell = 1}^n dx_\ell^0 \wedge  dx_\ell^3 + dx_\ell^1 \wedge  dx_\ell^2
\end{eqnarray*}

Identifying $\HH$ with $\RR^4$, and abusing notation and writing the complex structures $i, j$
and $k$
 on $\RR^4$ defined by left multiplication by $i, j$ and $k$, we see that $g$ is \kahler with respect to each complex structure $i, j$ and $k$.

 The subgroup of $GL(4n, \RR)$ that preserves $(g, \om_1, \om_2, \om_3)$ is $Sp(n)$,
 a compact, connected, simply-connected, semi-simple Lie Group of dimension
 $2n^2 +n$.

 \begin{exer}
Setting $x^0_\ell +   x^1_\ell i = z_{2\ell -1}$ and
 $x^2_\ell +   x^3_\ell i = z_{2\ell}$, we can define coordinates
 in $\CC^2$ that will allow an identification $\HH \cong \CC^2$.
 Using this identification, show that $Sp(n) \subset SU(2n)$.
\end{exer}

Now,  we can define the notion of a \hk structure on a $4n$- dimensional manifold $M$.
 Let $(M, g) $ be  a
$4n$-dimensional riemannian manifold with holonomy group
contained in $Sp(n)$.
Since $Sp(n)$ preserves 3 complex structures
$i,j$ and $k$ in $\RR^{4n}$, we can see that there are three corresponding
almost complex structures $J_1, J_2$ and $J_3$ on $M$,
whose covariant derivative with respect to the Levi-Civita connection of $g$ is zero, and
such that $J_1 J_2 = J_3$.
The metric $g$ is \kahler with respect to each of them {\footnote{In fact, $g$ is
\kahler with respect to $a_1J_1 + a_2 J_2 + a_3J_3$, for any real numbers $a_1,
a_2, a_3$ such that $a_1^2 + a_2^2 + a_3^2 = 1$.}}, and we will refer to it as
a {\em \hk metric}, and $(J_1, J_2, J_3, g)$ as a {\em \hk structure} on $M$.

Note that if we do not require that the almost complex structures are integrable
(or equivalently, as we will see in Proposition \ref{p.hkTFAE}, that they are not covariantly constant),
we will refer to $(J_1, J_2, J_3, g)$ as a {\em almost \hk structure} on $M$.
Each of those almost complex structures $J_i$ defines a hermitian form $\om_i$.

The following proposition is a characterization of \hk manifolds, and its proof can be found
in \cite{Salamon}.
\begin{prop}
\label{p.hkTFAE}
Let $M$ be a $4n$-dimensional manifold, and $(J_1, J_2, J_3, g)$ is an
almost  \hk structure on $M$. Let $\om_1, \om_2, \om_3$ be the associated hermitian forms.
Then, the following are equivalent:
\been
\item $(J_1, J_2, J_3, g)$ is a
 \hk structure on $M$.
\item
 $d\om_1 = d\om_2 = d\om_3 =0$.
\item
$\nabla \om_1 = \nabla \om_2 = \nabla \om_3 =0 $.
\item
$Hol(g) \subset Sp(n)$.
\eeen

\end{prop}

\subsection{Twistor Spaces}

Let $(M, J_1, J_2, J_3, g)$ be a \hk manifold of (real) dimension $4n$.
 We can associate to $M$
a $(2n+1)$-dimensional complex manifold $Z$, called the {\em twistor space of $M$},
in the following way.

Recall that the space of complex structures on $M$ is an $S_2$, that we can think of
$\CC\PP^1$, that carries a natural complex structure $J_0$.  Set $Z = \CC\PP^1 \times M$,
where each point in $z \in Z$ is of the form $z =(J, x)$, for $J \in \CC\PP^1$ and $x \in M$.

Then, we can define a complex structure $J_Z = J_0 \oplus J$ at the tangent space
$T_z Z = T_j \CC\PP^1 \oplus T_xM$. It turns out that $J_Z$ is integrable, and $(Z, J_Z)$
is a complex manifold, called the {\em twistor space of $M$}.

The twistor space $Z$ of $M$ carries two important structures:
\vspace {0.1in}

\noindent $1)$ If $p: Z \rightarrow \CC\PP^1$ and $\pi: Z \rightarrow M$ are the natural projections,
then $p$ is holomorphic. We can define a free holomorphic anti-involution $\si: Z \rightarrow Z$
by $(J, x) \mapsto (-J,x)$.  Also, for each $x\in M$, the fiber $\pi^{-1}(x)$ is  a holomorphic curve
in $Z$, isomorphic to $\CC\PP^1$, with normal bundle $2n \cO(1)$ and preserved by the involution $\si$ .

\vspace{0.1in}
\noindent $2)$ Let $\cD$ be the kernel of $dp: TZ \rar T\CC\PP^1$.
Then $\cD$ is a holomorphic subbundle of $TZ$, and so one can construct
a non-degenerated holomorphic section $\om$ of the holomorphic vector bundle
$p^*(\cO(2)) \otimes \La^2 \cD^*$ over $Z$.

What $2)$ actually means is that, for each fiber $p^{-1}(J)$, to choose a complex
symplectic form $\om$ on it.

\begin{exer}
What is $\om$ for each choice of $J_1$, $J_2$ and $J_3$?
\end{exer}

The nice feature about twistor spaces is that we can use $(Z, p, \si, \om)$ to reconstruct
$(M, g, J_1, J_2, J_3)$:

\begin{thm}
Let $(Z,J_Z)$ be a $(2n+1)$-dimensional complex manifold, equipped with
\been
\item
A holomorphic projection $p: Z \rightarrow \CC\PP^1$,
\item
A holomorphic section $\om$ of $p^*(\cO(2)) \otimes \La^2 \cD^*$,
where $\cD = \Ker(dp)$,
\item
A free antiholomorphic involution $\si: Z \rar Z$
such that $\si^*(\om)= \om$ and $p \circ \si = \si' \circ p$, where
$\si': \CC\PP^1 \rar \CC\PP^1$ is the antipodal map.
\eeen
Let $M'$ be the set of holomorphic curves in $Z$ which are isomorphic to
$\CC\PP^1$, preserved by $\si$, and with normal bundle $2n \cO(1)$.

Then $M'$ is a hyper-complex $4n$-manifold with a natural pseudo-riemannian metric $g$.
If $g$ is indeed riemannian, then $M'$ is \hk.

Furthermore, let $Z'$ be the twistor space of $M'$. Then, there is a natural map
$\iota: Z \rar Z'$ that is a local biholomorphism, and such that $p$, $\si$ and $\om$
are identified with their analogues in $Z'$,
\end{thm}

For the proof, we refer to \cite{J}.

\subsection{K3 Surfaces}

For the next few sections, the main references we suggest are \cite{J} and \cite{Salamon}.

The $4$-dimensional example of a \hk manifold are the K3 surfaces, named after
\kahler, Kummer and Kodaira, and as a reference to K2, the second highest mountain in the
world, located in South Asia.

\begin{defn}
A {\em K3 surface} $S$ is a compact, complex $2$-dimensional manifold
with $h^{1,0} = 0$ and trivial canonical line bundle.
\end{defn}
Note that our definition does not include the requirement that $S$ is \kahler.
The result that all K3 surfaces are indeed \kahler is due to Siu.

One example of a K3 surface is the {\em Fermat quartic}
$$
F = \{ [z_0, \cdots, z_3] \in \CC\PP^3 ; z_0^4 + \cdots z_3^4 = 0\}.
$$
We can use the adjunction formula $K_F = (K_{\CC\PP^3} \otimes L_F)|_F$
to conclude that the canonical line bundle of $F$
 is trivial.

 For this example, Riemann-Roch tells us that $\Xi(F) = 24$. Since $b^0 = 1$
 and $b^1 = 0$ (as $F$ is \kahler and $h^{1,0} = 0$), we learn that $b^2 = 22$.

 Also, since the canonical line bundle $K_F$ of $F$ is trivial, then $h^{2,0} = 1 = h^{0,2}$,
 which implies that $h^{1,1} = 20$. Further, the signature $\tau(F)$, defined by
 $\tau(F) = b_+^2 - b^2_-$, is given by
 $\tau(F) = \sum_{p, q = 0}^2 (-1)^p h^{p,q} = -16$, which tells us that $b_+^2 = 3$
and $b_-^2 = 19$.

The following result tells us that the same topological facts listed above hold true
for any K3 surface.

\begin{thm}
Let $S$ be a K3 surface. Then $S$ is simply connected, with
Betti numbers $b_+^2 = 3$ and $b_-^2 = 19$. Moreover,
$S$ is \kahler, with Hodge numbers $h^{2,0}= h^{0,2} = 1$, and $h^{1,1} = 20$.
All K3 are diffeomorphic.
\end{thm}

The proof of this theorem is a consequence of two main results: the first is the result
by Kodaira \cite{Kodaira} that proved that all K3s are diffeomorphic to the Fermat hypersurface
$S$ described above. Secondly, Siu \cite{Siu} showed that a K3 surface is necessarily \kahler.
The theorem follows trivially from the assertions above.

\section{$Sp(n)Sp(1)$}
\label{s.SpnSp1}

We begin this section by recalling what is the action of $Sp(n)Sp(1)$ on $\HH^n$: by looking
at $Sp(n)$ as the group of $n \times n$ matrices with quaternionic coefficients, we can describe
the action of $(A, q) \in Sp(n)\times Sp(1)$ on $\HH^n$ by
$$
(q_1, \cdots, q_n) \mapsto q (q_1, \cdots, q_n) A^H.
$$
Clearly, this action pushes down to
$$
Sp(n) Sp(1)   = \left( Sp(n)\times Sp(1) \right) / \{\pm (I_n, 1)  \},
$$
where $I_n$ is the $n\times n$ identity matrix.

As pointed out by Salamon \cite{Salamon}, one of the main interests in studying
quaternionic \kahler manifolds comes from an observation by Wolf \cite{Wolf} that
given any compact simple Lie Group $G$, there exists a quaternionic \kahler symmetric space
$G /H$.

\begin{defn}
We say that a $4n$-dimensional  Riemannian manifold $(M, g)$ ($n>1$) is
a {\em quaternionic \kahler manifold} if $Hol(g) \subset Sp(n)Sp(1)$.

\end{defn}

It should be pointed out that quaternionic \kahler manifolds may not be \kahler,
as $Sp(n)Sp(1) $ is not a subgroup of $U(2n)$. In fact, the quaternionic \kahler
manifolds should be seen as a generalization of \hk manifolds: the tangent space still
admits three complex structures $J_1, J_2, J_3 $ (satisfying the same composition rules
as the imaginary quaternions $i, j, k$), but there is no preferred choice of basis. It is
possible to make a smooth local choice of $J_1, J_2, J_3 $, but those will not be all
complex structures on $M$.

The following result is proved in Salamon \cite{Salamon}:
\begin{prop}
If $(M, g)$ is quaternionic \kahler, then it is necessarily Einstein. Furthermore,
if $g$ is also Ricci-flat, then $Hol^0(g) = Sp(n)$.
\end{prop}
Because of the proposition above, we can split the class of quaternionic \kahler manifolds
into positive or negative, according to the sign of the (constant) scalar curvature $s$ of g.
The case $s = 0$ corresponds to $(M, g)$ being locally \hk, according to the above result.

\bigskip

Now, for the $2$-forms  $\om_1, \om_2, \om_3$  defined in (\ref{e.omhk}), we see that
$Sp(n)Sp(1)$ preserved the $4$-form
\beeq
\Om_0 = \om_1 \wedge \om_1 + \om_2 \wedge \om_2 + \om_3 \wedge \om_3.
\eeeq
on $\HH^n$.

Then, by Prop \ref{prop:invariantelements}, we see that any quaternionic \kahler manifold
has a  constant $4$-form $\Om$ which is isomorphic to $\Om_0$ at each point.
The stabilizer of $\Om_0$ in $GL(4n, \RR)$ is $Sp(n)Sp(1)$, so $Sp(n)Sp(1)$-structures
on  a $4n$-dimensional manifold $M$ are equivalent to $4$-forms $\Om$ which are
 isomorphic to $\Om_0$ at each point.

By a result of Swann, if $n \geq 3$, the $Sp(n)Sp(1)$-structure is torsion-free
if and only if $d\Om = 0$ ({\em ie}, $\nabla \Om  = 0
\Leftrightarrow d\Om = 0$).

\subsection{Twistor Spaces}
Quaternionic \kahler manifolds $(M,g)$, as well as \hk, have also twistor spaces.
The manifold $M$ has a principal $Sp(n)Sp(1)$-bundle $P$, a subbundle of the
frame bundle of $M$.
We define
$$Z = P /Sp(n)U(1),$$
 so $Z$ is a bundle $\pi:Z \rar X$ with fiber
$$ Sp(n) Sp(1) / Sp(n) U(1) \cong S^2.$$ Let us make this identification more explicit.
For a point $z \in Z$, we can associate a complex structure on $T_{\pi(z)}M$
as follows: for any $x \in M$, due to the natural identification $T_{x}M \cong \RR^{4n}$,
$\Om \cong \Om_0$ and $g$ identified with the Euclidean metric,
we can naturally identify $\pi^{-1}x$ with $a_1 J_1 +a_2 J_2 +a_3 J_3 $
in $\RR^{4n}$, for $(a_1, a_2, a_3) \in S^2$.

\begin{exer}
Check that changing the choice of identification $T_{x}M \cong \RR^{4n}$
are related by the action of $Sp(n) Sp(1)$, so that the construction is well defined.
\end{exer}

We end this section  with the following theorem due to Salamon \cite{Salamon}
\begin{thm}
Let $(M, g)$ be a quaternionic \kahler manifold, and let $Z$ and $\pi$ be defined as above.
Then $Z$ has the following structures:
\beit
\item An integrable complex structure $J$, so that $Z$ is a complex manifold. Further,
each fiber $\pi^{-1}(x)$, $x \in M$, is a $\CC\PP^1$ in $Z$ with normal bundle $2n \cO(1)$;
\item An anti-holomorphic involution $\si: Z \rar Z$ such that $\pi \circ \si = \pi$, and
acting as $J \mapsto -J$;
\item If $g$ has nonzero scalar curve, then $Z$ carries a holomorphic
contact structure $\tau$.
\item If the scalar curvature of $g$ is positive, then $Z$ carries a \KE metric
$h$ with positive scalar curvature.
\eeit

Moreover, we can completely recover the quaternionic \kahler structure on $M$
(up to homothety) from the data $(Z, \pi, \si, \tau)$.
\end{thm}

\section{$G_2$}
\label{s.G2}

We have seen in Section \ref{ss.Berger} that riemannian metrics with exceptional holonomy $G_2$ and $Spin(7)$, 7 and 8-dimensional resp., are necessarily Ricci-flat.
These last two sections will be dedicated to those groups.

Manifolds with holonomy  $G_2$ were  first introduced by Edmond Bonan in 1966. The first complete, noncompact 7-manifolds with holonomy $G_2$ were constructed by Robert Bryant and Salamon in 1989. The first compact 7-manifolds with holonomy $G_2$ were constructed by Dominic Joyce in 1994, and compact  manifolds are sometimes known as "Joyce manifolds", especially in the physics literature.

These manifolds have been the subject of interest of theoretical physicists since the introduction of $M$-theory. In regular String Theory,
the model of Universe is $10$-dimensional: the usual Minkowsky space $\R^4$ times a Calabi-Yau manifold of very small diameter. M-Theory,
an extension of String Theory, replaces the previous model with  an 11-dimensional $\R^4$ times a $G_2$ manifold, also with diameter of
the order of $10^{-33}cm$. A good introductory reference  for the Physics-oriented student is \cite{BasicsofMTheory}.

\medskip

Let us recall the definition of the group $G_2$.

\begin{defn}
Let $(x_1, \cdots, x_7)$ be coordinates for $\R^7$. Define a $3$-form $\vp_0$
by
$$
\vp_0 = dx_{123} +  dx_{145} + dx_{167} + dx_{246} - dx_{257} - dx_{347} - dx_{356},
$$
where $dx_{ijk}$ denotes the $3$-form $dx_i \wedge dx_j \wedge dx_k$.

The subgroup of $GL(7, \R)$ which preserves the  form $\vp_0$ is the 14-dimensional exceptional Lie Group $G_2$.
It is compact, connected, semi-simple, and it fixes the Euclidean metric $g_0 = dx_1^2 + \cdots + dx_7^2$, the orientation,
and the $4$-form $* \vp_0$, the Hodge dual of the $3$-form $\vp_0$ with respect to the metric $g_0$.

\end{defn}

Once we understand the structure in $\R^7$, we can move on to manifolds. So, let $M$ be a $7$-dimensional
riemannian manifold, and let $p \in M$. We define
\begin{eqnarray*}
\cP^3_p M  & = & \{  \vp \in \wedge^3 T^*_p M ;  \text{ there exists an oriented  isomorphism}   
\\ &  & f: T_pM \rightarrow \R^7
\, \text{with} \, \vp  =  f^*\vp_0\}.
\end{eqnarray*}

Clearly, $\cP^3_p M$ is isomorphic to the quotient $GL(7, \R) / G_2$, which has dimension $49 - 14 - 35$. Since $35$ is also
the dimension of $\wedge^3 T^*_pM $, we see that  $\cP^3_p M$ is an open subset of $\wedge^3 T^*_pM $.

Let $\cP^3M $ be a fiber bundle over $M$, with fiber $\cP^3_p M$. Then
$\cP^3M $ is an open subbundle of $\wedge^3 T^*_pM $ with fiber $GL(7, \R) / G_2$.
We will say that a $3$-form $\vp$ is {\em positive} if $\vp|_p \in \cP^3M $.

Let $F$ be the frame bundle of $M$ (with fiber, at each $p \in M$, the set of all
isomorphisms from the tangent space at $p$ to $\R^7$), and let $\vp \in \cP^3M$ be a positive $3$-form.

Consider $Q$, a subbundle of $F$, such that its fiber at $p$ consists
of the isomorphisms $f: T_pM \rightarrow \R^7$ that preserve $\vp$.

\begin{exer}
Show that $Q$ is a principal subbundle of the frame bundle $F$, with fiber $G_2$.
\end{exer}

In other words, the exercise above completes the proof that fixing a positive
$3$-form $\vp$ determines a $G_2$-structure $Q$ on $M$.

The converse is also true: given a $G_2$-structure, we can define a metric $g$,
a $3$-form $\vp$ and a $4$-form $*\vp$ by requiring that $\vp_0, * \vp_0$
and $g_0$ are preserved under action of $G_2$. The form $\vp $ will be positive definite
if and only if $Q$ is an oriented $G_2$-structure.

Hence, we will from now on refer to  $(\vp, g)$ as a $G_2$ structure on $M$.

\begin{defn}
Let $M^7$ be  a  riemannian $7$-manifold, endowed with a $G_2$ structure $(\vp, g)$.
Let $\nabla_g$ be the Levi-Civita connection associated to $g$, and $\nabla_g \vp$ be the
{\em torsion} of the $G_2$ structure $(\vp, g)$.

We say that $(\vp, g)$ is {\em torsion-free} if  $\nabla_g \vp =0$.

A {\em $G_2$-manifold} is a triple $(M, \vp, g)$ such that $(\vp, g)$ is a torsion-free $G_2$-structure.

\end{defn}

The following proposition is  a characterization of $G_2$-manifolds, and its proof can be found at \cite{Salamon}.

\begin{prop}
\label{p.TFAEG2}
Let  $M$ be a riemannian $7$-manifold, and let $( \vp, g)$ be a $G_2$-structure on $M$. Then,
the following statements are equivalent:
\begin{enumerate}
\item
 $( \vp, g)$ is torsion-free.
 \item
$ hol(g) \subset G_2$, and $\vp$ is the induced $3$-form by the inclusion.
\item
$d\vp = 0 $ and $d^*\vp = 0$ on $M$.
\item
$d\vp = 0$ and  $ d(*\vp) = 0$ on $M$.

\end{enumerate}
\end{prop}

It should be noted that even though the condition $\nabla_g \vp  = 0$ looks
linear on $\vp$, it is not, because $\nabla_g$ depends on $g$, and so does the
Hodge star operator.

Recall that in Section \ref{s.stilltobemade}, we saw that a $G$-structure induces a splitting
of the bundle of forms on $M$ into irreducible components. For a manifold
with a $G_2$-structure, this take the shape described in the following result.

\begin{prop}
\label{p.wedgesplit}
For $M$ a manifold with $G_2$-structure $(\vp, g)$, each bundle $\wedge^k T^*M$  of
$k$-forms splits orthogonally as follows:
\begin{itemize}
\item
$
\wedge^1 T^*M = \wedge^1_7 $
\item
$\wedge^2 T^*M = \wedge^2_7  \oplus \wedge^2_{14}$
\item
$\wedge^3 T^*M = \wedge^3_1  \oplus \wedge^3_7  \oplus \wedge^3_{27}$,
\end{itemize}
where $\wedge^k_\ell  $ is an irreducible representation of $G_2$ of dimension $7$.

The Hodge star operator $*$ provides an isometry between
$\wedge^k_\ell  $ and $\wedge^{7-k}_\ell  $, and we have
 $\wedge^3_1 = \langle \vp \rangle  $ and
$\wedge^4_1 = \langle *\vp \rangle  $
\end{prop}

As discussed in Section \ref{ss.Berger}, the curvature tensor of a riemannian
metric is subject to very strong restrictions, once we fix the holonomy group
of the metric.

\begin{lemma}
If $(M, g)$ is a riemannian manifold with curvature tensor
$R_{ijkl}$, then for each $p \in M$, $R_{ijkl}$ lies on
 $S^2\hol \subset \wedge^2 T_x^*M \times \wedge^2 T_x^*M $.
\end{lemma}

This result implies

\begin{prop}
\label{p.G2isRicciFlat}
If  $(M, g)$ is a riemannian 7-manifold whose holonomy group
$Hol(g) \subset G_2$, then $g$ is Ricci-flat.

\end{prop}

\bigskip

In fact, having $Hol(g)$ contained in $G_2$ implies more. Since $G_2$
is a simply-connected subgroup of $SO(7)$, then any manifold with
a $G_2$-structure is spin, with a preferred spin structure.

Recall that $SU(2) \subset SU(3) \subset G_2$, and there are strong limits for the possibilities for the 
holonomy groups of a manifold with holonomy contained in $G_2$ :

\begin{prop}
\label{p.holcharactG2}
If $(\vp, g)$ is a torsion-free $G_2$ structure on $M$,
then $Hol^0(g)$ can only be one of the following subgroups of
$G_2$: $\{ 1\}$, $SU(2)$, $SU(3)$ or $G_2$.
\end{prop}

The proposition above tells us that from a $K3$ surface or
a Calabi-Yau $3$-fold, we can manufacture a $G_2$ manifold. This
simple result turns out to be a useful tool  to produce compact
examples of manifolds with holonomy group equal to $G_2$.

More explicitly, let $(M, J, h, \Om)$ be a $K3$ surface, where $h$ is
the (ricci-flat) \kahler metric with associated \kahler form $\om$,
and $\Om$ is the holomorphic $2$-form.
Let $(x, y, z)$ be coordinates for $\R^3$. We define
\begin{eqnarray*}
g &=& dx^2 + dy^2 + dz^2 + h; \\
\vp &=& dx \wedge dy\wedge dz + dx \wedge \om + dy \wedge \text{Re}(\Om) + \text{Im}(\Om) \wedge dz
\end{eqnarray*}
Then, $(\vp, g)$ is a torsion-free $G_2$ structure on $M \times \R^3$.

Similarly, for  $(M, J, h, \Om)$ a Calabi-Yau $3$-fold,
with $h$,
the (ricci-flat) \kahler metric with associated \kahler form $\om$,
and $\Om$, the holomorphic $3$-form, and $x \in \R$,
\begin{eqnarray*}
g &=& dx^2 + h; \\
\vp &=& dx \wedge \om + \text{Re}(\Om)
\end{eqnarray*}
define a torsion-free $G_2$ structure on $M \times \R$.

\subsection{Topological Properties of compact $G_2$ manifolds}

We list here a few out of the many topological properties of
compact $G_2$ manifolds.

\begin{prop}
Let $(M, \vp, g)$ be a compact $G_2$ manifold.

Then $\Hol(g) = G_2$ if and only if the fundamental group of $M$ is finite.

\end{prop}

\proof{
Note that $M$ must be Ricci-flat, which implies that
$M$ has a finite cover isomorphic to $T^k \times N$,
where $N$ is simply connected.

This splitting implies that $\pi_1(M) = F \ltimes \Z^k$,
where $F$ is finite. Hence, $\pi_1(M)$ is finite if and
only if $k = 0$.

Now, the only possibilities for $Hol^0(M)$ are
$\{ 1\}, SU(2), SU(3)$ or $G2$.
If $\Hol^0 = \{ 1\}$, then $k = 7$;
for the $SU(2)$ case, $k = 3$ and for
$SU(3)$, $k = 1$. Therefore, the only possibility is that
$\Hol^0(M) = G_2$.

}

Now, Proposition $\ref{p.wedgesplit}$ implies

\begin{prop}
Let $M$ be a compact $G_2$-manifold. Then,
\begin{eqnarray*}
H^2(M, \R) & = & H^2_7(M, \R) \otimes H^2_{14}(M, \R) \\
H^3(M, \R) & = & H^3_1(M, \R) \otimes H^3_{7}(M, \R) \otimes H^3_{27}(M, \R) \\
H^4(M, \R) & = & H^4_1(M, \R) \otimes H^4_{7}(M, \R) \otimes H^4_{27}(M, \R) \\
H^5(M, \R) & = & H^5_7(M, \R) \otimes H^5_{14}(M, \R)
\end{eqnarray*}

The Hodge star operator provides the isomorphism $H^k_\ell(M, \R) \cong H^{7-k}_{\ell}(M, \R)$,
which implies in particular that the refined Betti numbers satisfy $b^k_\ell = b^{7-k}_\ell$.

Also, $H^3_1(M, \R)$ is spanned by $\vp$, and $H^4_1(M, \R)$ is spanned by $*\vp.$
\end{prop}

\section{$Spin(7)$}
\label{s.Spin7}
The discussion on $Spin(7)$ will be analogous to the previous section. We start by giving the
definition of the group $Spin(7)$ as a subgroup of $GL(8, \R)$, as described in \cite{Bryant}.

\begin{defn}
Let $(x_1, \cdots, x_8)$ be coordinates for $\R^8$. Define a $4$-form $\Om_0$
by
\begin{eqnarray*}
\Om_0 &= & dx_{1234} +  dx_{1256} + dx_{1278} + dx_{1357}
- dx_{1368} - dx_{1458} - dx_{1467}  \\  & & -dx_{2358} - dx_{2367} -  dx_{2457}
+ dx_{2468} + dx_{3456}  + dx_{3478} + dx_{5678} ,
\end{eqnarray*}
where $dx_{ijkl}$ denotes the $4$-form $dx_i \wedge dx_j \wedge dx_k \wedge dx_l $.

The subgroup of $GL(8, \R)$ which preserves the  form $\Om_0$ is the
21-dimensional exceptional Lie Group $Spin(7)$.
It is compact, connected, semi-simple, and it fixes the Euclidean metric $g_0 = dx_1^2 + \cdots + dx_8^2$ and the orientation.

\end{defn}

We note that  the $4$-form $* \Om_0$, the Hodge dual of the $4$-form $\Om_0$ with respect to the metric $g_0$, is equal to $\Om_0$ itself, that is, $\Om_0$ is self-dual.

Let $M$ be an oriented $8$-manifold, and let $p \in M$.
In analogy to the $G_2$ case, we would like to define a $Spin(7)$ manifold by looking at the bundle of $4$-forms $\Om$ for which there exists an isomorphism $f$ from $T_pM$ to
$\R^8$ such that $f^*(\Om_0) = \Om$: define
$\cA_pM $, the {\em admissible set}, to be the set of all such forms where the isomorphism $f$ is oriented. We say that $\Om \in \wedge^4 T^*M$ is {\em admissible} if $\Om|_p \in \cA_pM$.

Let $\cA M $ be the bundle on $M$ with fiber $\cA_pM$. We note that, unlike the $G_2$
case, $\cA_pM$ is far from being open in $\wedge^4 T^*_pM$: being isomorphic to
$GL(8, \R)_+ / Spin(7)$, $\cA_pM$ has dimension $64 - 21 = 43$, while
$\wedge^4 T^*_pM$ has dimension $ {8 \choose 4} = 70$.

Nevertheless, there is still a one-to-one correspondence between oriented $Spin(7)$ structures on
$M$ and admissible $4$-forms $\Om \in \cA M$.
Therefore, we will abuse notation once again, and refer to the pair $(\Om, g)$
 as a $Spin(7)$-structure on $M$. Also, we say that a $Spin(7)$-structure $(\Om, g)$
 is {\em torsion-free} if $\nabla_g \Om = 0$.

\bigskip

The following results are analogous to Propositions \ref{p.TFAEG2}, \ref{p.wedgesplit}
and \ref{p.holcharactG2}, and their proofs can be done as an exercise by the reader,
using the similarities with the $G_2$ case.

\begin{prop}
\label{p.TFAESpin7}
Let  $M$ be an oriented riemannian $8$-manifold, and let $( \Om, g)$ be a
$Spin(7)$-structure on $M$. Then,
the following statements are equivalent:
\begin{enumerate}
\item
 $( \Om, g)$ is torsion-free.
 \item
$ \hol(g) \subset Spin(7)$, and $\Om$ is the induced $4$-form by the inclusion.
\item
$d\Om = 0$  on $M$.

\end{enumerate}
\end{prop}

It is true that the condition $d\Om = 0$ is linear, but the requirement that
$\Om$ is an admissible form is non-linear: so, we are dealing again with a non-linear
equation for $\Om$, as
in the $G_2$ case.

\begin{prop}
\label{p.wedgesplit2}
For $M$ a manifold with $Spin(7)$-structure $(\Om, g)$, each bundle $\wedge^k T^*M$  of
$k$-forms splits orthogonally as follows:
\begin{itemize}
\item
$
\wedge^1 T^*M = \wedge^1_8 $
\item
$\wedge^2 T^*M = \wedge^2_7  \oplus \wedge^2_{21}$
\item
$\wedge^3 T^*M = \wedge^3_8  \oplus \wedge^3_{48}$
\item
$\wedge^4 T^*M = \wedge_+^4 T^*M \oplus \wedge_-^4 T^*M$
\item
$\wedge^4_+ T^*M = \wedge^4_1  \oplus \wedge^4_7  \oplus \wedge^4_{27}$
\item
$\wedge^4_- T^*M = \wedge^4_{35}$,
\end{itemize}
where $\wedge^k_\ell  $ is an irreducible representation of $Spin(7)$ of dimension $\ell$.

The Hodge star operator $*$ provides an isometry between
$\wedge^k_\ell  $ and $\wedge^{8-k}_\ell  $, and we have
 $\wedge^4_1 = \langle \Om \rangle  $ .
 Also, the splitting on item $4$ on $\wedge^4_\pm T^*M$
 corresponds to the $\pm 1$-eigenspaces of the Hodge $*$ operator.

\end{prop}

\begin{prop}
\label{p.Spin7isRicciFlat}
If  $(M, g)$ is a riemannian 8-manifold whose holonomy group
$Hol(g) \subset Spin(7)$, then $g$ is Ricci-flat.

\end{prop}

\begin{thm}
If $(\Om, g)$ is a torsion-free $Spin(7)$ structure on $M$,
then $Hol^0(g)$ can only be one of the following subgroups of
$Spin(7)$: $\{ 1\}$, $SU(2)$,  $SU(2) \times SU(2)$, $SU(3)$,  $G_2$.,
$Sp(2)$, $SU(4)$ or $Spin(7)$.

\end{thm}

The theorem above allows us to construct $Spin(7)$ manifolds starting from a Calabi-Yau,
or a $G_2$ manifold. We now describe how to do this explicitly.

\begin{thm}
1)
If $(M, J, h, \theta )$ is a Calabi-Yau $2$-fold, with associated \kahler form
$\om$, and $x_1, \cdots, x_4$ are coordinates in $\R^4$, then we can define a
torsion-free $Spin(7)$ structure $(\Om, g)$ on $M$ by
\beeqnar
g & = & dx_1^2 + \cdots + dx^2_4 + h \\
\Om &=& dx_{1234}+ (dx_{12} + dx_{34})\wedge \om \\ & & 
(dx_{13} - dx_{24})\wedge \Re(\theta) - (dx_{14} + dx_{23})\wedge \Im(\theta)
+ \frac{1}{2}\om \wedge \om
\eeeqnar
2)
If $(M, J, h, \theta )$ is a Calabi-Yau $3$-fold, with associated \kahler form
$\om$, and $x_1, x_2$ are coordinates in $\R^2$. Now, we can define a
torsion-free $Spin(7)$ structure $(\Om, g)$  on $M \times \R^2$ by
\beeqnar
g & = & dx_1^2 + dx^2_2 + h \\
\Om &=&  dx_{12}\wedge \om +
dx_{1}\wedge \Re(\theta) - dx_{2} \wedge \Im(\theta)
+ \frac{1}{2}\om \wedge \om
\eeeqnar
3)
Now, if the Calabi-Yau manifold $(M, J, h, \theta )$ has dimension $4$,
the torsion-free $Spin(7)$-structure is given by $g = h$, and
\beeq
\Om = \frac{1}{2}\om\wedge\om + \Re(\theta)
\eeeq
4)
Finally, for the case where $(M, \vp, h)$ is a $G_2$ manifold (in which
case $\vp$ is torsion-free), we define a torsion-free $Spin(7)$ structure on
$M \times \R$ by $g =h + dx^2$
and $\Om = dx\wedge \vp + *\vp$, where $x$ is a coordinate for $\R$.

\end{thm}

%% file: Main.bbl
\providecommand{\bysame}{\leavevmode\hbox to3em{\hrulefill}\thinspace}
\providecommand{\MR}{\relax\ifhmode\unskip\space\fi MR }
\providecommand{\MRhref}[2]{%
  \href{http://www.ams.org/mathscinet-getitem?mr=#1}{#2}
}
\providecommand{\href}[2]{#2}
\begin{thebibliography}{GKM68}

\bibitem[Ale68]{Alexseevski}
D.~V. Alekseevski{\u\i}, \emph{Riemannian spaces with unusual holonomy groups},
  Funkcional. Anal. i Prilo\v zen \textbf{2} (1968), no.~2, 1--10. \MR{0231313
  (37 \#6868)}

\bibitem[Ast78]{Asterisque}
\emph{Premi\`ere classe de {C}hern et courbure de {R}icci: preuve de la
  conjecture de {C}alabi}, Ast\'erisque, vol.~58, Soci\'et\'e Math\'ematique de
  France, Paris, 1978, Lectures presented at a S{\'e}minaire held at the Centre
  de Math{\'e}matiques de l'{\'E}cole Polytechnique, Palaiseau, March--June
  1978. \MR{MR525893 (80j:53064)}

\bibitem[Ber60]{Berger}
Marcel Berger, \emph{Sur quelques vari\'et\'es riemanniennes suffisamment
  pinc\'ees}, Bull. Soc. Math. France \textbf{88} (1960), 57--71. \MR{MR0133781
  (24 \#A3606)}

\bibitem[Bes87]{Besse}
Arthur~L. Besse, \emph{Einstein manifolds}, Ergebnisse der Mathematik und ihrer
  Grenzgebiete (3) [Results in Mathematics and Related Areas (3)], vol.~10,
  Springer-Verlag, Berlin, 1987. \MR{867684 (88f:53087)}

\bibitem[BG72]{Brown-Gray}
Robert~B. Brown and Alfred Gray, \emph{Riemannian manifolds with holonomy group
  {${\rm S}pin$} (9)}, Differential geometry (in honor of {K}entaro {Y}ano),
  Kinokuniya, Tokyo, 1972, pp.~41--59. \MR{0328817 (48 \#7159)}

\bibitem[BGM71]{BGM}
Marcel Berger, Paul Gauduchon, and Edmond Mazet, \emph{Le spectre d'une
  vari\'et\'e riemannienne}, Lecture Notes in Mathematics, Vol. 194,
  Springer-Verlag, Berlin, 1971. \MR{0282313 (43 \#8025)}

\bibitem[Bry87]{Bryant}
Robert~L. Bryant, \emph{Metrics with exceptional holonomy}, Ann. of Math. (2)
  \textbf{126} (1987), no.~3, 525--576. \MR{916718 (89b:53084)}

\bibitem[Car26]{Cartan}
E.~Cartan, \emph{Sur une classe remarquable d'espaces de {R}iemann}, Bull. Soc.
  Math. France \textbf{54} (1926), 214--264. \MR{1504900}

\bibitem[CG72]{CheegerGromoll}
Jeff Cheeger and Detlef Gromoll, \emph{The splitting theorem for manifolds of
  nonnegative {R}icci curvature}, J. Differential Geometry \textbf{6}
  (1971/72), 119--128. \MR{0303460 (46 \#2597)}

\bibitem[Che99]{Chev}
Claude Chevalley, \emph{Theory of {L}ie groups. {I}}, Princeton Mathematical
  Series, vol.~8, Princeton University Press, Princeton, NJ, 1999, Fifteenth
  printing, Princeton Landmarks in Mathematics. \MR{1736269 (2000m:22001)}

\bibitem[DK00]{DuiKo}
J.~J. Duistermaat and J.~A.~C. Kolk, \emph{Lie groups}, Universitext,
  Springer-Verlag, Berlin, 2000. \MR{1738431 (2001j:22008)}

\bibitem[dR52]{deRham}
Georges de~Rham, \emph{Sur la reductibilit\'e d'un espace de {R}iemann},
  Comment. Math. Helv. \textbf{26} (1952), 328--344. \MR{0052177 (14,584a)}

\bibitem[GH81]{GH}
Marvin~J. Greenberg and John~R. Harper, \emph{Algebraic topology}, Mathematics
  Lecture Note Series, vol.~58, Benjamin/Cummings Publishing Co. Inc. Advanced
  Book Program, Reading, Mass., 1981, A first course. \MR{643101 (83b:55001)}

\bibitem[GKM68]{GKM}
D.~Gromoll, W.~Kingenberg, and W~Meyer, \emph{Riemannsche geometrie im
  grossen}, Lecture Notes in Math. Vol. 55, Springer, Berlin, 1968.

\bibitem[GT01]{GT}
David Gilbarg and Neil~S. Trudinger, \emph{Elliptic partial differential
  equations of second order}, Classics in Mathematics, Springer-Verlag, Berlin,
  2001, Reprint of the 1998 edition. \MR{MR1814364 (2001k:35004)}

\bibitem[Hel01]{Helg}
Sigurdur Helgason, \emph{Differential geometry, {L}ie groups, and symmetric
  spaces}, Graduate Studies in Mathematics, vol.~34, American Mathematical
  Society, Providence, RI, 2001, Corrected reprint of the 1978 original.
  \MR{1834454 (2002b:53081)}

\bibitem[HO56]{HO}
Jun-ichi Hano and Hideki Ozeki, \emph{On the holonomy groups of linear
  connections}, Nagoya Math. J. \textbf{10} (1956), 97--100. \MR{0082164
  (18,507e)}

\bibitem[Joy07]{J}
Dominic~D. Joyce, \emph{Riemannian holonomy groups and calibrated geometry},
  Oxford Graduate Texts in Mathematics, vol.~12, Oxford University Press,
  Oxford, 2007. \MR{2292510 (2008a:53050)}

\bibitem[KN96]{KN}
Shoshichi Kobayashi and Katsumi Nomizu, \emph{Foundations of differential
  geometry. {V}ol. {I}}, Wiley Classics Library, John Wiley \& Sons Inc., New
  York, 1996, Reprint of the 1963 original, A Wiley-Interscience Publication.
  \MR{1393940 (97c:53001a)}

\bibitem[Kod64]{Kodaira}
K.~Kodaira, \emph{On the structure of compact complex analytic surfaces. {I}},
  Amer. J. Math. \textbf{86} (1964), 751--798. \MR{0187255 (32 \#4708)}

\bibitem[Law85]{Laws}
H.~Blaine Lawson, Jr., \emph{The theory of gauge fields in four dimensions},
  CBMS Regional Conference Series in Mathematics, vol.~58, Published for the
  Conference Board of the Mathematical Sciences, Washington, DC, 1985.
  \MR{799712 (87d:58044)}

\bibitem[Lic62]{Lich}
Andr{\'e} Lichnerowicz, \emph{Th\'eorie globale des connexions et des groupes
  d'holonomie}, Consiglio Nazionale delle Ricerche Monografie Matematiche, Vol.
  2, Edizioni Cremonese, Rome, 1962. \MR{0165458 (29 \#2740)}

\bibitem[MS06]{BasicsofMTheory}
Andr{\'e} Miemiec and Igor Schnakenburg, \emph{Basics of {M}-theory}, Fortschr.
  Phys. \textbf{54} (2006), no.~1, 5--72. \MR{2192442 (2006k:81312)}

\bibitem[NN57]{NN}
A.~Newlander and L.~Nirenberg, \emph{Complex analytic coordinates in almost
  complex manifolds}, Ann. of Math. (2) \textbf{65} (1957), 391--404.
  \MR{0088770 (19,577a)}

\bibitem[Sal89]{Salamon}
Simon Salamon, \emph{Riemannian geometry and holonomy groups}, Pitman Research
  Notes in Mathematics Series, vol. 201, Longman Scientific \& Technical,
  Harlow, 1989. \MR{1004008 (90g:53058)}

\bibitem[San09]{Santoro}
Bianca Santoro, \emph{Introduction to evolution equations in geometry},
  Publica\c c\~oes Matem\'aticas do IMPA. [IMPA Mathematical Publications],
  Instituto Nacional de Matem\'atica Pura e Aplicada (IMPA), Rio de Janeiro,
  2009, 27${^{{}}{\rm{o}}}$ Col{\'o}quio Brasileiro de Matem{\'a}tica. [27th
  Brazilian Mathematics Colloquium]. \MR{2537139 (2010k:53102)}

\bibitem[Siu83]{Siu}
Y.~T. Siu, \emph{Every {$K3$} surface is {K}\"ahler}, Invent. Math. \textbf{73}
  (1983), no.~1, 139--150. \MR{707352 (84j:32036)}

\bibitem[Tia00]{Tianbook}
Gang Tian, \emph{Canonical metrics in {K}\"ahler geometry}, Lectures in
  Mathematics ETH Z\"urich, Birkh\"auser Verlag, Basel, 2000, Notes taken by
  Meike Akveld. \MR{MR1787650 (2001j:32024)}

\bibitem[Wel80]{Wells}
R.~O. Wells, Jr., \emph{Differential analysis on complex manifolds}, second
  ed., Graduate Texts in Mathematics, vol.~65, Springer-Verlag, New York, 1980.
  \MR{MR608414 (83f:58001)}

\bibitem[Wol65]{Wolf}
Joseph~A. Wolf, \emph{Complex homogeneous contact manifolds and quaternionic
  symmetric spaces}, J. Math. Mech. \textbf{14} (1965), 1033--1047. \MR{0185554
  (32 \#3020)}

\bibitem[Yau78]{Yau}
Shing~Tung Yau, \emph{On the {R}icci curvature of a compact {K}\"ahler manifold
  and the complex {M}onge-{A}mp\`ere equation. {I}}, Comm. Pure Appl. Math.
  \textbf{31} (1978), no.~3, 339--411. \MR{MR480350 (81d:53045)}

\end{thebibliography}
